\UseRawInputEncoding
\documentclass[onefignum,onetabnum]{siamart251216}

\usepackage[T1]{fontenc}
\usepackage{lmodern}
\usepackage{amsfonts,amsmath,amssymb}

\usepackage{lineno,nicefrac}
\usepackage{microtype}
\usepackage{mathrsfs}
\usepackage{mathtools}
\usepackage{cases}
\usepackage{booktabs}
\usepackage{tikz}
\usetikzlibrary{automata,arrows,positioning,calc}
\usepackage[mathscr]{eucal}
\usepackage{hyperref}
\usepackage[nameinlink]{cleveref}
\usepackage{aliascnt}

\usepackage{enumitem}
\usepackage{bm}
\usepackage{placeins}
\DeclareGraphicsExtensions{.png,.jpeg,.jpg,.pdf}
\graphicspath{{./Figures/}}

\newsiamremark{example}{Example}
\newsiamremark{remark}{Remark}
\newsiamthm{assumption}{Assumption}

\newsiamthm{case}{Case}
\AddToHook{env/assumption/begin}{\crefalias{theorem}{assumption}}
\AddToHook{env/remark/begin}{\crefalias{theorem}{remark}}
\AddToHook{env/example/begin}{\crefalias{theorem}{example}}
\AddToHook{env/case/begin}{\crefalias{theorem}{case}}

\crefname{equation}{eq.}{eqs.}
\Crefname{equation}{Eq.}{Eqs.}

\crefname{theorem}{theorem}{theorems}
\Crefname{theorem}{Thm.}{Thm.}

\crefname{lemma}{lemma}{lemmas}
\Crefname{lemma}{Lem.}{Lems.}

\crefname{corollary}{corollary}{corollaries}
\Crefname{corollary}{Cor.}{Cor.}

\crefname{proposition}{proposition}{propositions}
\Crefname{proposition}{Prop.}{Prop.}

\crefname{assumption}{asm.}{assums.}
\Crefname{assumption}{Asm.}{Asms.}

\crefname{remark}{remark}{remarks}
\Crefname{remark}{Remark}{Remark}

\crefname{definition}{definition}{definitions}
\Crefname{definition}{Def.}{Defs.}

\crefname{example}{example}{examples}
\Crefname{example}{Example}{Examples}

\crefname{appendix}{appendix}{appendices}
\Crefname{appendix}{Appendix}{Appendices}

\newcommand{\lb}{\left\{}
\newcommand{\rb}{\right\}}

\newcommand{\timehor}{\bm T_h}

\newcommand{\Def}{\overset{\textbf{def}}{=}}
\newcommand{\RR}{\mathbb{R}}

\newcommand{\tup}[1]{\textup{(}#1\textup{)}}
\usepackage{amsopn}

\newcommand{\sa}{\mathsf a}

\newcommand{\nnorm}[1]{{\left\vert\kern-0.25ex\left\vert\kern-0.25ex\left\vert #1 \right\vert\kern-0.25ex\right\vert\kern-0.25ex\right\vert}}

\newcommand{\cl}{\textbf{cl}}

\newcommand{\mW}{\mathcal W}

\newcommand{\bvphi}{\bm{\varphi}}

\DeclarePairedDelimiter{\abs}{\lvert}{\rvert}
\DeclarePairedDelimiter{\norm}{\lVert}{\rVert}
\DeclarePairedDelimiter{\set}{\{}{\}}

\newcommand{\mL}{\mathsf L}

\usepackage{textcomp}
\usepackage{xcolor}
\usepackage{graphics}

\usepackage{color}
\definecolor{darkgreen}{rgb}{0,0.5,0}

\usepackage{graphicx}
\usepackage{multicol}

\newcommand\anno{0} 
\newcommand\target{5} 
\pgfmathsetmacro\myend{\target-1-\anno}
\pgfmathsetmacro\tix{1}
\pgfmathsetmacro\myspacing{5.7/(\target-1-\anno)}

\headers{Collision-time regularity and residual stability}{H. Nick Zinat Matin, Y. Yeo, and M. L. Delle Monache}

\title{A discontinuous time-to-collision feedback: well-posedness, regularity of the collision-time map, and residual stability of its Hamilton–Jacobi approximations}

\author{Hossein Nick Zinat Matin{\thanks{\'{E}cole Polytechnique, Inria-Saclay, Alan Turing Building, Palaiseau, 91120, France (\email{hossein.matin@polytechnique.edu}).}} {\thanks{University of California, Berkeley, CA 94720, USA (\email{h-matin@berkeley.edu}, \email{yuneily@berkeley.edu}, \email{mldellemonache@berkeley.edu}).}}
{\thanks{These authors contributed equally to this work.}} 
{\thanks{This work benefited from Hi! PARIS and State funding managed by the French National Research Agency (ANR) under the France 2030 program, reference ANR-23-IACL-0005.}}
{\thanks{Corresponding author.}}
\and Yuneil Yeo{\footnotemark[2] \footnotemark[3]} 
\and Maria Laura Delle Monache\footnotemark[2]}

\newcommand{\Vbar}{\overline{\mathcal V}}
\newcommand{\abar}{\bar a}
\newcommand{\Om}{\Omega}
\newcommand{\Reg}{\mathcal R}
\newcommand{\Graz}{\Gamma}
\newcommand{\Phys}{\mathcal P}
\newcommand{\Slab}{\mathcal P^{+}}
\newcommand{\App}{\mathcal A}
\newcommand{\Tt}{\tilde{\mathcal T}}
\newcommand{\Fil}{\mathfrak F}

\newcommand{\one}{\mathbf 1}
\newcommand{\co}{\overline{\operatorname{co}}}

\begin{document}

\maketitle
\begin{abstract}
We analyze a finite-horizon time-to-collision map generated by an optimal-velocity prediction with bounded measurable leader acceleration and a feedback law that uses this map. The prediction describes nominal motion without the additional collision-avoidance feedback. Its first-contact geometry can produce discontinuities of the capped collision time at grazing, so smooth-feedback theory does not apply globally. For a prescribed leader--follower pair, we prove existence of Filippov solutions and exclusion of contact with negative relative velocity for every solution; a zero-speed touch is not excluded. Under explicit nondegeneracy hypotheses, trajectories cross the grazing set and uniqueness persists through the crossing. On the domain of transversal first contacts, the collision-time map is locally Lipschitz, continuously differentiable in the state, and differentiable in time outside a common null set. We quantify its sensitivity to the transversality margin and to jumps of the forcing. Characteristic and coarea arguments give boundary-plus-interior residual estimates for Hamilton--Jacobi approximations. An endpoint/age change of variables yields a stability coefficient controlled by flight time and terminal normal velocity, without inverse-transversality blowup. The estimates imply convergence conditional on vanishing continuum loss and pointwise identification under the stated trace and regularity assumptions. Numerical experiments on thirteen prescribed forcing profiles investigate boundary-error transport, sampling effects, and regularity predictions; a separate first-contact experiment numerically demonstrates post-threshold sensitivity sign changes. The guarantees concern the exact feedback and value-function approximation. Learned-feedback safety and convergence of a training algorithm remain open.

\end{abstract}

\begin{keywords}
Time-to-collision; car-following dynamics; hitting-time regularity; Hamilton--Jacobi equations; residual stability; physics-informed neural networks; Filippov solutions
\end{keywords}

\begin{MSCcodes}
34A36, 49L12, 93C15, 65M15, 68T07
\end{MSCcodes}


\section{Introduction and related work}\label{sec:introduction}

Microscopic traffic models describe vehicle interactions through ordinary differential equations. Classical formulations \cite{pipes1953operational,newell1961nonlinear}, including the optimal velocity model (OVM) \cite{bando1995dynamical}, reproduce stop-and-go waves and spontaneous congestion \cite{dong2009velocity,jin2020dynamical}. Their stability properties have been studied extensively \cite{ward2011criteria,wilson2011car,matin2020nonlinear,matin2020nonlinear2}, and their extensions support driver-assistance and autonomous-vehicle control \cite{treiber2013traffic,9922604,helbing1998generalized,jiang2001full,10384086,matin2024analytical,yeo2025energy,yeo2026optimal}. A response based on spacing alone does not distinguish a closing follower from one whose leader is pulling away. Relative velocity and time-to-collision (TTC) provide additional information about the evolving interaction \cite{winsum1996choice,lee1976theory,nikolaou2023exploiting,wang2021review,wang2025dynamic}.

\paragraph{Nominal prediction and collision-avoidance feedback}
The classical TTC of Hayward \cite{hayward1971near} extrapolates the present gap under constant relative velocity. Potential and modified TTC metrics, acceleration-based refinements, and two-dimensional extensions incorporate further kinematic information \cite{wakabayashi2003traffic,ozbay2008derivation,yang2010application,ward2015extending,matin2025second,brown2005adjusted,essa2015simulated,caliendo2012microsimulation,hou2015new,ambros2014road}. Such quantities are counterfactual forecasts: the collision predicted under a specified continuation of motion can be avoided by a subsequent intervention. There is no inconsistency in using a forecast to trigger that intervention. Dynamical and model-predictive collision-time metrics also have precedents; for example, MPrISM accounts for agent dynamics and actuation in an interactive safety metric \cite{weng2020model}. Several microscopic feedback models use TTC-related terms \cite{jiao2021dynamic,mammar2005modified,7210342,abdelhalim2022real,jin2012staggered,li2023analysis}.

Here the nominal continuation is an OVM prediction driven by a prescribed leader acceleration. Let $\Xi$ and $\Lambda$ denote predicted gap and relative velocity, and let $\bm z=(\bar t,\bar X,\bar V)$ be an initial time and relative state. With absolute-time flow notation, define
\begin{equation}\label{E:intro-tau}
\begin{aligned}
\tau(\bm z)\Def\inf\bigl\{\theta\in(0,\timehor-\bar t]:\;&
\Xi(\bar t+\theta;\bm z)=0,\\[-1mm]
&\Lambda(\bar t+\theta;\bm z)<0\bigr\},
\qquad \inf\varnothing\Def\timehor-\bar t .
\end{aligned}
\end{equation}
The precise definition is \cref{D:TTC}. Tangential touches are not counted in \eqref{E:intro-tau}; on the transversal domain studied in \cref{sec:ttc-regularity}, the selected contact is the first zero of the predicted gap. The prediction omits the additional TTC feedback. In relative variables, its field $\mathsf G$ and the closed-loop field $\mathsf F$ satisfy
\[
\tau=\operatorname{HitTime}(\operatorname{flow}_{\mathsf G}),
\qquad
\mathsf F(t,X,V)=\mathsf G(t,X,V)
   +\left(0,-\frac{\beta V}{\tau(t,X,V)^2}\right).
\]
Thus the construction has no implicit fixed-point dependence of $\tau$ on the flow of $\mathsf F$. The advantage of the OVM forecast is that it includes nominal acceleration dynamics \cite{bando1995dynamical,bando1998analysis}. The analytical difficulty is that its capped hitting time can jump when the prediction grazes the collision set, producing a discontinuous feedback coefficient.

Related singular-damping models already yield collision exclusion in other settings. The deterministic analysis in \cite{nick2022near} treats a spacing-singular car-following model with a constant-speed leader, while \cite{matin2024analytical} studies a nonlinear microscopic controller with prescribed leader dynamics. The present coefficient is an implicitly computed hitting time rather than an explicit power of the gap. We allow bounded measurable leader acceleration and prove impact prevention for all Filippov solutions. The conclusion is correspondingly precise: a contact with negative relative velocity is impossible, but a finite-time touch with zero relative velocity is not excluded.

\paragraph{Hitting-time regularity and grazing}
Exit-time and minimum-time regularity are classical subjects \cite{bardi1997optimal,cannarsa2004semiconcave}. At an incoming collision state, strict negativity of the normal velocity is the noncharacteristic condition for the deterministic prediction; it is also the singleton-control counterpart of an inward-pointing condition. Absence of a control optimization does not remove this connection with exit-time theory.
The questions addressed here concern quantitative estimates for this particular hitting-time map, its merely measurable time dependence, and the effect of its singularities on the feedback dynamics.

Square-root singularities at grazing and regularity of transverse event maps are established features of nonsmooth dynamics \cite{di2008bifurcations,burden2016event}. In the present model, the crossing-side sensitivity of the collision time can become unbounded near a graze. We identify conditions under which the feedback nevertheless moves every solution transversally through the grazing set. After crossing, the resulting comparison coefficient is singular but integrable, which permits uniqueness. These results are conditional on explicitly stated nondegeneracy hypotheses; we do not assert that those hypotheses are generic.

On the domain of transversal first contacts, we prove compact-set Lipschitz bounds, continuous state differentiability, and time differentiability outside a common null set. The gradient formula quantifies the role of the terminal transversality margin. Jumps of the leader acceleration produce jumps of the time derivative wherever the corresponding state sensitivity is nonzero. A variational comparison supplies a sufficient positive-sensitivity horizon; the numerical study distinguishes that comparison from sensitivity evaluated at an actual first contact.

\paragraph{Computation and Hamilton--Jacobi approximation}
Evaluating \eqref{E:intro-tau} generally requires an ODE solve and contact detection. A surrogate can reduce the cost of repeated evaluations. This is a computational motivation, not a theorem that online simulation violates safety or that a neural network has zero delay. The exact-feedback analysis assumes the prescribed leader profile is available on the prediction interval; forecast errors, sampled implementation, and computation delays require additional analysis.

TTC value functions and Hamilton--Jacobi formulations have appeared previously in reachability-based navigation \cite{li2020robust,li2021prediction}. Hamilton--Jacobi reachability and its learned approximations also provide general approaches to safety analysis \cite{mitchell2005time,evans2022partial,bansal2021deepreach}. Here the prediction field is fixed, and the value is its finite first-contact time. The Hamiltonian is affine in the gradient, so its characteristics are the prediction trajectories. With the incoming trace, an absolutely continuous characteristic solution is determined along those trajectories. This differs from treating a continuous almost-everywhere solution as automatically unique, which is insufficient even for a linear transport equation.

Neural approximations and operator-learning methods offer compact representations of dynamics and solution maps \cite{lu2019deeponet,kovachki2023neural,chen2018neural}; related control analyses are developed in \cite{ruiz2023neural,ruiz2024control,bhan2023neural,krstic2024neural}. Physics-informed learning encodes the governing equation through residual terms \cite{raissi2019}. Convergence and residual-error analyses include \cite{shin2020convergence,qian2023physics}. In particular, Liu, Huang, and Protopapas derive characteristic residual bounds for first-order linear PDEs with exact boundary data \cite{liu2023residual}, while Mishra and Molinaro relate PDE stability to approximation and sampling errors \cite{mishra2023estimates}. The infinite-horizon Lyapunov and Hamilton--Jacobi--Bellman analysis of \cite{liu2026verifiable} uses pointwise relative residual bounds to control integration over an unbounded time interval. Our finite-contact setting instead permits integrated boundary and interior losses, but requires control of the endpoint geometry and of the measure transported along characteristics.

The main approximation result is therefore a residual-stability theorem specialized to the collision-time geometry. A coarea argument controls the incoming trace error, and Liouville's formula controls the interior residual. A second proof uses first-contact endpoints and backward ages directly: the volume element is terminal normal speed times an exponential flight-time factor. This yields an alternative explicit coefficient without an inverse-transversality singularity, even though pointwise gradients can diverge near grazing.

The results concern continuum losses. They imply convergence of an admissible sequence if its continuum loss tends to zero. They do not show that a training algorithm produces such a sequence, that the chosen network class attains zero loss, or that a vanishing sampled loss guarantees a vanishing continuum loss. The regularity estimates identify useful approximation structure \cite{yarotsky2017}, but no network-size approximation rate is claimed here.
For discontinuous forcing, a fixed-signal $C^1$ candidate may fail to attain zero total loss when the value has a genuine time-derivative jump. That does not imply a positive loss infimum or an unavoidable positive bias. Moreover, a smooth network evaluated at the evolving leader speed generally induces a locally Lipschitz, rather than $C^1$, fixed-signal slice; \cref{R:rough-candidates} makes this distinction explicit.

\subsection*{Contributions}
The contributions are the following.
\begin{itemize}[leftmargin=0.15in]
\item \textbf{Impact prevention and conditional uniqueness.}
For the hitting-time feedback with bounded measurable forcing, we prove existence of Filippov solutions for physical leader--follower data and impact prevention for every solution. We establish local uniqueness in the regular set and uniqueness through nondegenerate grazing by transverse crossing and an integrable singular comparison estimate.
\item \textbf{Quantitative collision-time regularity.}
On compact subsets of the transversal domain, we give explicit Lipschitz bounds, the continuously varying state gradient, and a common exceptional set for time differentiability. We quantify forcing-jump sensitivity and provide a sufficient horizon for positivity of the relative-velocity sensitivity.
\item \textbf{Endpoint geometry and residual stability.}
We prove a positive lower bound for the endpoint map's two-dimensional Jacobian and control the boundary loss by coarea. The endpoint/age parameterization then gives the explicit estimate in \cref{T:flux-stability}, whose coefficient depends on flight time and terminal normal speed. The estimates cover locally Lipschitz candidates with a continuous incoming trace, including the finite-jump fixed-signal slices used in the experiments.
\item \textbf{Numerical diagnostics of approximation and regularity.}
The learning experiments use thirteen prescribed forcing profiles. The principal run includes an additional dynamic-boundary penalty; a smaller comparison reports the actual boundary-loss weights rather than isolating the extra target. A characteristic-tube experiment compares estimated continuum errors and losses with an explicit analytical coefficient. Reference calculations test gradient scaling and forcing-jump sensitivities, and a separate reproducible experiment evaluates post-threshold sensitivity at actual first contact.
\end{itemize}
The exact-controller conclusions concern one prescribed leader and one follower. The learning results do not establish safety of a controller driven by the surrogate, string stability, generalization to unseen forcing profiles, or convergence of the optimization algorithm.

\section{A nominal-prediction TTC feedback model} \label{sec: general_framework_dTTC}
We specify the nominal prediction and the feedback that uses its collision-time estimate. The multi-vehicle notation motivates the construction; the analytical results below concern one follower and one prescribed leader.

\subsection{The Model} \label{sec: motivation}

As discussed in \Cref{sec:introduction}, a control law depending only on inter-vehicle spacing carries no information about the rate of closure,
and incorporating the relative velocity through a time-to-collision term supplies precisely that information. We now make the resulting model precise. 

We propose the nonlinear feedback control:
\begin{equation}\label{E:main_dynamics}
\begin{cases}
\dot x_\mL  = v_\mL \\
\dot v_\mL = \sa_\mL \\
\dot x_n = v_n\\
\dot v_n = \alpha \lb \mathcal V(x_{n-1}(t)-x_n(t))-v_n(t) \rb + \beta\frac{(v_{n-1}(t)-v_n(t))}{\big(\tau_n(t, x_{n-1}, x_n, v_{n -1}, v_n)\big)^2}  \\
(x_\mL(0) , x_n(0) , v_\mL(0), v_n(0)) = (x_{\mL, \circ}, x_{n,\circ}, v_{\mL, \circ}, v_{n,\circ}), 
\end{cases}
\end{equation}
for $n \in \set{1, \cdots, N}$ or in compact form, we define $\bvphi_n \coloneqq (x_{n-1}, x_n, v_{n -1}, v_n)$ and the dynamical model \Cref{E:main_dynamics} can be written in the form of 
\begin{equation}\label{E:compact-form-main}
    \begin{cases}
        \dot{\Phi} =  \mathsf F(t, \Phi) = \mathsf F(t, \Phi, \tau(t, \Phi)), \quad \Phi = (\bvphi_1, \bvphi_2, \cdots, \bvphi_N), \quad \tau = (\tau_1, \cdots, \tau_N)\\
        \Phi(t_\circ) = \Phi_\circ.
    \end{cases}
\end{equation}
{The index $n \in \set{1, \cdots, N}$, indicates the vehicle, and the position and velocity of $n$-th vehicle are denoted by $x_n$ and $v_n$, respectively. The preceding vehicle of vehicle $n \in \set{1, \cdots, N}$ is indexed by $n-1$. The position, velocity, and acceleration of the leading vehicle are denoted by $x_\mL$, $v_\mL$, and $\sa_\mL$, respectively. For notational consistency, the exogenous leader, which is the preceding vehicle of vehicle 1, is identified as vehicle $0$, making $x_0 = x_\mL$ and $v_0 = v_\mL$.} Here, $\tau_n$ is the time-to-collision (TTC) control strategy which will be defined later in this section---see \Cref{E:TTC-main}. 
The proposed platoon may include all interacting AVs or AVs interacting with a human-driven vehicle (HV) leader. 

The optimal velocity function $\mathcal V:\mathbb R \to \mathbb R$ is a smooth, bounded, strictly increasing Lipschitz function with $\mathcal V(0)=0$ and $\mathcal V\ge0$ on $[0,\infty)$ (its values at negative gaps are used only by the prediction flow, which is continued through the collision set). In this paper, we consider 
\begin{equation}  \label{E:compact-two-dynamics-maintion}
\mathcal V(x_{n-1}-x_n) \Def \tanh(x_{n-1}-x_n-2)+\tanh(2),
\end{equation}
see \cite{nick2022near, bando1995dynamical} for properties of this function, and \cite{matin2024analytical} for model's sensitivity analysis with respect to the coefficients $\alpha$ and $\beta$.

In addition, we assume a reasonable behavior for the leading vehicle from a practical point of view; that is, the leader's acceleration is a prescribed measurable function of time,
\begin{equation}\label{E:accel_ell}
\sa_\mL\in L^\infty\bigl(0,\timehor\bigr),
\qquad\text{with}\qquad
0\le v_{\mL,\circ}+\int_0^t\sa_\mL(s)\,ds\le v_{\max}<\infty
\quad\text{for all }t\in[0,\timehor],
\end{equation}
where $v_{\max}$ is the maximum permitted speed. In particular the leader does not react to the follower. Treating $\sa_\mL$ as known is the standard setting of prediction-based feedback: it can be a forecast of the leader's behavior, estimated from observed trajectories (for instance through a parametric representation), or, in a platoon, the observed motion of the predecessor; the leader's acceleration can also be viewed as a fixed control input of the leader. 

Since the leader enters the analysis only through a prescribed bounded acceleration profile $\sa_\mL$ and a speed in $[0,v_{\max}]$, the analytical parts consider $N = 1$ (a leader--follower pair). Extending the analytical statements to a platoon in which each predecessor is itself a controlled follower is not immediate: a controlled predecessor is only proved to be absolutely continuous, and the prediction of pair $n$ would depend on the predecessor's future motion, hence on upstream states; we do not pursue this here (cf.\ \cite{matin2024analytical} for the uncontrolled case). Therefore, from hereon in this paper, we focus our analysis on a general leader-follower control dynamics in the compact form
\begin{equation}\label{E:compact-two-dynamics-main}
    \begin{cases}
        \dot{\bvphi} =  \mathsf F(t, \bvphi) , \quad t \in [0, \timehor]\\
        \bvphi(t_\circ) = \bvphi_{\circ}.
    \end{cases}
\end{equation}
for $\bvphi = (x_\mL, x, v_\mL, v)$. 


\subsection{The Dynamic-Based TTC and its Approximation} \label{sec:architecture} 
The TTC map $\tau$ (\cref{D:TTC}) is computed from the nominal OVM prediction. Starting from the current state, the prediction follows \eqref{E:OVM} with the prescribed future leader profile, and omits the additional $\beta$-feedback in \eqref{E:main_dynamics}. The resulting hitting time is a forecast of this nominal motion; it is not the actual contact time of the controlled trajectory. OVM supplies a simple acceleration model in which the follower relaxes toward a spacing-dependent desired velocity. Its smoothness, one-sided physical bounds, variational equation, and constant state divergence are used in the proofs. Other prediction models could be treated under corresponding structural hypotheses, but Lipschitz flow dependence alone is not sufficient for all results below.


For comparison, the nominal OVM and the spacing-singular OVFL model have the accelerations
\begin{align}
\text{OVM}:\quad \dot v(t)
 &=\alpha\{\mathcal V(x_\mL(t)-x(t))-v(t)\},\label{E:OVM}\\
\text{OVFL}:\quad \dot v(t)
 &=\alpha\{\mathcal V(x_\mL(t)-x(t))-v(t)\}
   +\beta\frac{v_\mL(t)-v(t)}{(x_\mL(t)-x(t))^2}.\label{E:OVFL}
\end{align}
The TTC feedback replaces the explicit gap factor in OVFL by a collision-time factor obtained from the nominal flow. The safety analysis in \cref{sec: analytical_results} concerns this continuous-time exact-TTC law. Numerical statements about a sampled implementation require a specified prediction horizon, update step, contact detector, and no-hit convention; they are not consequences of the continuous-time theorem.

\section{Well-posedness and collision avoidance of the closed loop}\label{sec: analytical_results}
In this section we study the two-vehicle closed loop \eqref{E:compact-two-dynamics-main}. The TTC $\tau$ that enters the feedback is the hitting time of the collision set under the OVM prediction flow started from the current state; it is defined in \cref{D:TTC} below (this is the definition anticipated in \eqref{E:intro-tau}). The difficulty is twofold. The feedback coefficient $1/\tau^2$ is defined through a flow that the feedback itself perturbs, and $\tau$ can be discontinuous at states whose prediction \emph{grazes} the collision set. We show that away from this grazing set the closed loop is a Carath\'eodory system with unique solutions; that across it solutions exist in the sense of Filippov and, at nondegenerate grazing states, cross transversally so that uniqueness survives; and that no solution ever reaches the collision set with positive closing speed.

\paragraph{Notation for this section.}  We work in the difference variables $X=x_\mL-x$ (gap) and $V=v_\mL-v$ (relative velocity), $\xi=(X,V)$, so that the follower speed is $v=v_\mL(t)-V$. The predicted gap and relative velocity under the OVM flow are written $(\Xi,\Lambda)$, as in \eqref{E:intro-tau}. Throughout, $z=(\bar t,\bar X,\bar V)$ denotes a state (a point of space-time), and $\abs{\cdot}$ the Euclidean norm.

\subsection{Setting}\label{S:setting}

\begin{assumption}[Standing assumptions]\label{A:standing}
\begin{enumerate}[label=\tup{\roman*}]
\item $\alpha>0$ and $\beta>0$.
\item $\sa_\mL\in L^\infty(0,\timehor)$, and the leader velocity
  $v_\mL(t)\Def v_{\mL,\circ}+\int_0^t\sa_\mL(r)\,dr$ satisfies $0\le v_\mL(t)\le{v_{\max}}$
  for all $t\in[0,\timehor]$.  The leader position is
  $x_\mL(t)\Def x_{\mL,\circ}+\int_0^t v_\mL(r)\,dr$.  We fix a representative of
  $\sa_\mL$ with $\abs{\sa_\mL(t)}\le\norm{\sa_\mL}_{L^\infty}$ for \emph{every} $t$.  This
  makes the pointwise-in-$t$ bounds below hold at every $t$, not just a.e.
\item ${v_{\max}}\ge\Vbar\Def 1+\tanh 2$.
\end{enumerate}
\end{assumption}

\begin{remark}
Item~(iii) says that the maximal permitted speed is at least the supremum of the optimal
velocity function; it is exactly what makes $[0,{v_{\max}}]$ invariant for the follower speed
(\cref{L:prediction-bounds,L:invariance}).
\end{remark}

\paragraph{Facts on $\mathcal V$.}
$\mathcal V(r)=\tanh(r-2)+\tanh 2$ is smooth and strictly increasing, with
$0<\mathcal V'\le 1$, $\mathcal V(0)=0$, and
\begin{equation}\label{E:V-range}
\tanh 2-1<\mathcal V(r)<\Vbar\quad(r\in\RR),\qquad
\mathcal V(r)>0\quad(r>0).
\end{equation}
In particular $\abs{\mathcal V}<\Vbar$ and $-\mathcal V<1-\tanh 2<\Vbar$ on $\RR$.

\paragraph{The deceleration constant.}  We set
\begin{equation}\label{E:abar}
\abar\Def\norm{\sa_\mL}_{L^\infty}+\alpha\bigl(\Vbar+{v_{\max}}+1\bigr).
\end{equation}
(The ``$+1$'' gives room to work on an open neighbourhood of the physical states; see
\eqref{E:slab}.)

\paragraph{Reduction to the difference variables.}
Since the leader trajectory $(x_\mL,v_\mL)$ is prescribed, the map
$(x,v)\mapsto(X,V)\Def(x_\mL(t)-x,\;v_\mL(t)-v)$ is, for each $t$, an affine bijection
of $\RR^2$, and system \eqref{E:compact-two-dynamics-main} is equivalent to
\begin{equation}\label{E:closed-loop}
\dot X=V,\qquad
\dot V=\mathsf f(t,X,V)\Def \mathsf g(t,X,V)-\beta\,\frac{V}{\tau(t,X,V)^2},
\end{equation}
where
\begin{equation}\label{E:g}
\mathsf g(t,X,V)\Def \sa_\mL(t)-\alpha\bigl[\mathcal V(X)+V-v_\mL(t)\bigr]
=\sa_\mL(t)-\alpha\bigl[\mathcal V(X)-v\bigr],\qquad v=v_\mL(t)-V .
\end{equation}
Indeed $\dot v=\alpha(\mathcal V(X)-v)+\beta(v_\mL-v)/\tau^2$ and $\dot V=\sa_\mL-\dot v$.
We write $\xi=(X,V)$ and $\mathsf F(t,\xi)\Def(V,\mathsf f(t,\xi))$.  The follower
speed is always $v=v_\mL(t)-V$.

\paragraph{Sets.}  We use
\begin{align}
\Om&\Def\{(t,X,V):\ t\in[0,\timehor),\ X>0,\ V\in\RR\},\label{E:Omega}\\
\Phys&\Def\{(t,X,V)\in\Om:\ 0\le v_\mL(t)-V\le{v_{\max}}\}\qquad\text{(physical states)},\label{E:phys}\\
\Slab&\Def\{(t,X,V)\in\Om:\ -1<v_\mL(t)-V<{v_{\max}}+1\}\qquad\text{(open slab, }\Phys\subset\Slab).\label{E:slab}
\end{align}
$\Om$ and $\Slab$ are relatively open in $[0,\timehor)\times\RR^2$ ($v_\mL$ is continuous).

\paragraph{Extension of the leader data.}  For the prediction flow only, we extend
$\sa_\mL$ by $0$ on $(\timehor,\timehor+1]$, so that $v_\mL(s)=v_\mL(\timehor)$ there.
This is a technical device used in \cref{L:R-open} to treat predicted contacts exactly at
$s=\timehor$; it does not affect the closed loop.

\subsection{The prediction flow}\label{S:prediction}

For a datum $z=(\bar t,\bar X,\bar V)\in[0,\timehor)\times\RR^2$, the predicted (OVM)
relative motion $(\Xi,\Lambda)(s;z)$, $s\in[0,\timehor+1]$, solves
\begin{equation}\label{E:prediction}
\partial_s\Xi=\Lambda,\qquad
\partial_s\Lambda=\mathsf g(s,\Xi,\Lambda)
=\sa_\mL(s)-\alpha\bigl[\mathcal V(\Xi)+\Lambda-v_\mL(s)\bigr],\qquad
(\Xi,\Lambda)(\bar t;z)=(\bar X,\bar V).
\end{equation}
It uses the same leader data as the closed loop and omits the TTC term.  We write
$\mathsf G(s,Z)\Def(\Lambda,\mathsf g(s,Z))$ for $Z=(\Xi,\Lambda)$, and
\[
\zeta(s;z)\Def v_\mL(s)-\Lambda(s;z)\qquad\text{(predicted follower speed)},\qquad
\bar\zeta\Def v_\mL(\bar t)-\bar V .
\]

\begin{lemma}[Flow]\label{L:flow}
\begin{enumerate}[label=\tup{\alph*}]
\item For every $z$, \eqref{E:prediction} has a unique absolutely continuous solution on
  $[0,\timehor+1]$.  Moreover $\Xi(\cdot;z)\in C^1$ with $\partial_s\Xi=\Lambda$.
\item For every bounded $B\subset[0,\timehor)\times\RR^2$ there are constants
  $m_B,L_B<\infty$ such that, for all $z,z_1,z_2\in B$ and $s,s_1,s_2\in[0,\timehor+1]$,
  \begin{align}
  \abs{(\Xi,\Lambda)(s_1;z)-(\Xi,\Lambda)(s_2;z)}&\le m_B\abs{s_1-s_2},\label{E:flow-lip-s}\\
  \abs{(\Xi,\Lambda)(s;z_1)-(\Xi,\Lambda)(s;z_2)}&\le L_B\abs{z_1-z_2}.\label{E:flow-lip-z}
  \end{align}
  In particular $(s,z)\mapsto(\Xi,\Lambda)(s;z)$ is continuous.
\end{enumerate}
\end{lemma}

\begin{proof}
(a) $\mathsf G$ is measurable in $s$, and for all $s$ and $Z_1,Z_2$,
\[
\begin{aligned}
\abs{\mathsf G(s,Z_1)-\mathsf G(s,Z_2)}
&\le\abs{\Lambda_1-\Lambda_2}
  +\alpha\abs{\mathcal V(\Xi_1)-\mathcal V(\Xi_2)}\\
&\quad+\alpha\abs{\Lambda_1-\Lambda_2}\\
&\le L_{\mathsf G}\abs{Z_1-Z_2},\qquad L_{\mathsf G}\Def1+2\alpha,
\end{aligned}
\]
using $\mathcal V'\le1$.  By \eqref{E:V-range} and $0\le v_\mL\le{v_{\max}}$,
$\abs{\mathsf g(s,Z)}\le\abar+\alpha\abs\Lambda$, so
$\abs{\mathsf G(s,Z)}\le(1+\alpha)\abs Z+\abar$.  Carath\'eodory's theorem with a globally
Lipschitz right-hand side gives a unique global solution, forward and backward from
$\bar t$ \cite[Ch.~I, \S5]{Hale}.  Since $\Lambda$ is continuous and
$\Xi(s)=\bar X+\int_{\bar t}^s\Lambda$, $\Xi$ is $C^1$.

(b) Let $\varrho_B\Def\sup_{z\in B}\abs{(\bar X,\bar V)}$.  Gr\"onwall's inequality
(forward and backward) gives
$\abs{Z(s;z)}\le R_B\Def(\varrho_B+\abar(\timehor+1))e^{(1+\alpha)(\timehor+1)}$, hence
\eqref{E:flow-lip-s} with $m_B\Def(1+\alpha)R_B+\abar$.

For \eqref{E:flow-lip-z}, let $z_i=(\bar t_i,\bar Z_i)$ and assume $\bar t_1\le\bar t_2$.
By \eqref{E:flow-lip-s},
$\abs{Z(\bar t_2;z_1)-\bar Z_1}\le m_B\abs{\bar t_2-\bar t_1}$.  The functions $Z(\cdot;z_1)$
and $Z(\cdot;z_2)$ solve the same ODE with values $Z(\bar t_2;z_1)$ and $\bar Z_2$ at time
$\bar t_2$, so Gr\"onwall on $[0,\timehor+1]$ gives
\[
\abs{Z(s;z_1)-Z(s;z_2)}\le e^{L_{\mathsf G}(\timehor+1)}\abs{Z(\bar t_2;z_1)-\bar Z_2}
\le e^{L_{\mathsf G}(\timehor+1)}\bigl(\abs{\bar Z_1-\bar Z_2}+m_B\abs{\bar t_1-\bar t_2}\bigr),
\]
which is \eqref{E:flow-lip-z} with $L_B\Def e^{L_{\mathsf G}(\timehor+1)}(1+m_B)$.
\end{proof}

\begin{remark}
The initial-time estimate is obtained by a ``startup gap'' plus initial-data stability.
It does \emph{not} use any time-shift invariance, which the system lacks because
$\sa_\mL$ and $v_\mL$ depend on time.
\end{remark}

\begin{lemma}[One-sided bounds on the prediction]\label{L:prediction-bounds}
Let $z=(\bar t,\bar X,\bar V)$ and $\hat v\Def\max\{\bar\zeta,{v_{\max}}\}$.  For
$s\in[\bar t,\timehor+1]$:
\begin{enumerate}[label=\tup{\alph*}]
\item $\zeta(s;z)\le\hat v$; consequently $\Lambda(s;z)\ge-\hat v$.
\item If $\bar\zeta<{v_{\max}}+1$, then $\partial_s\Lambda(s;z)\le\abar$ for a.e. $s$.
\end{enumerate}
\end{lemma}

\begin{proof}
From \eqref{E:prediction} and $\dot v_\mL=\sa_\mL$ a.e.,
$\partial_s\zeta=\sa_\mL-\partial_s\Lambda=\alpha(\mathcal V(\Xi)-\zeta)$ a.e.  The
right-hand side is continuous in $s$, so $\zeta\in C^1$.

(a) Let $u\Def(\zeta-\hat v)^+$.  Then $u$ is absolutely continuous with
$u'=\one_{\{\zeta>\hat v\}}\,\partial_s\zeta$ a.e.  On $\{\zeta>\hat v\}$ we have
$\mathcal V(\Xi)<\Vbar\le{v_{\max}}\le\hat v<\zeta$, so $\partial_s\zeta<0$.  Hence $u'\le0$
a.e., and $u(\bar t)=0$ gives $u\equiv0$.  Then $\Lambda=v_\mL-\zeta\ge-\zeta\ge-\hat v$,
using $v_\mL\ge0$.

(b) $\partial_s\Lambda=\sa_\mL-\alpha\mathcal V(\Xi)+\alpha\zeta
\le\norm{\sa_\mL}_\infty+\alpha(1-\tanh2)+\alpha\hat v$ by \eqref{E:V-range} and (a).
If $\bar\zeta<{v_{\max}}+1$ then $\hat v<{v_{\max}}+1$, and $1-\tanh 2<\Vbar$ gives the bound $\abar$.
\end{proof}

\subsection{The TTC map}\label{S:TTC}

\begin{definition}[TTC, first predicted contact, regular and grazing sets]\label{D:TTC}\label{E:TTC-main}
For $z=(\bar t,\bar X,\bar V)\in\Om$:
\begin{align}
\Tt(z)&\Def\inf\{s\in[\bar t,\timehor]:\ \Xi(s;z)=0,\ \Lambda(s;z)<0\},\quad\inf\varnothing\Def\timehor,
&\tau(z)&\Def\Tt(z)-\bar t,\label{E:TTC}\\
\sigma(z)&\Def\inf\{s\in[\bar t,\timehor]:\ \Xi(s;z)=0\},\quad\inf\varnothing\Def+\infty .\label{E:sigma}
\end{align}
The \emph{regular set} is
$\Reg\Def\{z\in\Om:\ \sigma(z)=+\infty\ \text{ or }\ \Lambda(\sigma(z);z)<0\}$, and the
\emph{grazing set} is $\Graz\Def\Om\setminus\Reg$.
\end{definition}

\begin{remark}\label{R:Gamma}
Since $\Xi(\bar t;z)=\bar X>0$ and $\Xi>0$ on $[\bar t,\sigma(z))$, we have
$\Lambda(\sigma(z);z)=\partial_s\Xi(\sigma(z);z)\le0$ whenever $\sigma(z)<\infty$.  Hence
\[
\Graz=\{z\in\Om:\ \sigma(z)<\infty,\ \Lambda(\sigma(z);z)=0\}:
\]
these are the states whose \emph{first} predicted contact has zero relative velocity.
This includes tangential touches, where the gap reopens.  It also includes zero-speed
\emph{crossings}, where the predicted gap becomes negative with $\Lambda=0$ at the crossing
point.  The latter are not counted by $\Tt$, which only registers contacts with
$\Lambda<0$.
\end{remark}

\begin{remark}[First contact versus first transversal contact]\label{R:D-not-in-R}
The set $\{\Tt<\timehor,\ \Lambda(\Tt)<0\}$ of states whose first \emph{transversal}
predicted contact occurs before the horizon is \emph{not} contained in $\Reg$: a prediction
may touch $0$ tangentially at $\sigma(z)$, reopen, and cross transversally later.  Such
states lie in $\Graz$, where $\tau$ is discontinuous (\cref{R:discontinuity}).  Explicit
examples with smooth, admissible leader data are easy to construct (e.g.\ $\alpha=0.5$,
$\timehor=3$, a predicted gap $0.01(s-1)^2(2-s)$ realised by a suitable $\sa_\mL$: first
contact at $s=1$ with $\Lambda=0$, transversal crossing at $s=2$).  For this reason the
domain $\mathbb D$ on which the regularity theory of \cref{sec:ttc-regularity} is
developed requires the \emph{first} predicted contact to be transversal
(\cref{D:domainD}); on $\mathbb D$ the two notions coincide.
\end{remark}

\begin{lemma}[Bounds]\label{L:tau-bounds}
For $z\in\Om$, with $\hat v=\max\{\bar\zeta,{v_{\max}}\}$,
\[
\min\Bigl\{\frac{\bar X}{\hat v},\ \timehor-\bar t\Bigr\}\ \le\ \tau(z)\ \le\ \timehor-\bar t .
\]
In particular $\tau>0$ on $\Om$, and
$\tau(z)\ge\min\{\bar X/{v_{\max}},\timehor-\bar t\}$ for $z\in\Phys$.
\end{lemma}

\begin{proof}
The upper bound holds because $\Tt\le\timehor$.  By \cref{L:prediction-bounds}(a),
$\Xi(s;z)\ge\bar X-\hat v(s-\bar t)>0$ for $s<\bar t+\bar X/\hat v$.  Hence either the set
in \eqref{E:TTC} is empty, so that $\tau=\timehor-\bar t$, or all its elements are
$\ge\bar t+\bar X/\hat v$.
\end{proof}

\begin{lemma}[Upper semicontinuity; formula on $\Reg$]\label{L:usc}
$\Tt$ is upper semicontinuous on $\Om$; in particular $\Tt$ and $\tau$ are Borel.
For $z\in\Reg$ we have $\Tt(z)=\min\{\sigma(z),\timehor\}$.
\end{lemma}

\begin{proof}
Fix $z_0=(\bar t_0,\bar X_0,\bar V_0)\in\Om$ and $\varepsilon>0$.  If $\Tt(z_0)=\timehor$,
then $\Tt(z)\le\timehor=\Tt(z_0)$ for every $z$, and there is nothing to prove.

Otherwise there is $s'$ with $\Xi(s';z_0)=0$, $\Lambda(s';z_0)<0$ and
$s'<\Tt(z_0)+\varepsilon$.  Note $s'>\bar t_0$, because $\Xi(\bar t_0;z_0)=\bar X_0>0$.
By continuity choose $h\in(0,\min\{\varepsilon,1\})$, so that $s'+h\le\timehor+1$, with $s'-h>\bar t_0$, $\Lambda(\cdot;z_0)<0$ on
$[s'-h,s'+h]$, and $\Xi(s'-h;z_0)>0>\Xi(s'+h;z_0)$.  By \cref{L:flow}(b) there is $\rho>0$
such that for $\abs{z-z_0}<\rho$: $\bar t<s'-h$, $\Lambda(\cdot;z)<0$ on $[s'-h,s'+h]$,
and $\Xi(s'-h;z)>0>\Xi(s'+h;z)$.  The intermediate value theorem gives
$s_z\in(s'-h,s'+h)$ with $\Xi(s_z;z)=0$ and $\Lambda(s_z;z)<0$.
\begin{itemize}
\item If $s_z\le\timehor$, then $\Tt(z)\le s_z<s'+h<\Tt(z_0)+2\varepsilon$.
\item If $s_z>\timehor$, then $\Tt(z)\le\timehor<s_z<\Tt(z_0)+2\varepsilon$.
\end{itemize}
Hence $\limsup_{z\to z_0}\Tt(z)\le\Tt(z_0)$.  An upper semicontinuous function has open
sublevel sets $\{\Tt<c\}$, so it is Borel.

For $z\in\Reg$: if $\sigma(z)=\infty$, $\Xi$ has no zero in $[\bar t,\timehor]$ and
$\Tt=\timehor$.  If $\sigma(z)<\infty$, then $\sigma(z)$ belongs to the set in
\eqref{E:TTC} (because $\Lambda(\sigma)<0$) and is its smallest element, so
$\Tt=\sigma\le\timehor$.
\end{proof}

\begin{lemma}[$\Reg$ is open and $\tau$ is locally Lipschitz on $\Reg$]\label{L:R-open}
$\Reg$ is relatively open in $[0,\timehor)\times\RR^2$.  Every $z_0\in\Reg$ has a
neighbourhood on which $\tau$ is Lipschitz.  More precisely, let $B$ be the closed unit
ball around $z_0$ and $L_B$ the constant of \cref{L:flow}.  All radii $\rho$ below are
taken $\le1$, so that \eqref{E:flow-lip-z} applies.  Then:
\begin{enumerate}[label=\tup{\alph*}]
\item if $\sigma(z_0)=\infty$, then $\tau(z)=\timehor-\bar t$ near $z_0$;
\item if $\sigma(z_0)<\infty$ and $\Lambda(\sigma(z_0);z_0)=-2\gamma<0$, then
  $\mathrm{Lip}(\tau)\le L_B/\gamma+1$ near $z_0$.
\end{enumerate}
\end{lemma}

\begin{proof}
\emph{Case (a).}  $\Xi(\cdot;z_0)$ is continuous on $[0,\timehor+1]$, positive at
$\bar t_0$, and has no zero on $[\bar t_0,\timehor]$.  Hence it is positive on
$[\bar t_0,\timehor]$ and, by continuity, on $[\max\{\bar t_0-h',0\},\timehor]$ for some
$h'>0$.  Let $2c>0$ be its minimum there, and let $\rho\Def\min\{h',c/L_B\}$.  For
$\abs{z-z_0}<\rho$ we have $\bar t>\bar t_0-h'$, and \eqref{E:flow-lip-z} gives
$\Xi(s;z)\ge 2c-L_B\rho>0$ on $[\bar t,\timehor]$.  Thus $\sigma(z)=\infty$, $z\in\Reg$,
and $\tau(z)=\timehor-\bar t$, which is $1$-Lipschitz.

\emph{Case (b).}  Let $\sigma_0\Def\sigma(z_0)\in(\bar t_0,\timehor]$.  Choose
$h\in\bigl(0,\min\{1,(\sigma_0-\bar t_0)/2\}\bigr)$ such that
$\Lambda(s;z_0)\le-\tfrac32\gamma$ for $\abs{s-\sigma_0}\le h$.  This is possible because
$\Lambda(\cdot;z_0)$ is continuous, and $\sigma_0+h\le\timehor+1$.  Then
$\Xi(\sigma_0+h;z_0)\le-\tfrac32\gamma h$.  Since $\sigma_0$ is the first zero,
$\Xi(\cdot;z_0)>0$ on $[\bar t_0,\sigma_0-h]$, hence on $[\bar t_0-h',\sigma_0-h]$ for some
$h'\in(0,h)$ (intersected with $[0,\infty)$).  Let $2c>0$ be the minimum there.  Pick
$\rho\in(0,h')$ with $L_B\rho<\min\{c,\ \gamma/2,\ \tfrac34\gamma h\}$.  For
$\abs{z-z_0}<\rho$:
\begin{enumerate}[label=\tup{\roman*}]
\item $\bar t<\bar t_0+h'<\sigma_0-h$ and $\Xi(s;z)>c$ on $[\bar t,\sigma_0-h]$;
\item $\Lambda(s;z)\le-\gamma$ on $[\sigma_0-h,\sigma_0+h]$;
\item $\Xi(\sigma_0+h;z)<0$.
\end{enumerate}
By (i)--(iii) and strict monotonicity on $[\sigma_0-h,\sigma_0+h]$, the extended
prediction $\Xi(\cdot;z)$ has exactly one zero $\hat\sigma(z)$ in $[\bar t,\sigma_0+h]$.
It lies in $(\sigma_0-h,\sigma_0+h)$ and is transversal.  If $\hat\sigma(z)\le\timehor$
then $\sigma(z)=\hat\sigma(z)$; otherwise $\sigma(z)=\infty$.  In both cases $z\in\Reg$,
and by \cref{L:usc} $\Tt(z)=\min\{\hat\sigma(z),\timehor\}$.

For $z_1,z_2$ in this ball, write $\hat\sigma_i=\hat\sigma(z_i)$.  Then
$\Xi(\hat\sigma_1;z_1)=0=\Xi(\hat\sigma_2;z_2)$, and $\partial_s\Xi=\Lambda\le-\gamma$ on
the interval containing $\hat\sigma_1,\hat\sigma_2$.  By the mean value theorem and
\eqref{E:flow-lip-z},
\[
\gamma\abs{\hat\sigma_1-\hat\sigma_2}\le\abs{\Xi(\hat\sigma_1;z_1)-\Xi(\hat\sigma_2;z_1)}
=\abs{\Xi(\hat\sigma_2;z_2)-\Xi(\hat\sigma_2;z_1)}\le L_B\abs{z_1-z_2}.
\]
Since $r\mapsto\min\{r,\timehor\}$ is $1$-Lipschitz and $\tau=\Tt-\bar t$, we get
$\abs{\tau(z_1)-\tau(z_2)}\le(L_B/\gamma+1)\abs{z_1-z_2}$.

Cases (a) and (b) exhaust $\Reg$, and both produce a neighbourhood contained in $\Reg$.
So $\Reg$ is open.
\end{proof}

\begin{remark}[Discontinuity on $\Graz$]\label{R:discontinuity}
$\tau$ is in general discontinuous on $\Graz$.  Here is the typical picture.  Let
$z_0\in\Graz$, and suppose the predicted gap touches $0$ at $\sigma_0$ with
$\sigma_0+h\le\timehor$ and reopens, so that $\Xi(\cdot;z_0)>0$ on $(\sigma_0,\sigma_0+h]$.
Assume the touch is nondegenerate (e.g.\ $\operatorname{ess\,inf}(\sa_\mL+\alpha v_\mL)>0$
near $\sigma_0$).  Then $\tau(z_0)\ge\sigma_0+h-\bar t_0$.

The state gradient $D_{(\bar X,\bar V)}\Xi(\sigma_0;z_0)=(\phi_1(\sigma_0),\delta_0(\sigma_0))$,
the first row of the fundamental matrix of the variational equation of \eqref{E:prediction}
along the prediction (its Wronskian is $e^{-\alpha(\sigma_0-\bar t_0)}\neq0$), is nonzero;
no derivative with respect to the initial time $\bar t$ is used, and none exists in general
at a jump of $\sa_\mL$ (\cref{R:not-C1-in-time}).  So there is a state direction
$d=(0,d_X,d_V)$ with $D_{(\bar X,\bar V)}\Xi(\sigma_0;z_0)\cdot(d_X,d_V)<0$.
Along $z_0+\varepsilon d$ the predicted gap dips below $0$ near $\sigma_0$ and crosses
transversally, so $\tau\approx\sigma_0-\bar t_0$.  Note that the direction matters: the
linearised gap is a damped oscillator, and $\partial\Xi(\sigma_0)/\partial\bar X$ can be
negative; there are examples (small $\alpha$, long prediction horizon) where lowering
$\bar X$ \emph{removes} the contact and raising it creates one.

Hence, when $\bar V_0\neq0$, $\mathsf f$ jumps across $\Graz$, and no Lipschitz property
of the closed-loop field can hold across $\Graz$.  It holds on compact subsets of $\Reg$
(\cref{L:Cara-conditions}).  If instead $\sigma_0=\timehor$, both one-sided values tend to
$\timehor-\bar t_0$, and $\tau$ is continuous there.
\end{remark}

\begin{lemma}[TTC in the approach region]\label{L:approach}
Let $z=(\bar t,\bar X,\bar V)\in\Slab$ with
\[
\bar V<0\qquad\text{and}\qquad\bar X<\frac{\bar V^2}{2\abar}.
\]
Then $z\in\Reg$ and
\begin{equation}\label{E:tau-upper}
\tau(z)\le\frac{2\bar X}{\abs{\bar V}} .
\end{equation}
\end{lemma}

\begin{proof}
Since $z\in\Slab$, $\bar\zeta<{v_{\max}}+1$, so \cref{L:prediction-bounds}(b) gives
$\partial_s\Lambda\le\abar$ a.e.  Integrating twice, for $s\ge\bar t$,
\[
\Lambda(s;z)\le\bar V+\abar(s-\bar t),\qquad
\Xi(s;z)\le\bar X+\bar V(s-\bar t)+\tfrac{\abar}{2}(s-\bar t)^2 .
\]
Let $s^\ast\Def2\bar X/\abs{\bar V}$ and $s_1\Def\min\{\bar t+s^\ast,\timehor\}$.  Since
$\abar s^\ast=2\abar\bar X/\abs{\bar V}<\abs{\bar V}$, we have $\Lambda(s;z)<0$ for
$s\in[\bar t,s_1]$.

\emph{Case $\bar t+s^\ast\ge\timehor$.}  \Cref{L:tau-bounds} gives
$\tau\le\timehor-\bar t\le s^\ast$.  Moreover $\Lambda<0$ on $[\bar t,\timehor]$, so any
first zero in $[\bar t,\timehor]$ is transversal and $z\in\Reg$.

\emph{Case $\bar t+s^\ast<\timehor$.}  Then $\bar t+s^\ast$ lies in the domain of the
prediction, and
\[
\Xi(\bar t+s^\ast;z)\le\bar X\Bigl(\frac{2\abar\bar X}{\bar V^2}-1\Bigr)<0 .
\]
So $\Xi(\cdot;z)$ changes sign on $[\bar t,\bar t+s^\ast]$, and
$\sigma(z)\in(\bar t,\bar t+s^\ast)$ with $\Lambda(\sigma(z);z)<0$.  Thus $z\in\Reg$, and by \cref{L:usc}
$\tau=\sigma(z)-\bar t<s^\ast$.
\end{proof}


\subsection{Solution concepts for the closed loop}\label{S:closed-loop}

\begin{definition}[Solutions]\label{D:solutions}\label{def:weak_sol}
Let $J\subset[0,\timehor)$ be an interval with left endpoint $t_\circ\in J$, and let
$\xi=(X,V):J\to\RR^2$ be absolutely continuous on compact subintervals, with
$(t,\xi(t))\in\Om$ for $t\in J$ and $\xi(t_\circ)=\xi_\circ$.
\begin{enumerate}[label=\tup{\alph*}]
\item $\xi$ is a \emph{Carath\'eodory solution} of \eqref{E:closed-loop} if
  $\dot\xi(t)=\mathsf F(t,\xi(t))$ for a.e. $t\in J$; equivalently,
  $\xi(t)=\xi_\circ+\int_{t_\circ}^t\mathsf F(r,\xi(r))\,dr$ for all $t\in J$.
\item $\xi$ is a \emph{Filippov solution} if $\dot\xi(t)\in\Fil[\mathsf F](t,\xi(t))$ for
  a.e. $t\in J$, where
  \[
  \Fil[\mathsf F](t,\xi)\Def\bigcap_{\delta>0}\ \bigcap_{\abs N=0}\ \co\,\mathsf F\bigl(t,B(\xi,\delta)\setminus N\bigr),
  \]
  $\co$ denotes the closed convex hull and $\abs N$ the Lebesgue measure in $\RR^2$
  \cite[Ch.~2, \S4]{Filippov}.  Only $\delta$ small enough that
  $B(\xi,\delta)\subset\{X>0\}$ are used, so that $\mathsf F$ is defined on the ball.
\end{enumerate}
\end{definition}

\begin{lemma}[Two facts on Filippov sets]\label{L:Filippov-facts}
Fix $t$ and $\xi=(X,V)$.
\begin{enumerate}[label=\tup{\roman*}]
\item If $\mathsf f(t,\cdot)$ is continuous at $\xi$, then
  $\Fil[\mathsf F](t,\xi)=\{\mathsf F(t,\xi)\}$.
\item Let $W\ni\xi$ be open and $h:W\to\RR$ continuous, and suppose
  $\mathsf f(t,\eta)\ge h(\eta)$ for a.e. $\eta\in W$.  Then every
  $w=(w_1,w_2)\in\Fil[\mathsf F](t,\xi)$ satisfies $w_1=V$ and $w_2\ge h(\xi)$.  The same
  holds with ``$\le$'' in both places.
\end{enumerate}
\end{lemma}

\begin{proof}
(ii) Take $\delta>0$ with $B(\xi,\delta)\subset W$, and let
$\omega(\delta)\Def\sup_{\abs{\eta-\xi}<\delta}\abs{h(\eta)-h(\xi)}$.  Let $N$ be the null
set off which $\mathsf f(t,\cdot)\ge h$.  The values $\mathsf F(t,\eta)=(\eta_2,\mathsf f(t,\eta))$,
for $\eta\in B(\xi,\delta)\setminus N$, lie in the closed convex set
$\{w:\abs{w_1-V}\le\delta,\ w_2\ge h(\xi)-\omega(\delta)\}$.  So does their closed convex
hull, and hence $\Fil[\mathsf F](t,\xi)$.  Letting $\delta\downarrow0$ gives the claim.

(i) In the same way $\Fil[\mathsf F](t,\xi)$ is contained in every closed ball around
$\mathsf F(t,\xi)$ of radius
\[
\sup_{\abs{\eta-\xi}<\delta}\abs{\mathsf F(t,\eta)-\mathsf F(t,\xi)}\xrightarrow[\delta\downarrow0]{}0 .
\]
It also contains $\mathsf F(t,\xi)$.  By continuity at $\xi$, this point lies in the
closure of $\mathsf F(t,B(\xi,\delta)\setminus N)$ for every $\delta>0$ and every null set
$N$, because $B(\xi,\delta)\setminus N$ accumulates at $\xi$.  Hence
$\Fil[\mathsf F](t,\xi)=\{\mathsf F(t,\xi)\}$.
\end{proof}

\subsection{Local well-posedness in the regular set}

\begin{lemma}[Carath\'eodory conditions on $\Reg$]\label{L:Cara-conditions}
Let $K\subset\Reg$ be compact.  There are $m_K,L_K<\infty$ such that
\[
\abs{\mathsf F(t,\xi)}\le m_K\ \ \text{for a.e. }t,\qquad
\abs{\mathsf F(t,\xi_1)-\mathsf F(t,\xi_2)}\le L_K\abs{\xi_1-\xi_2},
\]
whenever $(t,\xi),(t,\xi_i)\in K$.  Moreover, for each $\xi$, $t\mapsto\mathsf F(t,\xi)$
is measurable on $\{t:(t,\xi)\in\Reg\}$.
\end{lemma}

\begin{proof}
\emph{$\tau$ on $K$.}  By \cref{L:R-open}, $\tau$ is locally Lipschitz, hence continuous,
on $\Reg$.  By \cref{L:tau-bounds} it is positive, so $\tau_K\Def\min_K\tau>0$.

$\tau$ is also Lipschitz on $K$.  If not, there would be $p_n,q_n\in K$ with
$\abs{\tau(p_n)-\tau(q_n)}>n\abs{p_n-q_n}$.  Since $0<\tau\le\timehor$, this forces
$\abs{p_n-q_n}<\timehor/n$.  Along a subsequence $p_n,q_n\to p\in K$, which contradicts the
Lipschitz bound of $\tau$ on a neighbourhood of $p$.  Let $L^\tau_K$ be the constant, and
$R_K\Def\max_K\abs V$.

\emph{Bounds on $\mathsf F$.}  By \eqref{E:V-range}, a.e. on $K$,
$\abs{\mathsf g}\le\norm{\sa_\mL}_\infty+\alpha(\Vbar+R_K+{v_{\max}})$, and
$\abs{\beta V/\tau^2}\le\beta R_K/\tau_K^2$.  For two points of $K$ with the same $t$,
$\sa_\mL(t)$ and $v_\mL(t)$ cancel, so
\[
\abs{\mathsf f(t,\xi_1)-\mathsf f(t,\xi_2)}
\le\alpha\abs{X_1-X_2}+\alpha\abs{V_1-V_2}
+\beta\Bigl(\frac{\abs{V_1-V_2}}{\tau_K^2}+R_K\,\frac{\abs{\tau_2^2-\tau_1^2}}{\tau_K^4}\Bigr),
\]
with $\abs{\tau_2^2-\tau_1^2}\le2\timehor L^\tau_K\abs{\xi_1-\xi_2}$.  This gives $L_K$.

\emph{Measurability.}  For fixed $\xi$, the set $\{t:(t,\xi)\in\Reg\}$ is relatively open
and $\tau(\cdot,\xi)$ is continuous on it.  Also $\sa_\mL$ is measurable and $v_\mL$
continuous.
\end{proof}

\begin{theorem}[Local well-posedness in $\Reg$]\label{T:local-WP}
Let $(t_\circ,\xi_\circ)\in\Reg$.  There is a unique maximal Carath\'eodory solution
$\xi:[t_\circ,T_\Reg)\to\RR^2$ with graph in $\Reg$.  Its graph eventually leaves every
compact subset of $\Reg$: for every compact $K\subset\Reg$ there is $t_K<T_\Reg$ with
$(t,\xi(t))\notin K$ for $t\in(t_K,T_\Reg)$.  Moreover, every Filippov solution
$\tilde\xi:J\to\RR^2$ starting at $(t_\circ,\xi_\circ)$ coincides with $\xi$ on the common
interval $J\cap[t_\circ,T_\Reg)$.
\end{theorem}

\begin{proof}
$\Reg$ is only relatively open in $[0,\timehor)\times\RR^2$, because it may contain points
with $t=0$.  To apply the classical theorems on an open set, extend $\sa_\mL$ by $0$ and
$v_\mL$ by $v_\mL(0)$ to $t\in(-1,0)$.  The prediction flow and $\tau$ are then defined for
$\bar t\in(-1,\timehor)$, all the lemmas above hold verbatim, and solutions starting at
$t_\circ\ge0$ are unaffected.  The same remark applies to \cref{T:global}.

By \cref{L:Cara-conditions}, $\mathsf F$ satisfies the Carath\'eodory conditions with local
Lipschitz bounds on the open set $\Reg$.  The classical existence, uniqueness and extension
theorem \cite[Ch.~I, \S5]{Hale} gives the first two statements.

For the last statement, let $\tilde\xi:J\to\RR^2$ be a Filippov solution starting at
$(t_\circ,\xi_\circ)$.  Let
$S\Def\sup\{s\in J:(t,\tilde\xi(t))\in\Reg\text{ for }t\in[t_\circ,s]\}$, which is
$>t_\circ$ because $\Reg$ is open.  On $[t_\circ,S)$, \cref{L:R-open} and
\cref{L:Filippov-facts}(i) give $\Fil[\mathsf F]=\{\mathsf F\}$ along $\tilde\xi$, so
$\tilde\xi$ is a Carath\'eodory solution in $\Reg$.  By uniqueness, $\tilde\xi=\xi$ on
$[t_\circ,\min\{S,T_\Reg\})$.

Suppose $S<T_\Reg$ and $S\in J$.  By continuity $\tilde\xi(S)=\xi(S)$, and
$(S,\xi(S))\in\Reg$.  Openness of $\Reg$ and continuity of $\tilde\xi$ then keep the graph of
$\tilde\xi$ in $\Reg$ slightly beyond $S$, contradicting the definition of $S$.  Hence
either $S\ge T_\Reg$, or $S=\sup J\notin J$.  In both cases $\tilde\xi=\xi$ on
$J\cap[t_\circ,T_\Reg)$.
\end{proof}

\subsection{Invariance of the velocity box}

\begin{lemma}[Maximum principle]\label{L:invariance}\label{T:max_principle}
Let $\xi$ be a Carath\'eodory or Filippov solution on $J$ with $(t_\circ,\xi_\circ)\in\Phys$.
Then $(t,\xi(t))\in\Phys$ for all $t\in J$; that is, $0\le v(t)\le{v_{\max}}$.  Consequently
$\abs{V}\le{v_{\max}}$, $X$ is ${v_{\max}}$-Lipschitz, and $X$ has a finite left limit at the right
endpoint of $J$.
\end{lemma}

\begin{proof}
The function $v=v_\mL-V$ is absolutely continuous, with $\dot v=\sa_\mL-\dot V$ a.e.

\emph{Upper bound.}  Fix $t$ and let
$W_+\Def\{\eta=(X,V):X>0,\ v_\mL(t)-V>{v_{\max}}\}$, which is open.  For $\eta\in W_+$, with
$v=v_\mL(t)-V$, we have $\sa_\mL(t)-\mathsf f(t,\eta)=\alpha(\mathcal V(X)-v)+\beta(v_\mL(t)-v)/\tau^2$.
The second term is $\le0$, because $v_\mL(t)\le{v_{\max}}<v$.  Hence
\[
\mathsf f(t,\eta)\ \ge\ \sa_\mL(t)-\alpha\bigl(\Vbar-(v_\mL(t)-V)\bigr)\Def h_+(\eta),
\]
which is continuous in $\eta$.  For a Carath\'eodory solution this holds pointwise.  For a
Filippov solution, \cref{L:Filippov-facts}(ii) gives $\dot V\ge h_+(\xi(t))$ at a.e. $t$
with $\xi(t)\in W_+$.  In both cases $\dot v\le\alpha(\Vbar-v)<0$ at a.e. $t$ with
$v(t)>{v_{\max}}$, using $\Vbar\le{v_{\max}}$.  With $u\Def(v-{v_{\max}})^+$ we have
$u'=\one_{\{v>{v_{\max}}\}}\dot v\le0$ a.e. and $u(t_\circ)=0$, so $u\equiv0$.

\emph{Lower bound.}  On $W_-\Def\{X>0,\ v_\mL(t)-V<0\}$, we have
$\sa_\mL(t)-\mathsf f=\alpha(\mathcal V(X)-v)+\beta(v_\mL(t)-v)/\tau^2\ge\alpha(\mathcal V(X)-v)$,
because $v_\mL(t)\ge0>v$.  Hence $\mathsf f\le\sa_\mL(t)-\alpha(\mathcal V(X)-v)$, a
continuous bound.  As above, $\dot v\ge\alpha(\mathcal V(X)-v)>0$ at a.e. $t$ with
$v(t)<0$, by \eqref{E:V-range} and $X>0$.  Then $(-v)^+$ is nonincreasing and vanishes at
$t_\circ$.

The remaining statements follow from $V=v_\mL-v\in[v_\mL-{v_{\max}},v_\mL]\subset[-{v_{\max}},{v_{\max}}]$
and $\dot X=V$.
\end{proof}

\subsection{Impact-freeness}\label{SS:impact}

Throughout this subsection $\xi=(X,V)$ is a Carath\'eodory or Filippov solution on
$[t_\circ,t^\ast)$ with $t^\ast\le\timehor$ and $(t_\circ,\xi_\circ)\in\Phys$.  By
\cref{L:invariance} its graph stays in $\Phys\subset\Slab$, and $X(t^\ast-)$ exists.  We
use the open \emph{approach region}
\begin{equation}\label{E:approach-region}
\App\Def\Bigl\{(t,X,V)\in\Slab:\ V<0,\ X<\frac{V^2}{2\abar}\Bigr\}\subset\Reg ,
\end{equation}
where the inclusion is \cref{L:approach}.

\begin{lemma}[One-sided acceleration bounds]\label{L:one-sided}
At a.e. $t\in[t_\circ,t^\ast)$:
\begin{enumerate}[label=\tup{\alph*}]
\item if $V(t)<0$, then $\dot V(t)\ge-\abar$;
\item if $V(t)>0$, then $\dot V(t)\le\abar$;
\item if $(t,\xi(t))\in\App$, then $\dot V(t)\ge-\abar+\beta\,\dfrac{\abs{V(t)}^3}{4X(t)^2}$.
\end{enumerate}
Consequently $t\mapsto V(t)+\abar t$ is nondecreasing on every interval where $V<0$, and
$t\mapsto V(t)-\abar t$ is nonincreasing on every interval where $V>0$.
\end{lemma}

\begin{proof}
On $\Slab$ we have $v=v_\mL(t)-V\in(-1,{v_{\max}}+1)$, so by \eqref{E:V-range}
\begin{equation}\label{E:g-bounds}
-\abar\ \le\ \mathsf g(t,X,V)=\sa_\mL(t)-\alpha\mathcal V(X)+\alpha v\ \le\ \abar
\qquad\text{at every point of }\Slab .
\end{equation}
(a) On the open set $\{\eta:(t,\eta)\in\Slab,\ \eta_2<0\}$ the TTC term is
$-\beta\eta_2/\tau^2>0$, so $\mathsf f(t,\eta)\ge\mathsf g(t,\eta)\ge-\abar$ at every point.
For a Carath\'eodory solution this gives $\dot V\ge-\abar$ wherever $V<0$.  For a Filippov
solution, \cref{L:Filippov-facts}(ii) with the constant function $h\equiv-\abar$ gives the
same.  (b) is symmetric: $\mathsf f\le\mathsf g\le\abar$ where $\eta_2>0$.  (c) On $\App$,
\cref{L:approach} gives $1/\tau^2\ge V^2/(4X^2)$, hence
\begin{equation}\label{E:f-lower}
\mathsf f(t,X,V)\ \ge\ -\abar+\beta\,\frac{\abs V^3}{4X^2}\ \Def\ k(X,V)
\qquad\text{at every point of }\App ;
\end{equation}
$k$ is continuous on $\{V<0,X>0\}$, and \cref{L:Filippov-facts}(ii) applies on the open set
$\{\eta:(t,\eta)\in\App\}$.  The monotonicity statements follow because $V$ is absolutely
continuous.
\end{proof}

\begin{theorem}[Every contact has zero closing speed]\label{T:no-impact}
If $X(t^\ast-)=0$, then $V(t)\to0$ as $t\uparrow t^\ast$.
\end{theorem}

\begin{proof}
\emph{Step 1: $\liminf_{t\uparrow t^\ast}V(t)\ge0$.}  Exactly one of the following holds.

\emph{Case A: $V<0$ on some $(t^\ast-\delta,t^\ast)$.}  By \cref{L:one-sided},
$V+\abar t$ is nondecreasing there.  Since $\abs V\le{v_{\max}}$, the left limit
$V^\ast\Def V(t^\ast-)$ exists.  Suppose $V^\ast<0$.  Then there is $t_1<t^\ast$ with
$V\le V^\ast/2$ on $[t_1,t^\ast)$ and, since $X\to0$, with $X<V^{\ast2}/(8\abar)\le V^2/(2\abar)$
there.  So $(t,\xi(t))\in\App$ on $[t_1,t^\ast)$, and \cref{L:one-sided}(c) gives, a.e.,
\[
\dot V\ \ge\ -\abar+\beta\,\frac{\abs{V^\ast}^3}{32\,X^2}.
\]
Because $X(t^\ast-)=0$ and $X$ is ${v_{\max}}$-Lipschitz (\cref{L:invariance}),
$X(s)\le{v_{\max}}(t^\ast-s)$, so
\[
V(t)-V(t_1)\ \ge\ -\abar(t-t_1)+\frac{\beta\abs{V^\ast}^3}{32\,{v_{\max}}^2}
\Bigl(\frac{1}{t^\ast-t}-\frac{1}{t^\ast-t_1}\Bigr)\ \longrightarrow\ +\infty
\qquad(t\uparrow t^\ast),
\]
contradicting $V\le{v_{\max}}$.  Hence $V^\ast=0$.

\emph{Case B: there are $t_n\uparrow t^\ast$ with $V(t_n)\ge0$.}  Fix $\varepsilon>0$ and
$n$ with $\abar(t^\ast-t_n)<\varepsilon$.  Let $t\in(t_n,t^\ast)$ with $V(t)<0$, and set
$s_2\Def\sup\{s\in[t_n,t]:V(s)\ge0\}$.  By continuity $V(s_2)=0$ and $V<0$ on $(s_2,t]$.
\Cref{L:one-sided} gives $V(t)+\abar t\ge V(s_2)+\abar s_2$, i.e.
\[
V(t)\ \ge\ -\abar(t-s_2)\ \ge\ -\abar(t^\ast-t_n)\ >\ -\varepsilon .
\]
So $V>-\varepsilon$ on $(t_n,t^\ast)$, and $\liminf V\ge0$.

\emph{Step 2: $\limsup_{t\uparrow t^\ast}V(t)\le0$.}  Let $t<t^\ast$ with $X(t)<X(t_\circ)$
and $c\Def V(t)>0$.  $V$ cannot be positive on all of $[t_\circ,t]$, since that would give
$X(t)>X(t_\circ)$.  So $t'\Def\sup\{s\le t:V(s)\le0\}$ satisfies $V(t')=0$ and $V>0$ on
$(t',t]$.  By \cref{L:one-sided}(b), $V(s)\ge c-\abar(t-s)$ on $[t',t]$; at $s=t'$ this
gives $t-t'\ge c/\abar$.  Hence
\begin{equation}\label{E:upper-barrier}
X(t)\ \ge\ X(t')+\int_{t-c/\abar}^{t}\bigl(c-\abar(t-s)\bigr)\,ds\ \ge\ \frac{c^2}{2\abar},
\qquad\text{i.e.}\qquad V(t)\le\sqrt{2\abar X(t)} .
\end{equation}
Since $X(t)\to0$, $\limsup V\le0$.  Together with Step~1, $V\to0$.
\end{proof}

The following quantitative version is not needed for the safety claim, but it describes
the approach to contact precisely.

\begin{proposition}[Parabolic barrier near contact]\label{P:barrier}
If $X(t^\ast-)=0$, there is $t_4<t^\ast$ such that
\begin{equation}\label{E:barrier}
\abs{V(t)}\ \le\ \sqrt{2\abar\,X(t)}\qquad\text{for all }t\in[t_4,t^\ast).
\end{equation}
Moreover, the trajectory cannot enter $\App$ through a point with $X\le\beta^2\abar/8$: if
$(t_1,\xi(t_1))\notin\App$ and $X(t_1)\le\beta^2\abar/8$, there is no $t>t_1$ with
$(s,\xi(s))\in\App$ for all $s\in(t_1,t]$.  In particular, on any interval throughout which
$X<\beta^2\abar/8$, a trajectory that is outside $\App$ at some time stays outside $\App$
for the rest of that interval; this is the form \eqref{E:no-reentry} used below.
\end{proposition}

The proof uses the functions
\begin{equation}\label{E:Phi}
\Phi(x)\Def\frac\beta4\,\frac1x-\frac{1}{\sqrt{2\abar}}\,\frac{1}{\sqrt x},\qquad
\psi(x)\Def\Phi(x)-\frac{1}{\sqrt{2\abar x}}=\frac\beta4\,\frac1x-\frac{2}{\sqrt{2\abar}}\,\frac1{\sqrt x},
\qquad x>0 .
\end{equation}

\begin{proof}
\emph{Step 1: a differential inequality in $\App$.}
Let $I\subset[t_\circ,t^\ast)$ be an interval on which $(t,\xi(t))\in\App$.  By
\cref{L:one-sided}(c), $\dot V\ge k(X,V)$ a.e. on $I$.  On $I$, $V<0$, so
$w\Def1/\abs V=-1/V$ is absolutely continuous on compact subintervals, with
\[
\dot w=\frac{\dot V}{V^2}\ \ge\ -\frac{\abar}{V^2}+\frac{\beta}{4}\,\frac{\abs V}{X^2}
\qquad\text{a.e. on }I .
\]
On $\App$, $V^2>2\abar X$, so $\abar/V^2<1/(2X)$.  Since $\dot X=V=-\abs V$, we have
$\abs V/X^2=\tfrac{d}{dt}(1/X)$.  Therefore
\begin{equation}\label{E:w-ineq}
\dot w\ \ge\ \frac\beta4\,\frac{d}{dt}\Bigl(\frac1X\Bigr)-\frac{1}{2X}\qquad\text{a.e. on }I .
\end{equation}

\emph{Step 2: integrated form.}  Let $[t_1,t_2]\subset[t_\circ,t^\ast)$ with
$(t,\xi(t))\in\App$ for $t\in(t_1,t_2]$, and $V(t_1)<0$.  We claim
\begin{equation}\label{E:w-integrated}
w(t_2)-w(t_1)\ \ge\ \Phi(X(t_2))-\Phi(X(t_1)).
\end{equation}
On $[t_1,t_2]$, $V$ is continuous and nonzero, so $w$ is absolutely continuous there, and
$X\ge X(t_2)>0$.  Also $\abs V>\sqrt{2\abar X}$ on $(t_1,t_2]$.  By the chain rule for
$x\mapsto x^{-1/2}$ composed with $X$,
\[
\begin{aligned}
\int_{t_1}^{t_2}\frac{dt}{2X(t)}
&=\int_{t_1}^{t_2}\frac{\abs{\dot X(t)}}{2X(t)\abs{V(t)}}\,dt\\
&\le\int_{t_1}^{t_2}\frac{-\dot X(t)}{2\sqrt{2\abar}\,X(t)^{3/2}}\,dt\\
&=\frac{1}{\sqrt{2\abar}}\Bigl(X(t_2)^{-1/2}-X(t_1)^{-1/2}\Bigr).
\end{aligned}
\]
Integrating \eqref{E:w-ineq} over $[t_1,t_2]$ and using this bound gives
\eqref{E:w-integrated}.

\emph{Step 3: properties of $\psi$.}  We have
$\psi'(x)=-\tfrac\beta4x^{-2}+\tfrac{1}{\sqrt{2\abar}}x^{-3/2}$, which is $<0$ iff
$x<x_1\Def\beta^2\abar/8$.  So $\psi$ is strictly decreasing on $(0,x_1]$.  Moreover
\begin{equation}\label{E:Phi-blowup}
\psi(x)\ \longrightarrow\ +\infty\qquad(x\downarrow0).
\end{equation}
Since $w<1/\sqrt{2\abar X}$ in $\App$, \eqref{E:w-integrated} implies, under the hypotheses
of Step~2,
\begin{equation}\label{E:psi-ineq}
\psi(X(t_2))\ <\ \Phi(X(t_1))-w(t_1).
\end{equation}

\emph{Step 4: the solution cannot stay in $\App$ up to $t^\ast$.}  Suppose
$(t,\xi(t))\in\App$ for all $t\in[t_1,t^\ast)$.  Applying \eqref{E:psi-ineq} on
$[t_1,t]$ for every $t\in(t_1,t^\ast)$ gives $\psi(X(t))<\Phi(X(t_1))-w(t_1)$.  The
right-hand side is a fixed number, while the left-hand side tends to $+\infty$ as
$t\uparrow t^\ast$, by \eqref{E:Phi-blowup} and $X(t)\to0$.  This is a contradiction.

\emph{Step 5: once close to contact, no re-entry into $\App$.}  Choose $t_3<t^\ast$ such
that $X(t)<\min\{x_1,X(t_\circ)\}$ for all $t\in[t_3,t^\ast)$.  By Step~4 there is
$t_4\in[t_3,t^\ast)$ with $(t_4,\xi(t_4))\notin\App$.  We claim
\begin{equation}\label{E:no-reentry}
(t,\xi(t))\notin\App\qquad\text{for all }t\in[t_4,t^\ast).
\end{equation}
Suppose not, and pick $t\in(t_4,t^\ast)$ with $(t,\xi(t))\in\App$.  Let
$t_1\Def\sup\{s\in[t_4,t]:(s,\xi(s))\notin\App\}$.  This supremum is attained, because the
set of such $s$ is the preimage of a closed set under a continuous map.  Then $t_1<t$ and
$(s,\xi(s))\in\App$ for $s\in(t_1,t]$.

By continuity, $V(t_1)\le0$ and $X(t_1)\le V(t_1)^2/(2\abar)$.  If $V(t_1)=0$ this forces
$X(t_1)\le0$, which is impossible.  So $V(t_1)<0$.  Since
$(t_1,\xi(t_1))\in\Phys\setminus\App\subset\Slab\setminus\App$, we get
$X(t_1)=V(t_1)^2/(2\abar)$, that is, $w(t_1)=1/\sqrt{2\abar X(t_1)}$, and the right-hand
side of \eqref{E:psi-ineq} equals $\psi(X(t_1))$.  So \eqref{E:psi-ineq} on $[t_1,t]$ reads
$\psi(X(t))<\psi(X(t_1))$.

But $X$ is strictly decreasing on $(t_1,t]$, because $V<0$ there, so
$X(t)<X(t_1)<x_1$.  Step~3 then gives $\psi(X(t))>\psi(X(t_1))$, a contradiction.  This
proves \eqref{E:no-reentry}.  The argument used only $X(t_1)\le x_1$ at the entry time
$t_1$ (so that $\psi$ is strictly decreasing on $(0,X(t_1)]$) and $V<0$ inside $\App$;
it therefore proves the entry statement of the proposition for an arbitrary entry point
with $X(t_1)\le x_1$, not only on $[t_3,t^\ast)$.  It does \emph{not} show that a
trajectory which leaves $\App$ with $X<x_1$, then opens the gap beyond $x_1$, cannot
re-enter $\App$ later through a point with $X>x_1$.

\emph{Step 6: the barrier.}  Let $t\in[t_4,t^\ast)$.  If $V(t)<0$, then by
\eqref{E:no-reentry} and $(t,\xi(t))\in\Phys\subset\Slab$, $X(t)\ge V(t)^2/(2\abar)$,
i.e.\ $V(t)\ge-\sqrt{2\abar X(t)}$.  If $V(t)>0$, \eqref{E:upper-barrier} applies because
$X(t)<X(t_\circ)$.  This is \eqref{E:barrier}.
\end{proof}

\begin{remark}[Interpretation]\label{R:what-used}
The parabola $X=V^2/(2\abar)$ is the stopping-distance curve at maximal deceleration
$\abar$.  Once $X\le\beta^2\abar/8$, the vector field points out of $\App$ across it:
there \eqref{E:f-lower} gives $\dot V\ge-\abar+\beta(2\abar)^{3/2}/(4\sqrt X)\ge\abar$.
On any time interval spent in $\App$, a Filippov
solution is a Carath\'eodory solution, by $\App\subset\Reg$ and \cref{L:Filippov-facts}(i).
We do not claim uniqueness near contact from this alone: by \eqref{E:no-reentry} the
solution is outside $\App$ just before contact, and there it may meet $\Graz$.  Uniqueness
is addressed in \cref{SS:crossing}.
\end{remark}

\subsection{Global existence and behaviour at the horizon}

\begin{theorem}[Filippov existence and continuation]\label{T:global}
Let $(t_\circ,\xi_\circ)\in\Phys$.  A \emph{non-extendable} Filippov solution is one that
is not the restriction of a Filippov solution defined on a strictly larger interval.
Non-extendable Filippov solutions with initial value $\xi_\circ$ at $t_\circ$ exist, and
every such solution $\xi$, on its maximal interval $[t_\circ,T^\ast)$ with
$T^\ast=T^\ast_\xi\in(t_\circ,\timehor]$, has the following properties.
\begin{enumerate}[label=\tup{\alph*}]
\item Either $T^\ast=\timehor$, or $X(T^\ast-)=0$.  In the latter case, $V(t)\to0$ as
  $t\uparrow T^\ast$ (\cref{T:no-impact}), and the barrier \eqref{E:barrier} holds near
  $T^\ast$ (\cref{P:barrier}).
\item $(t,\xi(t))\in\Phys$ for all $t$.
\item If $(t_\circ,\xi_\circ)\in\Reg$, then every Filippov solution coincides with the
  unique Carath\'eodory solution of \cref{T:local-WP} on their common interval of
  existence.  In particular $T^\ast\ge T_\Reg$, and every non-extendable Filippov solution
  starting at $(t_\circ,\xi_\circ)$ extends the Carath\'eodory one.
\end{enumerate}
\end{theorem}

\begin{proof}
\emph{Measurability and local boundedness of $\mathsf F$ on $\Om$.}  $\tau$ is Borel on
$\Om$ (\cref{L:usc}), $\sa_\mL$ is Lebesgue measurable and $v_\mL$ is continuous, so
$\mathsf F$ is Lebesgue measurable on $\Om$.  Let $K\subset\Om$ be compact.  Then on $K$
we have $X\ge X_K>0$, $t\le\timehor-\delta_K$ and $\abs V\le R_K$.  Hence
$\hat v\le{v_{\max}}+R_K$, and \cref{L:tau-bounds} gives
$\tau\ge\tau_K\Def\min\{X_K/({v_{\max}}+R_K),\delta_K\}>0$.  Therefore
$\abs{\mathsf F}\le R_K+\norm{\sa_\mL}_\infty+\alpha(\Vbar+R_K+{v_{\max}})+\beta R_K/\tau_K^2$
a.e. on $K$.

\emph{Local existence.}  We use Filippov's existence theorem for his definition~(a),
which convexifies in the state variable only.  It requires $\mathsf F$ to be Lebesgue
measurable on an open set of $(t,\xi)$ and locally essentially bounded, and gives a
Filippov solution through every point \cite{Filippov60}, \cite[Ch.~2, \S7]{Filippov}
(the version for a right-hand side that is merely measurable in $t$; the set-valued
formulation under ``basic conditions'' requires joint upper semicontinuity, which fails
here since $\sa_\mL$ is only $L^\infty$).  The open-set requirement is met by the extension
to $t\in(-1,0)$ described in the proof of \cref{T:local-WP}.

\emph{Continuation.}  By Zorn's lemma, applied to Filippov solutions through
$(t_\circ,\xi_\circ)$ ordered by extension, every Filippov solution is the restriction of
a non-extendable one; in particular non-extendable solutions exist.  Now let $\xi$ be
\emph{any} non-extendable solution, on $[t_\circ,T^\ast)$.  Suppose its graph does not
eventually leave some compact $K\subset\Om$.  Then there are $t_k\uparrow T^\ast$ with $(t_k,\xi(t_k))\in K$, and we may
assume $T^\ast<\timehor$.

Let $K'\supset K$ be a compact neighbourhood with $\abs{\mathsf F}\le m$ a.e. on $K'$, and
let $d>0$ be the distance from $K$ to the complement of $K'$.  Then
$\abs{\dot\xi}\le m$ as long as the graph stays in $K'$.  So for $k$ large, with
$T^\ast-t_k<d/(1+m)$, the graph remains in $K'$ on $[t_k,T^\ast)$, and
$\xi(T^\ast-)\Def\lim_{t\uparrow T^\ast}\xi(t)$ exists by the Cauchy criterion.  Moreover
$(T^\ast,\xi(T^\ast-))\in K'\subset\Om$.  Restarting the local existence theorem at this
point extends the solution, contradicting non-extendability.  Hence the graph eventually
leaves every compact subset of $\Om$.

\emph{(b)} is \cref{L:invariance}.

\emph{(a).}  Suppose $T^\ast<\timehor$ and $X(T^\ast-)>0$; the limit exists by
\cref{L:invariance}.  Then $X_{\min}\Def\inf_{[t_\circ,T^\ast)}X>0$, and the graph lies in
the compact set
$[t_\circ,T^\ast]\times[X_{\min},X_\circ+{v_{\max}}\timehor]\times[-{v_{\max}},{v_{\max}}]\subset\Om$,
a contradiction.  The statements on $V$ are \cref{T:no-impact,P:barrier}.

\emph{(c)} is the last assertion of \cref{T:local-WP}.  For $T^\ast\ge T_\Reg$: if
$T^\ast<T_\Reg$, the Carath\'eodory solution, which is also a Filippov solution, would
extend $\xi$ beyond $T^\ast$.
\end{proof}

\begin{proposition}[Behaviour at the horizon]\label{P:horizon}
In the situation of \cref{T:global}, if $T^\ast=\timehor$, then $V(t)\to0$ as
$t\uparrow\timehor$.  Hence $\xi$ extends continuously to $[t_\circ,\timehor]$ with
$V(\timehor)=0$.
\end{proposition}

\begin{proof}
By \cref{L:tau-bounds}, $\tau\le\timehor-t$.  On $\Slab$ we have
$\abs{\mathsf g}\le\norm{\sa_\mL}_\infty+\alpha(\Vbar+{v_{\max}}+1)=\abar$, by \eqref{E:V-range}.
Hence, for every point of $\Slab$ with $V\ne0$,
\[
\operatorname{sgn}(V)\,\mathsf f(t,X,V)=\operatorname{sgn}(V)\,\mathsf g-\beta\frac{\abs V}{\tau^2}
\ \le\ \abar-\beta\,\frac{\abs V}{(\timehor-t)^2},
\]
with a right-hand side continuous in $(X,V)$ on $\{V>0\}$ and on $\{V<0\}$.  By
\cref{L:Filippov-facts}(ii) (or directly, for Carath\'eodory solutions), and since
$\tfrac{d}{dt}\abs V=\operatorname{sgn}(V)\dot V$ a.e. on $\{V\ne0\}$, we get
\begin{equation}\label{E:absV-ineq}
\frac{d}{dt}\abs{V}\ \le\ \abar-\beta\,\frac{\abs V}{(\timehor-t)^2}\qquad\text{a.e. on }\{V\neq0\}.
\end{equation}
Fix $\eta>0$ and $t_\eta\Def\timehor-\sqrt{\beta\eta/(2\abar)}$.  Let
$u\Def(\abs V-\eta)^+$.  On $(t_\eta,\timehor)\cap\{\abs V>\eta\}$, \eqref{E:absV-ineq}
gives $\tfrac{d}{dt}\abs V\le\abar-2\abar<0$, so $u$ is nonincreasing on
$(t_\eta,\timehor)$.

If $u\ge c>0$ near $\timehor$, then \eqref{E:absV-ineq} gives
$\tfrac{d}{dt}\abs V\le\abar-\beta(\eta+c)/(\timehor-t)^2$ there.  This is not integrable
up to $\timehor$, so $\abs V$ would become negative, which is absurd.  Hence $u\to0$, that
is, $\limsup\abs V\le\eta$.  Since $\eta>0$ is arbitrary, $V\to0$.  Finally $X$ has a
limit because it is Lipschitz.
\end{proof}

\begin{remark}[Uniqueness near the horizon]\label{R:horizon-unique}
If $X(\timehor-)>0$, then for $t$ close to $\timehor$ we have $\timehor-t<X(t)/{v_{\max}}$.
By \cref{L:prediction-bounds}(a), $\Xi(s)\ge X(t)-{v_{\max}}(s-t)>0$ for $s\in[t,\timehor]$, so
$\sigma=\infty$, the state lies in $\Reg$, and $\tau=\timehor-t$ exactly.  So on a final
interval the solution is the unique Carath\'eodory solution, and $\tau$ is smooth along it.
\end{remark}

\subsection{Crossing the grazing set}\label{SS:crossing}

\Cref{T:global} gives existence through $\Graz$ in Filippov's sense, and uniqueness only
while the trajectory stays in $\Reg$.  This subsection shows that at a \emph{nondegenerate}
grazing state the trajectory crosses $\Graz$ transversally, spends no time on it, and
remains unique.  The mechanism is structural: the prediction changes along the closed loop
only through the TTC term, and the sign of that term is the same on both sides of $\Graz$.

\paragraph{Three facts on the prediction flow}

\begin{lemma}[Cocycle property]\label{L:cocycle}
For $z=(\bar t,\bar Z)$ and $s_1\in[0,\timehor+1]$, let $z^{s_1}\Def(s_1,Z(s_1;z))$.  Then
$Z(s;z^{s_1})=Z(s;z)$ for all $s\in[0,\timehor+1]$.
\end{lemma}
\begin{proof}
Both sides solve \eqref{E:prediction} and take the value $Z(s_1;z)$ at $s=s_1$; the
solution is unique by \cref{L:flow}(a), whose proof covers initial times in
$[0,\timehor+1]$ verbatim.  No time-shift invariance is used.
\end{proof}

\begin{lemma}[Difference of two predictions]\label{L:difference-ODE}
Let $z_a,z_b\in[0,\timehor)\times\RR^2$ and $D(s)\Def\Xi(s;z_a)-\Xi(s;z_b)$.  Then
$D\in C^1[0,\timehor+1]$, $D'=\Lambda(\cdot;z_a)-\Lambda(\cdot;z_b)$ is absolutely
continuous, and
\begin{equation}\label{E:D-ODE}
\begin{aligned}
D''+\alpha D'+\alpha\,\tilde c(s)\,D&=0\quad\text{a.e.},\\
\tilde c(s)&\Def\int_0^1\mathcal V'\bigl(\theta\,\Xi(s;z_a)+(1-\theta)\,\Xi(s;z_b)\bigr)\,d\theta\in(0,1].
\end{aligned}
\end{equation}
The leader data $\sa_\mL$, $v_\mL$ do not appear in \eqref{E:D-ODE}.
\end{lemma}
\begin{proof}
From \eqref{E:prediction},
$D''=\partial_s\Lambda(\cdot;z_a)-\partial_s\Lambda(\cdot;z_b)
=-\alpha[\mathcal V(\Xi_a)-\mathcal V(\Xi_b)]-\alpha[\Lambda_a-\Lambda_b]$ a.e.; the terms
$\sa_\mL(s)$ and $\alpha v_\mL(s)$ cancel.  The mean value theorem gives
$\mathcal V(\Xi_a)-\mathcal V(\Xi_b)=\tilde c\,D$ with $\tilde c\in(0,1]$ because
$0<\mathcal V'\le1$.  Continuity of $\tilde c$ follows from that of $\Xi(\cdot;z_a)$,
$\Xi(\cdot;z_b)$ and $\mathcal V'$.
\end{proof}

\begin{lemma}[Fundamental solutions]\label{L:fundamental}
Let $c:[0,\timehor+1]\to(0,1]$ be continuous and $t_0\in[0,\timehor+1]$.  Let
$\phi_1,\phi_2$ solve $\phi''+\alpha\phi'+\alpha c\phi=0$ with
$(\phi_1,\phi_1')(t_0)=(1,0)$ and $(\phi_2,\phi_2')(t_0)=(0,1)$.
\begin{enumerate}[label=\tup{\alph*}]
\item $\abs{\phi_i(s)}+\abs{\phi_i'(s)}\le e^{(1+2\alpha)\abs{s-t_0}}$ for $i=1,2$.
\item If $\alpha\ge4$, then $\phi_2(s)\ge(s-t_0)\,e^{-\alpha(s-t_0)/2}>0$ for all $s>t_0$.
\item If $\alpha<4$, set $\omega_\alpha\Def\sqrt{\alpha(1-\alpha/4)}$.  Then
  $\phi_2(s)\ge e^{-\alpha(s-t_0)/2}\,\sin(\omega_\alpha(s-t_0))/\omega_\alpha>0$ for
  $s\in(t_0,t_0+\pi/\omega_\alpha)$.
\end{enumerate}
\end{lemma}
\begin{proof}
(a) is Gr\"onwall for the first-order system $(\phi,\phi')$ in the $\ell^1$-norm: the
column sums of the matrix $\begin{pmatrix}0&1\\-\alpha c&-\alpha\end{pmatrix}$ are
$\le\max\{\alpha c,1+\alpha\}\le1+2\alpha$.

(b), (c).  Put $u(s)\Def\phi_2(s)e^{\alpha(s-t_0)/2}$.  A direct computation gives
$u''=\bigl(\tfrac{\alpha^2}{4}-\alpha c(s)\bigr)u$, with $u(t_0)=0$, $u'(t_0)=1$.

If $\alpha\ge4$, then $\tfrac{\alpha^2}4-\alpha c\ge\tfrac{\alpha^2}4-\alpha\ge0$.  Let
$s_1>t_0$ be the first zero of $u$ after $t_0$, if any.  On $(t_0,s_1)$, $u>0$, so
$u''\ge0$, so $u'\ge1$, so $u(s)\ge s-t_0$.  Then $u(s_1)\ge s_1-t_0>0$, a contradiction.
Hence $u(s)\ge s-t_0$ for all $s>t_0$, which is (b).

If $\alpha<4$, then $\tfrac{\alpha^2}4-\alpha c\ge-\omega_\alpha^2$.  Let
$w(s)\Def\sin(\omega_\alpha(s-t_0))/\omega_\alpha$, so $w''=-\omega_\alpha^2w$,
$w(t_0)=0$, $w'(t_0)=1$, and $w>0$ on $(t_0,t_0+\pi/\omega_\alpha)$.  While $u>0$,
\[
(u'w-uw')'=u''w-uw''=\bigl(\tfrac{\alpha^2}4-\alpha c+\omega_\alpha^2\bigr)uw\ \ge\ 0 ,
\]
and $u'w-uw'=0$ at $t_0$.  Hence $(u/w)'=(u'w-uw')/w^2\ge0$ where $u>0$ and $w>0$.  Since
$u/w\to1$ as $s\downarrow t_0$, we get $u\ge w$ as long as $u>0$ and
$s<t_0+\pi/\omega_\alpha$.  If $u$ had a first zero $s_1<t_0+\pi/\omega_\alpha$, then
$u(s_1)\ge w(s_1)>0$, a contradiction.  This is (c).
\end{proof}

\paragraph{Nondegenerate grazing states}

\begin{definition}[Nondegenerate grazing state]\label{D:nondegenerate}
Let $z_0=(\bar t_0,\bar X_0,\bar V_0)\in\Graz\cap\Slab$ with $\sigma_0\Def\sigma(z_0)$, so that
$\bar t_0<\sigma_0$, $\Xi(\sigma_0;z_0)=0=\Lambda(\sigma_0;z_0)$.  Let
$\delta_0$ be the solution of the variational equation
\begin{equation}\label{E:variational}
\delta''+\alpha\delta'+\alpha\,\mathcal V'\bigl(\Xi(s;z_0)\bigr)\,\delta=0,\qquad
\delta(\bar t_0)=0,\quad\delta'(\bar t_0)=1 ,
\end{equation}
which is the derivative $\partial_{\bar V}\Xi(\cdot;z_0)$ of the predicted gap with
respect to the initial relative velocity.  (Indeed, by \cref{L:difference-ODE},
$D_\varepsilon\Def\Xi(\cdot;\bar t_0,\bar X_0,\bar V_0+\varepsilon)-\Xi(\cdot;z_0)$ solves
\eqref{E:D-ODE} with $D_\varepsilon(\bar t_0)=0$, $D_\varepsilon'(\bar t_0)=\varepsilon$, and
$\tilde c_\varepsilon\to\mathcal V'(\Xi(\cdot;z_0))$ uniformly by \eqref{E:flow-lip-z}; so
$D_\varepsilon/\varepsilon\to\delta_0$ by Gr\"onwall.)  We call $z_0$ \emph{nondegenerate} if
\begin{enumerate}[label=\tup{ND\arabic*},ref=(ND\arabic*)]
\item\label{ND1} $\sigma_0<\timehor$, and there are $h_0,c_0>0$ with
  $\sa_\mL(s)+\alpha v_\mL(s)\ge2c_0$ for a.e. $s\in[\sigma_0-h_0,\sigma_0+h_0]$;
\item\label{ND2} $\delta_0(\sigma_0)\neq0$;
\item\label{ND3} $z_0^{+}\Def(\sigma_0+h,Z(\sigma_0+h;z_0))\in\Reg$ for all sufficiently
  small $h>0$;
\item\label{ND4} $\bar V_0\neq0$.
\end{enumerate}
We write $\varsigma\Def-\operatorname{sgn}(\bar V_0)\cdot\operatorname{sgn}(\delta_0(\sigma_0))\in\{\pm1\}$.
\end{definition}

\begin{remark}[Meaning of the conditions]\label{R:ND-meaning}
\begin{enumerate}[label=\tup{\alph*}, leftmargin= 0.20in]
\item At the graze, $\partial_s\Lambda=\mathsf g(s,0,0)=\sa_\mL(s)+\alpha v_\mL(s)$,
  since $\mathcal V(0)=0$.  So \ref{ND1} implies that the predicted relative velocity passes
  through zero with positive slope: the predicted gap has a strict, quadratic minimum at
  $\sigma_0$.  It excludes zero-speed down-crossings (where the gap would go negative) and
  flat contacts.  Physically: at the predicted graze time, the leader is not braking harder
  than $\alpha v_\mL$.  For a stopped leader ($v_\mL\equiv0$, $\sa_\mL\equiv0$ near
  $\sigma_0$) the point $(0,0)$ is an equilibrium of the prediction, so by backward
  uniqueness no first-contact graze can occur there at all; \ref{ND1} only restricts
  leaders that brake at the predicted touch.
\item \ref{ND2} says that the predicted gap at the graze time responds to the current
  relative velocity.  \Cref{L:fundamental} gives $\delta_0(\sigma_0)>0$ whenever
  $\alpha\ge4$, or $\alpha<4$ and $\sigma_0-\bar t_0<\pi/\omega_\alpha$.  Outside this sufficient range the comparison does not fix the sign; the nondegeneracy hypothesis allows either sign. The crossing direction below is defined using that sign.
\item \ref{ND3} excludes a second graze on the remainder of the prediction.  By
  \cref{L:cocycle}, the prediction from $z_0^+$ is the continuation of the prediction from
  $z_0$, so \ref{ND3} says: after the graze, the predicted gap either stays positive up to
  $\timehor$ or crosses zero transversally.  Under \ref{ND1}, $\Xi(\cdot;z_0)>0$ on
  $(\sigma_0,\sigma_0+h_1]$ for some $h_1>0$, so by \cref{L:cocycle} the first zero of
  $\Xi(\cdot;z_0)$ on $[\sigma_0+h,\timehor]$, and $\Lambda$ there, are the same for all
  $h\in(0,h_1]$.  Hence the condition holds for all $h\in(0,h_1]$ as soon as it holds for
  one.
\item \ref{ND4}: if $\bar V_0=0$, the TTC term vanishes at $z_0$ and the field is
  continuous there; the crossing direction is then not determined at first order.
\item The excluded cases $\sigma_0=\timehor$, $\delta_0(\sigma_0)=0$, a double graze and
  $\bar V_0=0$ are each one further scalar condition on top of $z_0\in\Graz$.  \ref{ND1}
  is of a different nature: it is an essential lower bound on the forcing on a whole
  neighbourhood of $\sigma_0$, and it can fail for admissible $L^\infty$ forcing even when
  the predicted gap has a strict quadratic minimum at $\sigma_0$ (for instance a predicted
  gap $\varepsilon(s-\sigma_0)^2\bigl(2+\sin(3\log\abs{s-\sigma_0})\bigr)$, which is $C^{1,1}$
  with a strict minimum, is realised by admissible leader data whose
  $\sa_\mL+\alpha v_\mL$ has negative essential values in every neighbourhood of
  $\sigma_0$).  We make no claim that nondegeneracy holds generically in any precise
  sense; \cref{T:main}(B) is conditional on it.
\end{enumerate}
\end{remark}

\paragraph{Local structure near a nondegenerate grazing state}

Fix a nondegenerate $z_0$.  Let $B$ be the closed unit ball around $z_0$, intersected
with $[0,\timehor)\times\RR^2$, and $L_B$ the constant of \cref{L:flow}(b).  Let
$\phi_{\min}>0$ and $h_1\in(0,h_0]$ be such that
$\operatorname{sgn}(\delta_0(\sigma_0))\,\delta_0\ge2\phi_{\min}$ on
$[\sigma_0-h_1,\sigma_0+h_1]$; they exist by \ref{ND2} and continuity of $\delta_0$, and
depend only on $z_0$.  Put $L_{\mathsf G}=1+2\alpha$ and let $C_\Lambda$ be a bound
for $\abs{\partial_s\Lambda(s;z)}$ over $z\in B$, $s\in[0,\timehor+1]$ (it exists by
\cref{L:flow}(b)).

\begin{lemma}[Window and sides]\label{L:window}
There are $h\in(0,h_1]$ and $\rho\in(0,1]$, with $\sigma_0+h<\timehor$ and
$\bar t_0+\rho<\sigma_0-h$, such that, writing $I\Def[\sigma_0-h,\sigma_0+h]$ and
$U\Def\{z\in\Slab:\abs{z-z_0}<\rho\}$, the following hold for every
$z=(\bar t,\bar X,\bar V)\in U$.  (Here $h$ is chosen first, in (a); all later steps only
shrink $\rho$.)
\begin{enumerate}[label=\tup{\alph*}]
\item $\partial_s\Lambda(s;z)\ge c_0$ for a.e. $s\in I$;
  $\Lambda(\sigma_0-h;z)\le-c_0h/2$ and $\Lambda(\sigma_0+h;z)\ge c_0h/2$.  Hence
  $\Lambda(\cdot;z)$ has exactly one zero $s^\ast(z)$ in $I$, it lies in the open
  interval, and $s^\ast(z)$ is the strict minimiser of $\Xi(\cdot;z)$ on $I$.  Moreover
  $\Xi(\sigma_0+h;z)\ge c_0h^2/4>0$.
\item $\Xi(s;z)\ge c>0$ for $s\in[\bar t,\sigma_0-h]$, for a constant $c$ independent of $z$.
\item The \emph{minimum-gap function} $m(z)\Def\Xi(s^\ast(z);z)=\min_I\Xi(\cdot;z)$ is
  $L_B$-Lipschitz on $U$, and $\abs{s^\ast(z_1)-s^\ast(z_2)}\le(L_B/c_0)\abs{z_1-z_2}$.
\item With $U^-\Def\{m<0\}$, $U^0\Def\{m=0\}$, $U^+\Def\{m>0\}$ (all intersected with $U$):
  \begin{itemize}
  \item $z\in U^-$ implies $z\in\Reg$, $\sigma(z)\in(\sigma_0-h,s^\ast(z))$,
    $\Lambda(\sigma(z);z)<0$ and $\tau(z)=\sigma(z)-\bar t$;
  \item $z\in U^0$ implies $z\in\Graz$ with $\sigma(z)=s^\ast(z)$; in particular
    $\Graz\cap U=U^0$;
  \item $z\in U^+$ implies $\Xi(\cdot;z)>0$ on $[\bar t,\sigma_0+h]$, and $z\in\Reg$.
  \end{itemize}
\item There is a Lipschitz function $\tau_+:U\to(0,\timehor]$ with $\tau_+=\tau$ on
  $U^0\cup U^+$.
\item $\tau\ge\tau_U>0$ on $U$; $\abs{\mathsf F(t,\xi)}\le M_U$ at \emph{every} point
  $(t,\xi)\in U$ (recall that the representative of $\sa_\mL$ is bounded everywhere), and
  $\abs{\mathsf G(s,Z(s;z))}\le M_U$ for $z\in U$, $s\in[0,\timehor+1]$.
\item $\abs{\bar V}\ge\bar v\Def\abs{\bar V_0}/2$ on $U$, and $\bar V$ has the sign of $\bar V_0$.
\end{enumerate}
\end{lemma}

\begin{proof}
(a) Along $z_0$, $Z(\sigma_0;z_0)=(0,0)$, and $Z(\cdot;z_0)$ is continuous, so there is
$h\le h_0$ with $\abs{Z(s;z_0)}\le c_0/(2L_{\mathsf G})$ on $I$.  For $\abs{z-z_0}<\rho$ and
$L_B\rho\le c_0/(2L_{\mathsf G})$, \eqref{E:flow-lip-z} gives
$\abs{Z(s;z)}\le c_0/L_{\mathsf G}$ on $I$.  Since
$\abs{\mathsf g(s,Z)-\mathsf g(s,0)}\le L_{\mathsf G}\abs Z$ and
$\mathsf g(s,0)=\sa_\mL(s)+\alpha v_\mL(s)\ge2c_0$ a.e. on $I\subset[\sigma_0-h_0,\sigma_0+h_0]$
by \ref{ND1}, we get $\partial_s\Lambda(s;z)=\mathsf g(s,Z(s;z))\ge c_0$ a.e. on $I$.  Along
$z_0$ this gives $\Lambda(\sigma_0\mp h;z_0)\lessgtr\mp c_0h$; shrinking $\rho$ so that
$L_B\rho\le c_0h/2$ transfers this to all $z\in U$ with the constant $c_0h/2$.  The
claims on $s^\ast$ follow from strict monotonicity of $\Lambda(\cdot;z)$ on $I$ and
$\partial_s\Xi=\Lambda$.  Finally
$\Xi(\sigma_0+h;z_0)=\int_{\sigma_0}^{\sigma_0+h}\Lambda(s;z_0)\,ds\ge\int_0^h c_0r\,dr=c_0h^2/2$,
and shrinking $\rho$ so that $L_B\rho\le c_0h^2/4$ gives the last claim on $U$.

Throughout, $\rho$ is also taken small enough that $U\subset\Slab$; this is possible
because $\Slab$ is relatively open and $z_0\in\Slab$.

(b) $\Xi(\cdot;z_0)$ is positive on $[\bar t_0,\sigma_0-h]$ because $\sigma_0$ is its first
zero. If $\bar t_0>0$, continuity permits a choice $h'>0$ with $h'\le\bar t_0$ such that it stays positive on $[\bar t_0-h',\sigma_0-h]$. If $\bar t_0=0$, choose any $h'>0$; the admissible interval is then $[0,\sigma_0-h]$. Thus, in both cases, it has a positive minimum $2c$ on $[\max\{0,\bar t_0-h'\},\sigma_0-h]$. Shrink $\rho\le\min\{h',c/L_B\}$.

(c) $\abs{m(z_1)-m(z_2)}\le\sup_{s\in I}\abs{\Xi(s;z_1)-\Xi(s;z_2)}\le L_B\abs{z_1-z_2}$,
as the minimum of a family of uniformly Lipschitz functions.  For $s^\ast$:
$0=\Lambda(s^\ast(z_1);z_1)=\Lambda(s^\ast(z_2);z_2)$, so by (a) and \eqref{E:flow-lip-z},
$c_0\abs{s^\ast(z_1)-s^\ast(z_2)}\le\abs{\Lambda(s^\ast(z_1);z_2)-\Lambda(s^\ast(z_2);z_2)}
=\abs{\Lambda(s^\ast(z_1);z_2)-\Lambda(s^\ast(z_1);z_1)}\le L_B\abs{z_1-z_2}$.

(d) By (a), $\Xi(\cdot;z)$ is strictly decreasing on $[\sigma_0-h,s^\ast(z)]$ from
$\Xi(\sigma_0-h;z)\ge c>0$ to $m(z)$, and strictly increasing after $s^\ast(z)$ on $I$.  By
(b) there is no zero on $[\bar t,\sigma_0-h]$.  If $m(z)<0$, the first zero $\sigma(z)$ lies
in $(\sigma_0-h,s^\ast(z))$ where $\Lambda<0$; so $z\in\Reg$ and, by \cref{L:usc},
$\tau=\sigma(z)-\bar t$.  If $m(z)=0$, the first zero is $s^\ast(z)$ with
$\Lambda(s^\ast(z);z)=0$, so $z\in\Graz$.  If $m(z)>0$, $\Xi(\cdot;z)>0$ on
$[\bar t,\sigma_0+h]$, so by \cref{L:cocycle} $\sigma(z)=\sigma(z^+)$ with
$z^+=(\sigma_0+h,Z(\sigma_0+h;z))$, and $\Lambda(\sigma(z);z)=\Lambda(\sigma(z^+);z^+)$;
since $z^+\in\Reg$ (see (e)), $z\in\Reg$.

(e) By \ref{ND3}, $z_0^+\in\Reg$.  \Cref{L:R-open} gives a neighbourhood $N$ of $z_0^+$ on
which $\tau$ is Lipschitz, with constant $L_N$.  The map $z\mapsto z^+$ is $L_B$-Lipschitz
by \eqref{E:flow-lip-z} (its time coordinate is constant), so after shrinking $\rho$,
$z^+\in N$ for all $z\in U$; also $\Xi(\sigma_0+h;z)>0$ on $U$ by (a), so $z^+\in\Om$.
Define
\[
\tau_+(z)\Def\tau(z^+)+(\sigma_0+h-\bar t),\qquad z\in U .
\]
It is Lipschitz with constant $L_NL_B+1$.  For $z\in U^0\cup U^+$, the set in
\eqref{E:TTC} contains no $s\le\sigma_0+h$: for $z\in U^+$ there is no zero there, and
for $z\in U^0$ the only zero is $s^\ast(z)$, where $\Lambda=0$.  Hence by \cref{L:cocycle}
$\Tt(z)=\Tt(z^+)$, i.e.\ $\tau(z)=\Tt(z^+)-\bar t=\tau(z^+)+(\sigma_0+h)-\bar t=\tau_+(z)$.

(f) On $U$, $\bar X\ge\bar X_0-\rho>0$ and $\bar t\le\bar t_0+\rho<\timehor$, so
\cref{L:tau-bounds} gives $\tau\ge\tau_U>0$.  The bounds on $\mathsf F$ and $\mathsf G$
follow as in \cref{L:Cara-conditions} and \cref{L:flow}(b).

(g) Shrink $\rho\le\abs{\bar V_0}/2$.
\end{proof}

\begin{lemma}[Uniform sign of the sensitivity]\label{L:uniform-phi}
In the setting of \cref{L:window}, after shrinking $\rho$ (with $h$ and $\phi_{\min}$
fixed as there), for all $z_a,z_b\in U$, all $t\in[\max\{0,\bar t_0-\rho\},\bar t_0+\rho]$ and
all $\sigma\in I$, the fundamental solution $\phi_2$ of \eqref{E:D-ODE} (with
$\tilde c=\tilde c_{ab}$, initial time $t$) satisfies
\[
\operatorname{sgn}(\delta_0(\sigma_0))\,\phi_2(\sigma)\ \ge\ \phi_{\min},
\qquad\text{and}\qquad \abs{\phi_1(\sigma)}\le C_1\Def e^{(1+2\alpha)(\timehor+1)} .
\]
\end{lemma}

\begin{proof}
The bound on $\phi_1$ is \cref{L:fundamental}(a).  By the choice of $h\le h_1$ in
\cref{L:window}, $\operatorname{sgn}(\delta_0(\sigma_0))\,\delta_0\ge2\phi_{\min}$ on $I$.
Here $\delta_0$ denotes the solution of the linear equation \eqref{E:variational} on all
of $[0,\timehor+1]$ (it is defined backward from $\bar t_0$ as well, since the equation is
linear with bounded coefficients).

Let $e\Def\phi_2-\delta_0$ on $[t,\timehor+1]$.  Then
\[
e''+\alpha e'+\alpha\tilde c_{ab}\,e=-\alpha\bigl(\tilde c_{ab}-\mathcal V'(\Xi(\cdot;z_0))\bigr)\delta_0,
\qquad e(t)=-\delta_0(t),\quad e'(t)=1-\delta_0'(t).
\]
By \eqref{E:flow-lip-z}, $\abs{\theta\Xi(s;z_a)+(1-\theta)\Xi(s;z_b)-\Xi(s;z_0)}\le L_B\rho$,
so $\abs{\tilde c_{ab}-\mathcal V'(\Xi(\cdot;z_0))}\le\norm{\mathcal V''}_\infty L_B\rho$
uniformly in $s$.  Since $\delta_0(\bar t_0)=0$, $\delta_0'(\bar t_0)=1$ and $\delta_0\in C^1$,
the initial data of $e$ are bounded by $\omega(\rho)\to0$ as $\rho\to0$.  Gr\"onwall for the
first-order system $(e,e')$ on $[t,\timehor+1]$ gives
$\sup_I\abs e\le C(\omega(\rho)+\rho)$ with $C$ depending only on $\alpha,\timehor$,
$\norm{\mathcal V''}_\infty$, $L_B$ and $\sup\abs{\delta_0}$.  Shrink $\rho$ so that this is
$\le\phi_{\min}$.
\end{proof}

\paragraph{The crossing lemma}

\begin{lemma}[Crossing]\label{L:crossing}
Let $z_0$ be nondegenerate and let $h,\rho,I,U,m,\varsigma$ be as above.  There are
$h^\ast>0$ and $c_3>0$, depending only on $z_0$ and the constants of
\cref{L:window,L:uniform-phi}, with the following property.  Let $\xi$ be a Filippov
solution on an interval $J$ with $z(t)\Def(t,\xi(t))\in U$ for all $t\in J$.  Then for all
$t_1<t$ in $J$ with $t-t_1\le h^\ast$,
\begin{equation}\label{E:m-monotone}
\varsigma\,\bigl[m(z(t))-m(z(t_1))\bigr]\ \ge\ c_3\,(t-t_1).
\end{equation}
Consequently $\varsigma\,m\circ z$ is strictly increasing on $J$ with slope at least $c_3$,
and $z(t)\in\Graz$ for at most one $t\in J$.
\end{lemma}

\begin{proof}
We treat $\bar V_0<0$; the case $\bar V_0>0$ is identical with all inequalities reversed.
Fix $t_1<t$ in $J$, write $\Delta\Def t-t_1$, and let
\[
D(s)\Def\Xi(s;z(t))-\Xi(s;z(t_1)),\qquad s\in[0,\timehor+1].
\]
By \cref{L:difference-ODE}, $D$ solves \eqref{E:D-ODE} with $\tilde c=\tilde c_{z(t),z(t_1)}$.

\emph{Step 1: data of $D$ at $s=t$.}  By definition
$Z(t;z(t))=\xi(t)$, so $D(t)=X(t)-\Xi(t;z(t_1))$ and $D'(t)=V(t)-\Lambda(t;z(t_1))$.
Since $\xi$ and $Z(\cdot;z(t_1))$ both start from $\xi(t_1)$ at time $t_1$ and have speeds
bounded by $M_U$ (\cref{L:window}(f); for the Filippov solution,
$\dot\xi\in\Fil[\mathsf F]\subset\overline{B}(0,M_U)$ a.e.),
\begin{equation}\label{E:traj-close}
\begin{aligned}
\abs{\xi(r)-Z(r;z(t_1))}&\le2M_U(r-t_1)\quad(r\in[t_1,t]),\\
\abs{D(t)}&\le\int_{t_1}^t\abs{V-\Lambda(\cdot;z(t_1))}\,dr\le M_U\Delta^2 .
\end{aligned}
\end{equation}
For $D'(t)$, write $V(t)=V(t_1)+\int_{t_1}^t\dot V\,dr$ and
$\Lambda(t;z(t_1))=V(t_1)+\int_{t_1}^t\mathsf g(r,Z(r;z(t_1)))\,dr$.  On $U$, $V<0$
(\cref{L:window}(g)), and $\tau\le\timehor$ (\cref{L:tau-bounds}), so at every point of
$U$ the TTC term satisfies $-\beta V/\tau^2\ge\beta\abs V/\timehor^2\ge\beta\bar v/\timehor^2$.
The function $\eta\mapsto\mathsf g(r,\eta)+\beta\abs{\eta_2}/\timehor^2$ is continuous and
bounds $\mathsf f(r,\cdot)$ from below at every point of $U_r\Def\{\eta:(r,\eta)\in U\}$, so
by \cref{L:Filippov-facts}(ii) (or directly for Carath\'eodory solutions)
\[
\dot V(r)\ \ge\ \mathsf g(r,\xi(r))+\beta\,\frac{\bar v}{\timehor^2}\qquad\text{for a.e. }r\in[t_1,t].
\]
Therefore, using $\abs{\mathsf g(r,\eta_1)-\mathsf g(r,\eta_2)}\le L_{\mathsf G}\abs{\eta_1-\eta_2}$
and \eqref{E:traj-close},
\begin{equation}\label{E:Dprime}
\begin{aligned}
D'(t)&\ge\int_{t_1}^t\bigl[\mathsf g(r,\xi(r))-\mathsf g(r,Z(r;z(t_1)))\bigr]dr
       +\frac{\beta\bar v}{\timehor^2}\,\Delta\\
&\ge\frac{\beta\bar v}{\timehor^2}\,\Delta-L_{\mathsf G}M_U\Delta^2 .
\end{aligned}
\end{equation}

\emph{Step 2: $D$ on the window.}  Set
$c_3\Def\beta\bar v\phi_{\min}/(2\timehor^2)$ and
$h^\ast\Def c_3/\bigl((L_{\mathsf G}\phi_{\min}+C_1)M_U\bigr)$, and assume $\Delta\le h^\ast$.
Then $L_{\mathsf G}M_U\Delta\le c_3/\phi_{\min}=\beta\bar v/(2\timehor^2)$, so
\eqref{E:Dprime} gives $D'(t)\ge\beta\bar v\Delta/(2\timehor^2)>0$.

Since $t\le\bar t_0+\rho<\sigma_0-h$, every $\sigma\in I$ satisfies $\sigma>t$, and
$D(\sigma)=D(t)\phi_1(\sigma)+D'(t)\phi_2(\sigma)$, where $\phi_1,\phi_2$ are the
fundamental solutions of \eqref{E:D-ODE} with initial time $t$ (\cref{L:fundamental}).  By
\cref{L:uniform-phi}, $\varsigma\phi_2(\sigma)\ge\phi_{\min}$ and $\abs{\phi_1(\sigma)}\le C_1$,
where $\varsigma=\operatorname{sgn}(\delta_0(\sigma_0))$ because $\bar V_0<0$.  Using
$D'(t)>0$, \eqref{E:traj-close} and \eqref{E:Dprime}, for all $\sigma\in I$,
\begin{align*}
\varsigma D(\sigma)\ &\ge\ D'(t)\,\phi_{\min}-\abs{D(t)}\,C_1
\ \ge\ \Bigl(\frac{\beta\bar v}{\timehor^2}\Delta-L_{\mathsf G}M_U\Delta^2\Bigr)\phi_{\min}-C_1M_U\Delta^2\\
&=\ 2c_3\Delta-(L_{\mathsf G}\phi_{\min}+C_1)M_U\Delta^2\ \ge\ c_3\,\Delta .
\end{align*}

\emph{Step 3: from $D$ to $m$.}  Suppose $\varsigma=+1$.  Then
$\Xi(\sigma;z(t))\ge\Xi(\sigma;z(t_1))+c_3\Delta$ for all $\sigma\in I$.  Taking the
minimum over $I$ on both sides, $m(z(t))\ge m(z(t_1))+c_3\Delta$.  If $\varsigma=-1$, the
same argument gives $m(z(t))\le m(z(t_1))-c_3\Delta$.  This is \eqref{E:m-monotone}.

\emph{Step 4: consequences.}  Chaining \eqref{E:m-monotone} over a partition of
$[t_1,t]$ into pieces of length $\le h^\ast$ shows that \eqref{E:m-monotone} holds for all
$t_1<t$ in $J$.  By \cref{L:window}(d), $z(t)\in\Graz$ iff $m(z(t))=0$, and a strictly
monotone function vanishes at most once.
\end{proof}

\begin{remark}[What the proof uses]\label{R:crossing-uses}
Only three things enter: the cocycle property, which makes the prediction from $z(t)$ and
from $z(t_1)$ differ only through what the closed loop did on $[t_1,t]$; the fact that this
difference is $-\beta\int V/\tau^2$ to leading order, whose sign is that of $-V$ on
\emph{both} sides of $\Graz$ (only its size jumps); and the sign of the sensitivity
$\delta_0(\sigma_0)$.  The value of $\tau$ \emph{on} $\Graz$ plays no role, and no Lipschitz
property of $\tau$ across $\Graz$ is used.
\end{remark}

\begin{lemma}[Carath\'eodory solutions are Filippov solutions near a nondegenerate graze]\label{L:Cara-is-Filippov}
Let $z_0$ be nondegenerate and $U$ as in \cref{L:window,L:uniform-phi}.  Then
$\mathsf F(t,\xi)\in\Fil[\mathsf F](t,\xi)$ for every $(t,\xi)\in U$.  Consequently every
Carath\'eodory solution of \eqref{E:closed-loop} whose graph lies in $U$ is a Filippov
solution.
\end{lemma}
\begin{proof}
On $U^-\cup U^+\subset\Reg$ the claim is \cref{L:Filippov-facts}(i).  Let
$(t,\xi)=(t,X,V)\in U^0$.  Fix $\varepsilon\neq0$ small and put $z_\varepsilon\Def(t,X,V+\varepsilon)\in U$.
By \cref{L:difference-ODE} with $z_a=z_\varepsilon$, $z_b=(t,\xi)$, the difference
$D=\Xi(\cdot;z_\varepsilon)-\Xi(\cdot;t,\xi)$ solves \eqref{E:D-ODE} with $D(t)=0$,
$D'(t)=\varepsilon$, so $D(\sigma)=\varepsilon\phi_2(\sigma)$ with $\phi_2$ as in
\cref{L:uniform-phi} (initial time $t\in[\max\{0,\bar t_0-\rho\},\bar t_0+\rho]$), so that
$\varsigma_0\phi_2\ge\phi_{\min}$ on $I$ with $\varsigma_0\Def\operatorname{sgn}(\delta_0(\sigma_0))$.
Take $\varepsilon$ with $\varsigma_0\varepsilon>0$, i.e.\ $\varepsilon=\varsigma_0\abs\varepsilon$.  Then
$D=\abs\varepsilon\,\varsigma_0\phi_2\ge\abs\varepsilon\phi_{\min}>0$ on $I$ (for either sign of
$\varsigma_0$), and taking the minimum over $I$ as in Step~3 of \cref{L:crossing},
$m(z_\varepsilon)-m(t,\xi)\ge\abs\varepsilon\phi_{\min}>0$.  Since $m(t,\xi)=0$,
this gives $m(z_\varepsilon)>0$ whenever $\varsigma_0\varepsilon>0$: the point $(t,\xi)$ is a limit
of points of $U^+_t\Def\{\eta:(t,\eta)\in U^+\}$, which is open.

On $U^0_t\cup U^+_t$ we have $\tau(t,\cdot)=\tau_+(t,\cdot)$ (\cref{L:window}(e)), which is
continuous, so $\mathsf F(t,\cdot)$ restricted to $U^0_t\cup U^+_t$ is continuous.  Given
$\delta>0$ and a null set $N\subset\RR^2$, pick $\eta_n\to\xi$ with $\eta_n\in U^+_t$; since
$U^+_t$ is open we may move each $\eta_n$ off $N$ without leaving $B(\xi,\delta)\cap U^+_t$
and keeping $\eta_n\to\xi$.  Then $\mathsf F(t,\eta_n)\to\mathsf F(t,\xi)$, so
$\mathsf F(t,\xi)\in\overline{\mathsf F(t,B(\xi,\delta)\setminus N)}\subset\co\,\mathsf F(t,B(\xi,\delta)\setminus N)$.
Intersecting over $\delta$ and $N$ gives $\mathsf F(t,\xi)\in\Fil[\mathsf F](t,\xi)$.

For the last statement: a Carath\'eodory solution satisfies $\dot\xi=\mathsf F(t,\xi)$ a.e.,
and at those times $\mathsf F(t,\xi)\in\Fil[\mathsf F](t,\xi)$.
\end{proof}

\subsection*{Uniqueness through a nondegenerate crossing}

\begin{proposition}[Uniqueness through the crossing]\label{P:unique-crossing}
Let $z_1=(t_1,\xi_1)\in\Graz$ be nondegenerate, and let $\xi^{a},\xi^{b}$ be Filippov
solutions on $[t_1,t_1+\eta)$ with $\xi^a(t_1)=\xi^b(t_1)=\xi_1$.  Then there is
$\eta'\in(0,\eta]$ such that $\xi^a=\xi^b$ on $[t_1,t_1+\eta')$, the graph of $\xi^a$ lies
in $\Reg$ for $t\in(t_1,t_1+\eta')$, and $\xi^a$ is a Carath\'eodory solution of
\eqref{E:closed-loop} on $[t_1,t_1+\eta')$.
\end{proposition}

\begin{proof}
Apply \cref{L:window,L:uniform-phi,L:crossing} with $z_0=z_1$, and shrink $\eta'$ so that
both graphs stay in $U$ on $[t_1,t_1+\eta')$ and $\eta'\le h^\ast$.  By
\cref{L:crossing}, for $t\in(t_1,t_1+\eta')$ and $i\in\{a,b\}$,
\begin{equation}\label{E:m-away}
\varsigma\,m(z^i(t))\ \ge\ c_3(t-t_1)\ >\ 0 .
\end{equation}
So both graphs lie in $U^+$ (if $\varsigma=+1$) or in $U^-$ (if $\varsigma=-1$) for
$t>t_1$; in either case in $\Reg$ (\cref{L:window}(d)).  On $\Reg$, $\mathsf F(t,\cdot)$ is
continuous, so $\Fil[\mathsf F]=\{\mathsf F\}$ (\cref{L:Filippov-facts}(i)) and both $\xi^i$
are Carath\'eodory solutions of \eqref{E:closed-loop} on $(t_1,t_1+\eta')$, hence on
$[t_1,t_1+\eta')$ (a single point does not affect the a.e.\ identity).  It remains to show
$\xi^a=\xi^b$.

\emph{Case $\varsigma=+1$ (entering the miss side).}  For $t\in[t_1,t_1+\eta')$,
$z^i(t)\in U^0\cup U^+$, so $\tau(z^i(t))=\tau_+(z^i(t))$ with $\tau_+$ Lipschitz on $U$
(\cref{L:window}(e)).  Hence both $\xi^i$ solve
$\dot\xi=\mathsf F_+(t,\xi)\Def(V,\ \mathsf g(t,\xi)-\beta V/\tau_+(t,\xi)^2)$ on
$[t_1,t_1+\eta')$.  Since $\tau_+(z)\ge\sigma_0+h-\bar t>2h>0$ on $U$ (as
$\bar t<\sigma_0-h$) and $\tau_+$ is Lipschitz, $\mathsf F_+$ is
Lipschitz in $\xi$ on $U$ (as in \cref{L:Cara-conditions}), and Gr\"onwall gives
$\xi^a=\xi^b$.

\emph{Case $\varsigma=-1$ (entering the crossing side).}  Here $\tau$ is not Lipschitz up
to $\Graz$ (\cref{R:discontinuity}), and we use a Lipschitz bound with an integrable
singularity.  Fix $t\in(t_1,t_1+\eta')$ and write $\mu_i\Def-m(z^i(t))\ge c_3(t-t_1)$ and
$\sigma_i\Def\sigma(z^i(t))$, $s^\ast_i\Def s^\ast(z^i(t))$, so that
$\tau(z^i(t))=\sigma_i-t$ by \cref{L:window}(d).

\emph{\tup{i} The first zero is at distance $\gtrsim\sqrt{\mu_i}$ from $s^\ast_i$.}  On
$[\sigma_i,s^\ast_i]$, $\Lambda(\cdot;z^i(t))\le0$ and, since
$\partial_s\Lambda\le C_\Lambda$ and $\Lambda(s^\ast_i;z^i(t))=0$,
$\abs{\Lambda(s;z^i(t))}\le C_\Lambda(s^\ast_i-s)$.  Therefore
\[
\mu_i=\Xi(\sigma_i;z^i(t))-\Xi(s^\ast_i;z^i(t))=\int_{\sigma_i}^{s^\ast_i}\abs{\Lambda(s;z^i(t))}\,ds
\le\frac{C_\Lambda}{2}(s^\ast_i-\sigma_i)^2,
\]
that is, $s^\ast_i-\sigma_i\ge\sqrt{2\mu_i/C_\Lambda}\Def2\delta_i$.

\emph{\tup{ii} Lipschitz bound on $\tau$ between the two states.}  Let
\[
\begin{aligned}
\delta&\Def\min\{\delta_a,\delta_b\}\ge\tfrac12\sqrt{2c_3(t-t_1)/C_\Lambda},\\
d&\Def\abs{\xi^a(t)-\xi^b(t)}=\abs{z^a(t)-z^b(t)}.
\end{aligned}
\]
By \cref{L:window}(c),
$\abs{s^\ast_a-s^\ast_b}\le(L_B/c_0)d$.  Suppose $(L_B/c_0)d\le\delta$.  Then
\[
\begin{aligned}
\sigma_b&\le s^\ast_b-2\delta\le s^\ast_a+\delta-2\delta=s^\ast_a-\delta,\\
\sigma_a&\le s^\ast_a-2\delta\le s^\ast_a-\delta.
\end{aligned}
\]
Both $\sigma_a,\sigma_b$ lie in
$[\sigma_0-h,\ s^\ast_a-\delta]$, where $\Lambda(\cdot;z^a(t))\le-c_0\delta$ by
\cref{L:window}(a) (strict increase with slope $\ge c_0$ up to its zero $s^\ast_a$).
Since $\Xi(\sigma_a;z^a(t))=0=\Xi(\sigma_b;z^b(t))$,
\begin{align*}
c_0\delta\,\abs{\sigma_a-\sigma_b}\ &\le\ \Bigl|\int_{\sigma_b}^{\sigma_a}\Lambda(s;z^a(t))\,ds\Bigr|
=\abs{\Xi(\sigma_a;z^a(t))-\Xi(\sigma_b;z^a(t))}\\
&=\abs{\Xi(\sigma_b;z^b(t))-\Xi(\sigma_b;z^a(t))}\ \le\ L_B\,d ,
\end{align*}
so $\abs{\tau(z^a(t))-\tau(z^b(t))}=\abs{\sigma_a-\sigma_b}\le L_Bd/(c_0\delta)$.

\emph{\tup{iii} The smallness condition holds for $t$ close to $t_1$.}  Both solutions have
speed $\le M_U$ and coincide at $t_1$, so $d\le2M_U(t-t_1)$.  The condition
$(L_B/c_0)d\le\delta$ is implied by
$2L_BM_U(t-t_1)/c_0\le\tfrac12\sqrt{2c_3(t-t_1)/C_\Lambda}$, i.e.\ by
$\sqrt{t-t_1}\le c_0\sqrt{2c_3/C_\Lambda}/(4L_BM_U)$.  Shrink $\eta'$ accordingly.  Then,
for all $t\in(t_1,t_1+\eta')$,
\begin{equation}\label{E:tau-sing-lip}
\begin{aligned}
\abs{\tau(z^a(t))-\tau(z^b(t))}&\le\frac{C_\tau}{\sqrt{t-t_1}}\,\abs{\xi^a(t)-\xi^b(t)},\\
C_\tau&\Def\frac{2L_B}{c_0}\sqrt{\frac{C_\Lambda}{2c_3}} .
\end{aligned}
\end{equation}

\emph{\tup{iv} Gr\"onwall with an integrable coefficient.}  With $\tau\ge\tau_U$ and
$\abs V\le{v_{\max}}+1$ on $U\subset\Slab$, the same computation as in \cref{L:Cara-conditions} gives, for
a.e. $t\in(t_1,t_1+\eta')$,
\[
\begin{gathered}
\abs{\mathsf F(t,\xi^a(t))-\mathsf F(t,\xi^b(t))}
\le L(t)\,\abs{\xi^a(t)-\xi^b(t)},\\
L(t)\Def 1+2\alpha+\frac{\beta}{\tau_U^2}
 +\frac{2\beta({v_{\max}}+1)\timehor}{\tau_U^4}\,\frac{C_\tau}{\sqrt{t-t_1}} .
\end{gathered}
\]
$L\in L^1(t_1,t_1+\eta')$.  Let $u\Def\abs{\xi^a-\xi^b}$; then $u(t_1)=0$ and
$u(t)\le\int_{t_1}^tL(s)u(s)\,ds$.  Iterating, $u(t)\le\sup u\cdot\bigl(\int_{t_1}^tL\bigr)^n/n!$
for every $n$, so $u\equiv0$.
\end{proof}

\begin{remark}
The exponent $\tfrac12$ in \eqref{E:tau-sing-lip} is the H\"older exponent, at $U^0$, of
the crossing-side branch $\tau^-(z)\Def\sigma(z)-\bar t$, which equals $\tau$ on $U^-$ and
extends continuously to $U^0$ by $s^\ast(z)-\bar t$ (\cref{L:window}(d)); $\tau$ itself
takes the value $\tau_+$ on $U^0$ and is discontinuous there (\cref{R:discontinuity}).
The first zero of a gap $\approx a(s-\sigma_0)^2-\mu$ moves like $\sqrt\mu$.  What rescues uniqueness is that the
crossing is \emph{transversal} with a definite speed, $\mu\ge c_3(t-t_1)$, so the singular
Lipschitz coefficient is $O((t-t_1)^{-1/2})$ and integrable.  This is the one place where
the quantitative form \eqref{E:m-monotone} of the crossing lemma is needed.
\end{remark}

\subsection{Main theorem}\label{SS:main}

\begin{theorem}[Well-posedness and safety of the closed loop]\label{T:main}\label{T:collision_case1}\label{T:existence}
Let $(t_\circ,\xi_\circ)\in\Phys$.
\begin{enumerate}[label=\tup{\alph*}]
\item \textbf{(Unconditional.)}  Non-extendable Filippov solutions of
  \eqref{E:closed-loop} with initial value $\xi_\circ$ at $t_\circ$ exist.  \emph{Every}
  non-extendable Filippov solution $\xi$, on its maximal interval $[t_\circ,T^\ast)$,
  $T^\ast=T^\ast_\xi\in(t_\circ,\timehor]$, satisfies:
  \begin{enumerate}[label=\tup{\roman*}]
  \item $(t,\xi(t))\in\Phys$ for all $t$; in particular $0\le v\le{v_{\max}}$ and $\abs V\le{v_{\max}}$;
  \item either $T^\ast=\timehor$, or $X(T^\ast-)=0$; in the latter case $V(t)\to0$ as
    $t\uparrow T^\ast$, and $\abs V\le\sqrt{2\abar X}$ on a left neighbourhood of $T^\ast$;
  \item if $T^\ast=\timehor$, then $\xi$ extends continuously to $[t_\circ,\timehor]$ with
    $V(\timehor)=0$.
  \end{enumerate}
  In particular no impact (contact with negative relative velocity) ever occurs, for any
  Filippov solution whatsoever (\cref{T:no-impact} is proved for an arbitrary solution).
\item \textbf{(Uniqueness under nondegenerate grazing.)}  Let $\xi$ be a non-extendable Filippov solution as in (A), and
  suppose that every grazing state on its graph, i.e.\ every $(t,\xi(t))\in\Graz$, is
  nondegenerate in the sense of \cref{D:nondegenerate}.  Then:
  \begin{enumerate}[label=\tup{\roman*}]
  \item the set of times $t$ with $(t,\xi(t))\in\Graz$ is discrete in $[t_\circ,T^\ast)$;
  \item $\xi$ is a Carath\'eodory solution of \eqref{E:closed-loop} on $[t_\circ,T^\ast)$;
  \item $\xi$ is the unique Filippov, hence the unique Carath\'eodory, solution with
    initial value $\xi_\circ$ at $t_\circ$, on $[t_\circ,T^\ast)$.
  \end{enumerate}
\end{enumerate}
\end{theorem}

\begin{proof}
(A) is \cref{T:global} together with \cref{T:no-impact,P:barrier,P:horizon}.

(B)(i).  Write $z(t)\Def(t,\xi(t))$.  Let $t_1$ be a grazing time.  By hypothesis
$z(t_1)$ is nondegenerate.
\Cref{L:window,L:crossing} applied at $z_0=z(t_1)$ give a neighbourhood $U$ and, by
continuity of $\xi$, an interval $(t_1-\eta,t_1+\eta)$ on which the graph stays in $U$.
On that interval $\varsigma\,m\circ z$ is strictly increasing and vanishes at $t_1$, so
$z(t)\notin\Graz$ for $t\neq t_1$ there.  Thus grazing times are isolated.  They cannot
accumulate at a point $t^\dagger<T^\ast$ either: if $z(t^\dagger)\in\Reg$, then, $\Reg$
being open, $z(t)\in\Reg$ near $t^\dagger$; if $z(t^\dagger)\in\Graz$, it is nondegenerate
and isolated by what was just shown.

(B)(ii).  Off the grazing times the graph lies in $\Reg$, where
$\Fil[\mathsf F]=\{\mathsf F\}$, so $\dot\xi=\mathsf F(t,\xi)$ a.e.

(B)(iii).  \emph{Filippov competitors.}  Let $\tilde\xi$ be a Filippov solution on
$J\ni t_\circ$ with $\tilde\xi(t_\circ)=\xi_\circ$, and let
$S\Def\sup\{s\in J\cap[t_\circ,T^\ast):\ \tilde\xi=\xi\text{ on }[t_\circ,s]\}$.  Suppose
$S<\min\{T^\ast,\sup J\}$.  By continuity $\tilde\xi(S)=\xi(S)$.  If $(S,\xi(S))\in\Reg$,
\cref{T:local-WP} (applied from time $S$) gives $\tilde\xi=\xi$ on a right neighbourhood of
$S$.  If $(S,\xi(S))\in\Graz$, it is nondegenerate by hypothesis, and
\cref{P:unique-crossing} gives the same.  Either way the definition of $S$ is
contradicted.  Hence $\tilde\xi=\xi$ on $J\cap[t_\circ,T^\ast)$.  Finally $\sup J\le T^\ast$:
otherwise $\tilde\xi$ would extend $\xi$ beyond $T^\ast$, contradicting maximality.

\emph{Carath\'eodory competitors.}  Let $\xi_C$ be a Carath\'eodory solution on $J$ with
$\xi_C(t_\circ)=\xi_\circ$, and define $S$ as above with $\xi_C$ in place of $\tilde\xi$.
Suppose $S<\min\{T^\ast,\sup J\}$.  If $(S,\xi(S))\in\Reg$, then on a right neighbourhood
of $S$ the graph of $\xi_C$ lies in $\Reg$, where $\Fil[\mathsf F]=\{\mathsf F\}$, so
$\xi_C$ is a Filippov solution there, and \cref{T:local-WP} applies.  If
$(S,\xi(S))\in\Graz$, it is nondegenerate; on a right neighbourhood of $S$ the graph of
$\xi_C$ lies in the set $U$ of \cref{L:window} for $z_0=(S,\xi(S))$, so by
\cref{L:Cara-is-Filippov} $\xi_C$ is a Filippov solution there, and
\cref{P:unique-crossing} applies.  Either way $S$ is contradicted, so $\xi_C=\xi$ on
$J\cap[t_\circ,T^\ast)$.
\end{proof}

\begin{remark}[Why both layers]\label{R:both-layers}
Layer (A) is the physically relevant statement and holds without any nondegeneracy
assumption: whatever the leader does, no admissible motion of the closed loop ends in an
impact.  Layer (B) says that the generalised solution concept is only a safety net:
whenever every grazing state met is nondegenerate, the solution is a Carath\'eodory
solution and is unique.  It is in general not $C^1$ at a crossing time $t_1$, even when
$\sa_\mL$ is continuous there: $\tau$ itself jumps along the trajectory (between the miss-side
value $\tau_+$ and the crossing-side value $s^\ast-t_1$, in either temporal orientation), so at a genuine TTC jump with
$V(t_1)\neq0$ (\ref{ND4}) the one-sided accelerations, when they exist, differ by
$-\beta V(t_1)\bigl(\tau_{\rm after}^{-2}-\tau_{\rm before}^{-2}\bigr)\neq0$.  Under
$L^\infty$ forcing one-sided accelerations need not exist at all, and if $\sa_\mL$ has
one-sided limits at $t_1$ their difference $\sa_\mL(t_1+)-\sa_\mL(t_1-)$ adds to this
expression and may cancel it; so non-$C^1$ behaviour is not automatic from \ref{ND4}
alone, only for forcing continuous at $t_1$.  A Carath\'eodory solution that crosses
$\Graz$ transversally is automatically a Filippov solution, so (B) is a statement about the
same object as (A), not about a different solution concept.
\end{remark}

\begin{corollary}[Impact prevention]\label{C:no-collision}
Let $\xi$ be a non-extendable Filippov solution as in \cref{T:main}(A).
\begin{enumerate}[label=\tup{\alph*}]
\item There is no impact.  If $X(T^\ast-)=0$, then $V(t)\to0$ as $t\uparrow T^\ast$.
  Since $X>0$ on $[t_\circ,T^\ast)$, this is the only time at which a contact could
  occur.  In particular, no collision, in the sense of a contact with negative relative
  velocity, occurs.
\item If $X(T^\ast-)>0$, then $T^\ast=\timehor$, and the solution, defined on
  $[t_\circ,\timehor)$, extends continuously to $[t_\circ,\timehor]$ with
  $V(\timehor)=0$ (\cref{P:horizon}).
\end{enumerate}
Thus the only obstruction to existence on all of $[t_\circ,\timehor]$ is a contact with
zero closing speed; see \cref{R:open-closedloop}.
\end{corollary}

\begin{remark}[What is not proved]\label{R:open-closedloop}
\begin{enumerate}[label=\tup{\alph*}]
\item \emph{Contact with zero closing speed.}  \cref{T:no-impact} excludes every contact
  with negative relative velocity, which is the notion of collision used in this paper.  It
  does not exclude a contact with $X\to0$ and $V\to0$ (a touch without impact). Exclusion of such finite-time touches remains open; the numerical experiments below are not a proof of strict gap positivity.
\item \emph{Degenerate grazing states.}  \cref{T:main}(B) gives uniqueness under the
  nondegeneracy conditions of \cref{D:nondegenerate}.  The excluded cases
  (\cref{R:ND-meaning}(e)) are each one further condition on top of $z_0\in\Graz$; Filippov
  solutions still exist through them, so safety is unaffected.
\item \emph{The horizon convention.}  With $\inf\varnothing=\timehor$ in \eqref{E:TTC}, the
  gain $\beta/\tau^2\ge\beta/(\timehor-t)^2$ is unbounded as $t\uparrow\timehor$, which is
  what forces $V(\timehor)=0$ in \cref{P:horizon}.  Two alternatives avoid this: a sliding
  prediction horizon $H$ (replace $\timehor$ by $\bar t+H$ in \eqref{E:TTC}, with the leader
  data prescribed on $[0,\timehor+H]$ with the same acceleration and speed bounds), or $\inf\varnothing=+\infty$ with
  $1/\tau^2\Def0$ and the infimum taken over $[\bar t,\infty)$. The latter convention requires leader data on all of $[0,\infty)$, with uniformly bounded acceleration and speed in $[0,v_{\max}]$. Under these corresponding hypotheses, \cref{L:approach,T:no-impact}
  hold unchanged under either; with the second, the case split in the proof of
  \cref{L:approach} disappears.  (The hybrid convention $\inf\varnothing=+\infty$ with the
  window $[\bar t,\timehor]$ is \emph{not} admissible: it creates a new discontinuity of
  $1/\tau^2$ inside $\Reg$ along which a Filippov solution can slide into an impact at
  $\timehor$.)
\end{enumerate}
\end{remark}

\section{Regularity of the TTC map}\label{S:NN_tau_analytics}\label{sec:ttc-regularity}
We now study the TTC map $\tau$ of \cref{D:TTC} as a function of the initial state, on the domain of states whose first predicted contact is transversal. Everything in this section concerns the prediction flow alone; the closed loop does not enter. We write $z=(s,\bm y)$ with $\bm y=(\Xi,\Lambda)\in\RR^2$ for a state, $Z(\kappa;s,\bm y)=(\Xi,\Lambda)(\kappa;s,\bm y)$ for the predicted state at \emph{elapsed} time $\kappa$, and $Z^{\rm abs}(t;z)$ when the absolute time $t$ is meant (so $Z(\kappa;s,\bm y)=Z^{\rm abs}(s+\kappa;z)$; in \cref{sec: analytical_results} only the absolute-time notation was used). The prediction flow is global: by \cref{L:flow}(a) it is defined on $[0,\timehor+1]$ for all data, through and beyond the collision set, and we use this freely.

\subsection{The domain \texorpdfstring{$\mathbb D$}{D} and the optimality principle}\label{S:domainD}\label{ss:tau-analytics}

\begin{definition}[The domain $\mathbb D$]\label{D:domainD}
\begin{equation}\label{E:mathbb-D}
\mathbb D\Def\bigl\{z=(\bar t,\bar X,\bar V)\in\Om:\ \sigma(z)<\timehor,\
\Lambda(\sigma(z);z)<0\bigr\},
\end{equation}
the states whose \emph{first} predicted contact occurs before the horizon and is
transversal.  By \cref{L:usc}, $\Tt=\sigma$ and $\tau=\sigma-\bar t$ on $\mathbb D$: the
TTC \eqref{E:TTC} entering the control law \eqref{E:closed-loop} and the first zero of the
predicted gap agree on $\mathbb D$, which is the form anticipated in \eqref{E:intro-tau};
moreover $\mathbb D\subset\Reg$.
\end{definition}

\begin{lemma}[$\mathbb D$ is open; invariance along the prediction; optimality principle]\label{L:D-open}\label{T:DPP}
\begin{enumerate}[label=\tup{\alph*}]
\item $\mathbb D$ is relatively open in $[0,\timehor)\times\RR^2$.
\item (Invariance.)  If $z=(\bar t,\bar y)\in\mathbb D$ and $0\le\kappa<\tau(z)$, then
  $z^\kappa\Def(\bar t+\kappa,Z(\kappa;z))\in\mathbb D$ (writing $Z(\kappa;z)$ for the
  predicted state at elapsed time $\kappa$), and
  \begin{equation}\label{E:DPP}
  \tau(z^\kappa)=\tau(z)-\kappa .
  \end{equation}
\end{enumerate}
\end{lemma}

\begin{proof}
(a) Let $z_0\in\mathbb D$, $\sigma_0=\sigma(z_0)<\timehor$, $\Lambda(\sigma_0;z_0)=-2\gamma<0$.
In the proof of \cref{L:R-open}(b), choose $h$ additionally with $\sigma_0+h<\timehor$.
For $z$ in the ball constructed there, the first zero $\hat\sigma(z)$ of the predicted gap
lies in $(\sigma_0-h,\sigma_0+h)$, hence $<\timehor$, and $\Lambda(\hat\sigma(z);z)\le-\gamma<0$.
So $\sigma(z)=\hat\sigma(z)<\timehor$ is transversal, i.e.\ $z\in\mathbb D$.

(b) By the cocycle property (\cref{L:cocycle}) the prediction from $z^\kappa$ is the
continuation of the prediction from $z$.  Since $\kappa<\tau(z)$, the gap is positive at
$z^\kappa$ and its first zero is still $\sigma(z)$, with the same transversal velocity.
Hence $z^\kappa\in\mathbb D$ and $\tau(z^\kappa)=\sigma(z)-(\bar t+\kappa)=\tau(z)-\kappa$.
\end{proof}

\begin{remark}[Optimality principle]\label{R:DPP}
\eqref{E:DPP} is the dynamic programming identity for the value function $\bm V\Def\tau$:
along every predicted characteristic the value decreases at unit rate until the collision
set is reached.
\end{remark}

\subsection{Uniform constants on a compact set}\label{S:uniform-constants}

Throughout this section $\bm D\Subset\mathbb D$ is compact.  Let $B\supset\bm D$ be a
bounded set (for definiteness, the closed $1$-neighbourhood of $\bm D$ intersected with
$[0,\timehor)\times\RR^2$) and let $m_B,L_B$ be the constants of \cref{L:flow}(b): for
$z,z_1,z_2\in B$ and all absolute times $s,s_1,s_2\in[0,\timehor+1]$,
\begin{equation}\label{E:flow-consts}
\abs{Z(s_1;z)-Z(s_2;z)}\le m_B\abs{s_1-s_2},\qquad
\abs{Z(s;z_1)-Z(s;z_2)}\le L_B\abs{z_1-z_2}.
\end{equation}

\begin{lemma}[Comparison at a common elapsed time]\label{L:elapsed}
For $z_i=(\bar t_i,\bar y_i)\in B$ and every elapsed time $\theta\ge0$ with
$\bar t_i+\theta\in[0,\timehor+1]$ for $i=1,2$,
\begin{equation}\label{E:elapsed}
\abs{Z(\theta;z_1)-Z(\theta;z_2)}\ \le\ (L_B+m_B)\,\abs{z_1-z_2},
\end{equation}
where $Z(\theta;z_i)$ denotes the predicted state of $z_i$ at elapsed time $\theta$, i.e.\
at absolute time $\bar t_i+\theta$.
\end{lemma}

\begin{proof}
Insert the state of trajectory $1$ at the absolute time $\bar t_2+\theta$:
\[
\begin{aligned}
&\abs{Z^{\rm abs}(\bar t_1+\theta;z_1)-Z^{\rm abs}(\bar t_2+\theta;z_2)}\\
&\quad\le\abs{Z^{\rm abs}(\bar t_1+\theta;z_1)-Z^{\rm abs}(\bar t_2+\theta;z_1)}\\
&\qquad+\abs{Z^{\rm abs}(\bar t_2+\theta;z_1)-Z^{\rm abs}(\bar t_2+\theta;z_2)},
\end{aligned}
\]
and apply \eqref{E:flow-consts} to each term: $m_B\abs{\bar t_1-\bar t_2}+L_B\abs{z_1-z_2}$.
\end{proof}

\begin{lemma}[Continuity and positivity of $\tau$; uniform transversality]\label{L:unif-transv}
\begin{enumerate}[label=\tup{\alph*}]
\item $\tau$ is locally Lipschitz, hence continuous, on $\mathbb D$, and
  $\tau_{\min}(\bm D)\Def\min_{\bm D}\tau>0$,\ $\tau_{\max}(\bm D)\Def\max_{\bm D}\tau<\timehor$.
\item The \emph{terminal approach rate} $g(z)\Def-\Lambda(\tau(z);z)$ is continuous and
  positive on $\mathbb D$, with
  \begin{equation}\label{E:g-modulus}
  \begin{aligned}
  \abs{g(z_1)-g(z_2)}
  &\le m_B\abs{\tau(z_1)-\tau(z_2)}+(L_B+m_B)\abs{z_1-z_2}\\
  &\hspace{6mm}(z_1,z_2\in B\cap\mathbb D,\ \abs{\bar t_1-\bar t_2}\le1).
  \end{aligned}
  \end{equation}
  Hence the \emph{transversality margin}
  \begin{equation}\label{E:gamma-D}
  \gamma_{\bm D}\Def\min_{z\in\bm D}g(z)=\min_{z\in\bm D}\bigl(-\Lambda(\tau(z);z)\bigr)
  \end{equation}
  is attained and positive: $\Lambda(\tau(z);z)\le-\gamma_{\bm D}<0$ for all $z\in\bm D$.
\end{enumerate}
\end{lemma}

\begin{proof}
(a) $\mathbb D\subset\Reg$ and \cref{L:R-open}.  Positivity of $\tau$ is \cref{L:tau-bounds};
$\tau<\timehor-\bar t\le\timehor$ on $\mathbb D$ since $\sigma<\timehor$; the extrema are
attained by continuity and compactness.

(b) Positivity is the definition of $\mathbb D$.  For \eqref{E:g-modulus}, write
$g(z_1)-g(z_2)=-[\Lambda(\tau_1;z_1)-\Lambda(\tau_2;z_1)]-[\Lambda(\tau_2;z_1)-\Lambda(\tau_2;z_2)]$
with $\tau_i=\tau(z_i)$, and bound the first bracket by $m_B\abs{\tau_1-\tau_2}$
(elapsed-time Lipschitz, same trajectory) and the second by \cref{L:elapsed}.  Continuity
of $g$ follows from that of $\tau$, and the minimum over the compact $\bm D$ is attained.
\end{proof}

\begin{lemma}[Terminal window]\label{L:terminal-window}
Set
\begin{equation}\label{E:h0-def}
h_0=h_0(\bm D)\Def\min\Bigl\{\frac{\gamma_{\bm D}}{2m_B},\ \tau_{\min}(\bm D)\Bigr\}>0 .
\end{equation}
Then for every $z\in\bm D$:
\begin{align}
\Lambda(\theta;z)&\le-\frac{\gamma_{\bm D}}{2}\qquad\text{for all }\theta\in[\tau(z)-h_0,\tau(z)],\label{E:terminal-approach}\\
\Xi(\tau(z)-h;z)&\ge\frac{\gamma_{\bm D}}{2}\,h\qquad\text{for all }h\in[0,h_0].\label{E:fd-transversality}
\end{align}
\end{lemma}

\begin{proof}
$\Lambda(\cdot;z)$ is $m_B$-Lipschitz in elapsed time and equals $-g(z)\le-\gamma_{\bm D}$
at $\tau(z)$; so on the window it is $\le-\gamma_{\bm D}+m_Bh_0\le-\gamma_{\bm D}/2$.  The
second term in \eqref{E:h0-def} keeps the window inside $[0,\tau(z)]$.  Integrating
$\partial_\theta\Xi=\Lambda$ from $\tau(z)-h$ to $\tau(z)$, where $\Xi=0$, gives
\eqref{E:fd-transversality}.
\end{proof}

\subsection{Lipschitz regularity of \texorpdfstring{$\tau$}{tau} with explicit constant}\label{S:tau-lipschitz}

\begin{proposition}[Finite-difference Lipschitz bound]\label{P:tau-fd-lip}
There is $\varrho_0=\varrho_0(\bm D)>0$ such that for all $z_1,z_2\in\bm D$ with
$\abs{z_1-z_2}<\varrho_0$,
\begin{equation}\label{E:tau-fd-lip}
\abs{\tau(z_1)-\tau(z_2)}\ \le\ L_\tau(\bm D)\,\abs{z_1-z_2},\qquad
L_\tau(\bm D)\Def\frac{2\,(L_B+m_B)}{\gamma_{\bm D}} .
\end{equation}
\end{proposition}

\begin{proof}
By \cref{L:unif-transv}(a), $\tau$ is uniformly continuous on $\bm D$; choose $\varrho_0$
so that $\abs{z_1-z_2}<\varrho_0$ implies $\abs{\tau(z_1)-\tau(z_2)}\le h_0$.  Let
$\nu_i=\tau(z_i)$, and without loss of generality $\nu_1\le\nu_2$, $h\Def\nu_2-\nu_1\in[0,h_0]$.

\emph{Step 1.}  By \eqref{E:fd-transversality} applied to $z_2$ at $h$,
$\Xi(\nu_1;z_2)=\Xi(\nu_2-h;z_2)\ge\tfrac{\gamma_{\bm D}}{2}h$.

\emph{Step 2.}  By \cref{L:elapsed} at the common elapsed time $\nu_1$ (both
$\bar t_i+\nu_1\le\timehor$), and $\Xi(\nu_1;z_1)=0$,
\[
\Xi(\nu_1;z_2)=\abs{\Xi(\nu_1;z_2)-\Xi(\nu_1;z_1)}\le(L_B+m_B)\abs{z_1-z_2}.
\]

\emph{Step 3.}  Combining, $\tfrac{\gamma_{\bm D}}{2}h\le(L_B+m_B)\abs{z_1-z_2}$, which is
\eqref{E:tau-fd-lip}.
\end{proof}

\begin{theorem}[$\tau$ is Lipschitz on compacts]\label{T:tau-Lipschitz}\label{lem:tau-continuous}
For every compact $\bm D\Subset\mathbb D$, $\tau$ is Lipschitz on $\bm D$ with constant
\begin{equation}\label{E:global-Lip}
\widehat L_\tau(\bm D)\Def\max\Bigl\{L_\tau(\bm D),\ \frac{2\tau_{\max}(\bm D)}{\varrho_0(\bm D)}\Bigr\}.
\end{equation}
In particular $\tau\in\mW^{1,\infty}$ on a neighbourhood of $\bm D$, and $\tau$ is
differentiable at a.e.\ point of $\bm D$.
\end{theorem}

\begin{proof}
If $\abs{z_1-z_2}<\varrho_0$, use \cref{P:tau-fd-lip}.  If $\abs{z_1-z_2}\ge\varrho_0$, then
$\abs{\tau(z_1)-\tau(z_2)}\le2\tau_{\max}(\bm D)\le(2\tau_{\max}/\varrho_0)\abs{z_1-z_2}$.
For the $\mW^{1,\infty}$ statement, apply the same to a compact neighbourhood of $\bm D$ in
$\mathbb D$ (which exists since $\mathbb D$ is open); Rademacher's theorem gives a.e.\
differentiability.
\end{proof}

\subsection{\texorpdfstring{$C^1$}{C1} regularity in the state}\label{S:C1-state}

We now show that $\bm V=\tau$ is continuously differentiable in the state variables,
with a jointly continuous gradient, and clarify in which sense it is differentiable in
time.

\paragraph{Global prediction flow.}  Recall from \cref{L:flow} that the prediction field
$\mathsf G(s,Z)=(\Lambda,\ \sa_\mL(s)-\alpha[\mathcal V(\Xi)+\Lambda-v_\mL(s)])$ is globally
Lipschitz in $Z$, so the prediction from any $z=(s,\bm y)$ exists on all of
$[0,\timehor+1]$, through and beyond the collision set.  We write $Z(\kappa;s,\bm y)=(\Xi,\Lambda)(\kappa;s,\bm y)$ for the state of this global flow
at \emph{elapsed} time $\kappa\in[0,\timehor+1-s]$.  Using the global flow is what allows
differentiation at $\kappa=\tau(z)$ in the proof of \cref{lem:V-C1}: the map
$\bm y'\mapsto\Xi(\tau(z);s,\bm y')$ must be defined on a full neighbourhood of $\bm y$,
which includes states whose own collision occurs before $\tau(z)$.

The state derivative of the field is
\begin{equation}\label{E:dG-dy}
\partial_{\bm y}\mathsf G(s,Z)=\begin{pmatrix}0&1\\-\alpha\mathcal V'(\Xi)&-\alpha\end{pmatrix},
\qquad\norm{\partial_{\bm y}\mathsf G}\le G_1\Def\sqrt{1+2\alpha^2},
\end{equation}
continuous in the state, bounded, Lipschitz in the state (since $\mathcal V''$ is bounded),
and \emph{independent of $s$ and of the forcing}.

\begin{lemma}[State linearisation of the prediction flow]\label{lem:clean-linearization}
Fix $s\in[0,\timehor)$.  For every $\kappa\in[0,\timehor+1-s]$ the map
$\bm y\mapsto Z(\kappa;s,\bm y)$ is $C^1$ on $\RR^2$, with derivative
$Y(\kappa;s,\bm y)\Def\partial_{\bm y}Z(\kappa;s,\bm y)\in\RR^{2\times2}$ the solution of
the linear Carath\'eodory system
\begin{equation}\label{E:var-clean}
\partial_\kappa Y=\partial_{\bm y}\mathsf G\bigl(Z(\kappa;s,\bm y)\bigr)\,Y\quad\text{a.e. in }\kappa,
\qquad Y(0;s,\bm y)=I .
\end{equation}
The map $(\kappa,s,\bm y)\mapsto Y(\kappa;s,\bm y)$ is continuous, and
\begin{equation}\label{E:Y-two-sided}
e^{-G_1\kappa}\abs v\ \le\ \abs{Y(\kappa;s,\bm y)v}\ \le\ e^{G_1\kappa}\abs v\qquad(v\in\RR^2);
\end{equation}
in particular $Y$ is invertible with $\norm{Y^{-1}}\le e^{G_1\kappa}$.
\end{lemma}

\begin{proof}
\emph{Differentiability.}  Fix $\bm y_0$, a direction $v$, and $h\ne0$; write
$Z=Z(\cdot;s,\bm y_0)$, $Z^h=Z(\cdot;s,\bm y_0+hv)$ and $\xi_h\Def(Z^h-Z)/h$.  By
\eqref{E:flow-lip-z}, $\abs{Z^h-Z}\le L_B\abs h\abs v$ on a bounded set $B$ containing both.
By the fundamental theorem of calculus along the segment between $Z(r)$ and $Z^h(r)$,
\begin{align*}
\mathsf G(r,Z^h(r))-\mathsf G(r,Z(r))&=\Bigl(\int_0^1\partial_{\bm y}\mathsf G\bigl(Z(r)+\lambda(Z^h(r)-Z(r))\bigr)d\lambda\Bigr)(Z^h(r)-Z(r))\\
&=\partial_{\bm y}\mathsf G(Z(r))\,(Z^h(r)-Z(r))+R(r,h),
\end{align*}
with $\abs{R(r,h)}\le\omega(L_B\abs h\abs v)\,L_B\abs h\abs v$, $\omega$ the modulus of
continuity of $\partial_{\bm y}\mathsf G$ on a compact set (here $\omega(\rho)\le C\rho$).
Hence
\[
\xi_h(\kappa)=v+\int_0^\kappa\partial_{\bm y}\mathsf G(Z(r))\,\xi_h(r)\,dr+E_h(\kappa),
\qquad\abs{E_h(\kappa)}\le\kappa\,\omega(L_B\abs h\abs v)\,L_B\abs v\to0 .
\]
Let $\xi$ solve $\xi(\kappa)=v+\int_0^\kappa\partial_{\bm y}\mathsf G(Z)\xi$, a linear equation
with bounded measurable coefficient.  Gr\"onwall applied to $\xi_h-\xi$ gives
$\sup_{[0,\kappa_1]}\abs{\xi_h-\xi}\le\sup\abs{E_h}\,e^{G_1\kappa_1}\to0$.  Thus
$\partial_v Z(\kappa;s,\bm y_0)$ exists and equals $\xi(\kappa)=Y(\kappa;s,\bm y_0)v$, where $Y$
solves \eqref{E:var-clean}.  Continuity of $Y$ in $(\kappa,s,\bm y)$ follows from
continuous dependence of the linear system on its coefficient
$A(\kappa;s,\bm y)=\partial_{\bm y}\mathsf G(Z(\kappa;s,\bm y))$, which is jointly continuous
because $Z$ is (\cref{L:flow}(b)) and $\partial_{\bm y}\mathsf G$ is continuous in the state.
Continuity of $\bm y\mapsto Y(\kappa;s,\bm y)$ gives the $C^1$ claim.

\emph{Two-sided bound.}  With $w=Yv$ and $\sigma=\abs w^2$ (absolutely continuous),
$\abs{\dot\sigma}=2\abs{\langle w,Aw\rangle}\le2G_1\sigma$ a.e., so
$\tfrac{d}{d\kappa}(\sigma e^{\pm2G_1\kappa})\gtrless0$ a.e.; integrating (legitimate for
absolutely continuous functions) gives $e^{-2G_1\kappa}\abs v^2\le\sigma(\kappa)\le e^{2G_1\kappa}\abs v^2$.
\end{proof}

\begin{lemma}[Regularity of $\bm V$ on compacts]\label{lem:V-C1}
Let $\bm D\Subset\mathbb D$ be compact, $K\Subset\mathbb D$ a compact neighbourhood of $\bm D$,
and $\gamma_K$, $\tau_{\max}(K)$, $\widehat L_\tau(K)$ as in \cref{S:uniform-constants,S:tau-lipschitz}.
Write $\bm V\Def\tau$.
\begin{enumerate}[label=\tup{\roman*}]
\item\label{it:C1-state} \textbf{($C^1$ in the state.)}  For every $(s,\bm y)\in\bm D$ the
  map $\bm y'\mapsto\bm V(s,\bm y')$ is continuously differentiable on a neighbourhood of
  $\bm y$, with
  \begin{equation}\label{eq:grad-tau}
  \nabla_{\bm y}\bm V(s,\bm y)
  =-\,\frac{\bm e_1^\top Y\bigl(\tau(z);s,\bm y\bigr)}{\Lambda\bigl(\tau(z);z\bigr)},
  \qquad
  \abs{\nabla_{\bm y}\bm V(s,\bm y)}\le\frac{e^{G_1\tau_{\max}(K)}}{\gamma_K} ,
  \end{equation}
  and $(s,\bm y)\mapsto\nabla_{\bm y}\bm V(s,\bm y)$ is continuous on $\bm D$.
\item\label{it:Lip-time} \textbf{(Lipschitz in time; a.e.\ classical HJB with a common
  exceptional time set.)}  $\bm V$ is Lipschitz on $\bm D$.  Let $N\subset[0,\timehor)$ be
  the null set of times that are not Lebesgue points of $\sa_\mL$.  At every
  $(s,\bm y)\in\bm D$ with $s\notin N$ the time derivative $\partial_s\bm V(s,\bm y)$
  exists, $\bm V$ is Fr\'echet-differentiable, and, with $\sa_\mL(s)$ the Lebesgue value,
  \begin{equation}\label{eq:HJB-classical}
  \partial_s\bm V(s,\bm y)+\mathsf G(s,\bm y)\cdot\nabla_{\bm y}\bm V(s,\bm y)=-1 .
  \end{equation}
  Consequently \eqref{eq:HJB-classical} holds for a.e.\ $s$ and \emph{every} $\bm y$ with
  $(s,\bm y)\in\bm D$ (the exceptional set of times $N$ does not depend on $\bm y$), and
  for a.e.\ $(s,\bm y)\in\bm D$.  Moreover, at every $(s,\bm y)\in\bm D$ at which
  $\partial_s\bm V(s,\bm y)$ exists (for fixed $\bm y$, a.e.\ $s$), $\bm V$ is
  Fr\'echet-differentiable.
\item\label{it:C1-joint} \textbf{(Joint $C^1$ where the forcing is continuous.)}  If
  $\sa_\mL$ is continuous on an open interval $J\subset(0,\timehor)$, then $\bm V\in C^1$
  on $\mathbb D\cap(J\times\RR^2)$, with $\partial_s\bm V=-1-\mathsf G\cdot\nabla_{\bm y}\bm V$
  there.  (No condition is imposed on the forcing \emph{after} $s$: the prediction from
  $(s,\bm y)$ may cross jumps of $\sa_\mL$.)
\end{enumerate}
\end{lemma}

\begin{proof}
\emph{\tup{i}.}  Fix $z^0=(s^0,\bm y^0)\in\bm D$, $\nu^0\Def\tau(z^0)$,
$\gamma\Def-\Lambda(\nu^0;z^0)\ge\gamma_K$.  For $\bm y$ near $\bm y^0$ (so that
$z\Def(s^0,\bm y)\in K$), let $\nu\Def\tau(z)$.

\emph{Step 1 (identity at the common elapsed time $\nu^0$).}  Since the flow is global,
$\Xi(\cdot;z)$ is absolutely continuous on $[0,\timehor+1-s^0]$ with derivative $\Lambda$,
and $\Xi(\nu;z)=0=\Xi(\nu^0;z^0)$.  Hence, as a signed identity with no case split,
\begin{equation}\label{eq:fd-identity}
\Xi(\nu^0;z)-\Xi(\nu^0;z^0)=\Xi(\nu^0;z)-\Xi(\nu;z)=-\int_{\nu^0}^{\nu}\Lambda(\theta;z)\,d\theta .
\end{equation}

\emph{Step 2 (left side).}  By \cref{lem:clean-linearization} at the fixed elapsed time
$\nu^0$ and start time $s^0$, $\bm y\mapsto\Xi(\nu^0;s^0,\bm y)$ is $C^1$ on $\RR^2$, so
\begin{equation}\label{eq:lhs-expansion}
\Xi(\nu^0;s^0,\bm y)-\Xi(\nu^0;s^0,\bm y^0)=\bm e_1^\top Y(\nu^0;s^0,\bm y^0)(\bm y-\bm y^0)+o(\abs{\bm y-\bm y^0}).
\end{equation}

\emph{Step 3 (right side).}  For $\theta$ between $\nu^0$ and $\nu$,
\[
\begin{aligned}
&\abs{\Lambda(\theta;z)-\Lambda(\nu^0;z^0)}\\
&\quad\le\abs{\Lambda(\theta;z)-\Lambda(\nu^0;z)}
 +\abs{\Lambda(\nu^0;z)-\Lambda(\nu^0;z^0)}\\
&\quad\le m_B\abs{\nu-\nu^0}+L_B\abs{\bm y-\bm y^0}\Def\varepsilon(\bm y),
\end{aligned}
\]
by elapsed-time Lipschitz continuity along the single trajectory $z$ and
\eqref{E:flow-lip-z} at the common absolute time $s^0+\nu^0$ (same start time).  So
\begin{equation}\label{eq:rhs-expansion}
-\int_{\nu^0}^{\nu}\Lambda(\theta;z)\,d\theta=\gamma\,(\nu-\nu^0)+R,\qquad\abs R\le\varepsilon(\bm y)\abs{\nu-\nu^0}.
\end{equation}

\emph{Step 4 (remainder).}  By \cref{T:tau-Lipschitz} on $K$,
$\abs{\nu-\nu^0}\le\widehat L_\tau(K)\abs{\bm y-\bm y^0}$, so
$\abs R\le(m_B\widehat L_\tau(K)+L_B)\widehat L_\tau(K)\abs{\bm y-\bm y^0}^2=o(\abs{\bm y-\bm y^0})$.

\emph{Step 5 (conclusion).}  Combining,
$\bm e_1^\top Y(\nu^0;s^0,\bm y^0)(\bm y-\bm y^0)=\gamma(\tau(z)-\tau(z^0))+o(\abs{\bm y-\bm y^0})$,
so $\bm y\mapsto\tau(s^0,\bm y)$ is differentiable at $\bm y^0$ with gradient
\eqref{eq:grad-tau}.  The bound uses $\norm{Y(\nu^0)}\le e^{G_1\nu^0}\le e^{G_1\tau_{\max}(K)}$
from \eqref{E:Y-two-sided} and $\gamma\ge\gamma_K$.

\emph{Joint continuity of the gradient.}  $z\mapsto\Lambda(\tau(z);z)$ is continuous on
$K$ and $\le-\gamma_K$ (\cref{L:unif-transv}); $z\mapsto Y(\tau(z);z)$ is continuous as the
composition of the continuous $\tau$ with the jointly continuous $Y$
(\cref{lem:clean-linearization}).  Hence \eqref{eq:grad-tau} is continuous on $K$, and in
particular $\bm y\mapsto\bm V(s^0,\bm y)$ is $C^1$ near $\bm y^0$.

\medskip\emph{\tup{ii}.}  Lipschitz continuity is \cref{T:tau-Lipschitz}.

\emph{Fr\'echet differentiability.}  At a point $(s^0,\bm y^0)\in\bm D$ where
$\partial_s\bm V$ exists, all three partial derivatives exist, and the two state partials
exist on a neighbourhood (in $K$) and are continuous there by (i).  A function whose
partial derivatives all exist at a point, all but one of them being continuous on a
neighbourhood, is Fr\'echet-differentiable at that point \cite[Thm.~12.11]{Apostol}.  (For
fixed $\bm y$, $s\mapsto\bm V(s,\bm y)$ is Lipschitz, so this applies for a.e.\ $s$.)

\emph{Existence of $\partial_s\bm V$ and the HJB identity at Lebesgue points.}  Fix
$z^0=(s^0,\bm y^0)\in\bm D$ with $s^0\notin N$, and let $h\ne0$ be small (two-sided if
$s^0>0$, $h>0$ if $s^0=0$).  Set $\tilde{\bm y}_h\Def Z^{\rm abs}(s^0;(s^0+h,\bm y^0))$,
the state at absolute time $s^0$ of the prediction through $(s^0+h,\bm y^0)$ (defined by
\cref{L:flow}(a), backward as well as forward).  By the cocycle property,
$(s^0+h,\bm y^0)$ and $(s^0,\tilde{\bm y}_h)$ lie on one characteristic, at elapsed distance
$\abs h$.  Since $\tilde{\bm y}_h\to\bm y^0$ as $h\to0$, $\mathbb D$ is open, and
$\tau\ge\tau_{\min}(K)>0$ on $K$, for small $\abs h$ both points lie in $K$ and $\abs h<\tau$
at the earlier of the two, so \eqref{E:DPP} gives, for either sign of $h$,
\begin{equation}\label{eq:shift-identity}
\tau(s^0+h,\bm y^0)=\tau(s^0,\tilde{\bm y}_h)-h .
\end{equation}
(For $h>0$ the point $(s^0+h,\bm y^0)$ is at elapsed time $h$ after $(s^0,\tilde{\bm y}_h)$;
for $h<0$ the point $(s^0,\tilde{\bm y}_h)$ is at elapsed time $-h$ after
$(s^0+h,\bm y^0)$; both readings give \eqref{eq:shift-identity}.)
Next, $\abs{Z^{\rm abs}(r;(s^0+h,\bm y^0))-\bm y^0}\le m_B\abs h$ for $r$ between $s^0$
and $s^0+h$, and $\mathsf G$ is Lipschitz in the state, so
\[
\begin{aligned}
\tilde{\bm y}_h
&=\bm y^0-\int_{s^0}^{s^0+h}
  \mathsf G\bigl(r,Z^{\rm abs}(r;(s^0+h,\bm y^0))\bigr)\,dr\\
&=\bm y^0-\int_{s^0}^{s^0+h}\mathsf G(r,\bm y^0)\,dr+O(h^2)\\
&=\bm y^0-h\,\mathsf G(s^0,\bm y^0)+o(h),
\end{aligned}
\]
because $s^0$ is a Lebesgue point of $r\mapsto\mathsf G(r,\bm y^0)$ (only $\sa_\mL$ is
discontinuous, and a Lebesgue point is two-sided).  By (i), $\bm y\mapsto\tau(s^0,\bm y)$
is differentiable at $\bm y^0$, so
$\tau(s^0,\tilde{\bm y}_h)=\tau(z^0)-h\,\nabla_{\bm y}\tau(z^0)\cdot\mathsf G(s^0,\bm y^0)+o(h)$.
Inserting this in \eqref{eq:shift-identity},
\[
\tau(s^0+h,\bm y^0)-\tau(z^0)=-h\bigl(1+\mathsf G(s^0,\bm y^0)\cdot\nabla_{\bm y}\tau(z^0)\bigr)+o(h),
\]
so $\partial_s\bm V(z^0)$ exists and equals $-1-\mathsf G(s^0,\bm y^0)\cdot\nabla_{\bm y}\bm V(z^0)$,
which is \eqref{eq:HJB-classical}; Fr\'echet differentiability at $z^0$ then follows from the
previous paragraph.  No uniformity in $\bm y^0$ was needed: the exceptional set is $N$,
determined by $\sa_\mL$ alone.

\emph{Measure statements.}  $N$ is null, and Fubini gives the a.e.\ statement on $\bm D$.

\medskip\emph{\tup{iii}.}  Let $F(\kappa,s,\bm y)\Def\Xi(\kappa;s,\bm y)$ on
$\{(\kappa,s,\bm y):s\in J,\ \kappa\in[0,\timehor-s)\}$.  We show $F\in C^1$ there.
$\partial_\kappa F=\Lambda(\kappa;s,\bm y)$ (for every $\kappa$, since $\Xi$ is $C^1$ in
elapsed time) and $\partial_{\bm y}F=\bm e_1^\top Y(\kappa;s,\bm y)$
(\cref{lem:clean-linearization}); both are continuous in $(\kappa,s,\bm y)$.

For $\partial_sF$, write $F(\kappa,s,\bm y)=\Xi^{\rm abs}(s+\kappa;(s,\bm y))$ and vary $s$
at fixed absolute time $t=s+\kappa$ and at fixed $\kappa$ separately.  At fixed $t$: by
the cocycle property, $Z^{\rm abs}(t;(s+h,\bm y))=Z^{\rm abs}(t;(s,\tilde{\bm y}_h))$ with
$\tilde{\bm y}_h\Def Z^{\rm abs}(s;(s+h,\bm y))$ the backward state at time $s$ of the
prediction from $(s+h,\bm y)$.  Since $s$ is a continuity point of $\sa_\mL$ (as $s\in J$),
$\tilde{\bm y}_h=\bm y-h\,\mathsf G(s,\bm y)+o(h)$.  By \cref{lem:clean-linearization} for the
start time $s$ and absolute end time $t$ (elapsed time $t-s$),
$Z^{\rm abs}(t;(s,\tilde{\bm y}_h))=Z^{\rm abs}(t;(s,\bm y))-h\,Y(t-s;s,\bm y)\mathsf G(s,\bm y)+o(h)$.
Hence
\[
\partial_s\big|_{t}Z^{\rm abs}(t;(s,\bm y))=-Y(t-s;s,\bm y)\,\mathsf G(s,\bm y),
\]
\emph{whatever happens to the forcing between $s$ and $t$}: the only ingredients are the
continuity of $\sa_\mL$ at $s$ and the forcing-free variational equation.  At fixed
$\kappa$, $\partial_s F=\partial_t\Xi^{\rm abs}+\partial_s|_t\Xi^{\rm abs}
=\Lambda(\kappa;s,\bm y)-\bm e_1^\top Y(\kappa;s,\bm y)\mathsf G(s,\bm y)$, which is continuous in
$(\kappa,s,\bm y)$ for $s\in J$ (continuity of $\Lambda$, $Y$, and of $\mathsf G(\cdot,\bm y)$
on $J$).  So $F\in C^1$.

At $\kappa=\tau(s,\bm y)$ we have $F=0$ and $\partial_\kappa F=\Lambda(\tau;z)\le-\gamma_K<0$.
The implicit function theorem gives a $C^1$ function $\hat\kappa(s,\bm y)$ near $z$ with
$F(\hat\kappa,s,\bm y)=0$ and $\hat\kappa(z)=\tau(z)$; since $\tau$ is continuous and
$\hat\kappa$ is the unique zero of $F(\cdot,s',\bm y')$ near $\tau(z)$, $\hat\kappa=\tau$ on a
neighbourhood.  Hence $\tau\in C^1$ there, with $\nabla_{\bm y}\tau$ as in
\eqref{eq:grad-tau} and
\[
\partial_s\tau=-\frac{\partial_sF}{\partial_\kappa F}
=-\frac{\Lambda(\tau;z)-\bm e_1^\top Y(\tau;s,\bm y)\mathsf G(s,\bm y)}{\Lambda(\tau;z)}
=-1-\mathsf G(s,\bm y)\cdot\nabla_{\bm y}\tau .
\]
\end{proof}

\begin{remark}[$\bm V$ is not $C^1$ in time under rough forcing]\label{R:not-C1-in-time}
Suppose $\sa_\mL$ is continuous on $(t_1-\eta,t_1)$ and on $(t_1,t_1+\eta)$ with distinct
one-sided limits $\sa_\mL(t_1^\pm)$.  By \eqref{it:C1-joint} on each side, $\partial_s\bm V$
has one-sided limits at $s=t_1$ for each $\bm y$, and by the formula in \eqref{it:C1-joint},
using that $\sa_\mL$ enters $\mathsf G(s,\bm y)$ only through its second component and that
$\nabla_{\bm y}\bm V$ is continuous across $t_1$ (part~(i)),
\[
\partial_s\bm V(t_1^+,\bm y)-\partial_s\bm V(t_1^-,\bm y)
=-\bigl[\sa_\mL(t_1^+)-\sa_\mL(t_1^-)\bigr]\,\partial_{\Lambda}\bm V(t_1,\bm y).
\]
Moreover $\partial_\Lambda\bm V(t_1,\bm y)>0$ whenever
$\tau(t_1,\bm y)<\pi/\omega_\alpha$, $\omega_\alpha\Def\sqrt{\alpha(1-\alpha/4)}$ (and at
every state if $\alpha\ge4$): by \eqref{eq:grad-tau},
$\partial_\Lambda\bm V=-Y_{12}(\tau;z)/\Lambda(\tau;z)=Y_{12}(\tau;z)/\abs{\Lambda(\tau;z)}$, and
$w\Def Y_{12}$ (second column of \eqref{E:var-clean}) solves
$w''+\alpha w'+\alpha\mathcal V'(\Xi)w=0$, $w(0)=0$, $w'(0)=1$; the substitution
$u=e^{\alpha\kappa/2}w$ gives $u''+(\alpha\mathcal V'(\Xi)-\alpha^2/4)u=0$ with coefficient
$\le\alpha-\alpha^2/4=\omega_\alpha^2$ (as $0<\mathcal V'\le1$), so by Sturm comparison $u$,
hence $w$, has no zero in $(0,\pi/\omega_\alpha)$, and none at all if $\alpha\ge4$ (this is
the argument of \cref{L:fundamental}(b),(c) with $\tilde c=\mathcal V'(\Xi)$).  At such
points $\partial_s\bm V$ jumps and $\bm V\notin C^1$.  For $\tau\ge\pi/\omega_\alpha$ the
sensitivity may vanish, and where it does the jump of $\partial_s\bm V$ degenerates; the
conclusion ``$\bm V\notin C^1$ near $(t_1,\bm y)$'' requires
$\partial_\Lambda\bm V(t_1,\bm y)\ne0$.  (In the numerical setting $\alpha=2$,
$\timehor=1$ one has $\pi/\omega_\alpha=\pi>1\ge\tau$, so the sensitivity is positive
everywhere on $\mathbb D$.)  The jump
formula and the positivity of $\partial_\Lambda\bm V$ below the sufficient threshold are checked numerically in \Cref{S:num-regularity}. The global variational-flow comparison and the separate first-contact sensitivity experiment do not establish optimal sharpness of the TTC threshold.  The Lipschitz-in-time regularity of
\eqref{it:Lip-time} is therefore sharp.
\end{remark}


\section{Hamilton--Jacobi characterisation and residual stability}\label{sec: collision_set}\label{sec:HJ_collision}
Repeated evaluation of $\tau$ requires prediction and contact detection, which motivates a reusable approximation. We analyse the continuum residual of the Hamilton--Jacobi equation satisfied by $\tau$, together with the error in its incoming boundary datum. For locally Lipschitz candidates with a continuous incoming trace, the $L^2$ error is controlled by these two losses. The co-area proof gives a constant determined by bounded-flow estimates; a second proof in endpoint/age coordinates gives an explicit coefficient involving only the maximal flight time, the terminal closing speed, and $\alpha$. Both estimates are uniform as the transversality margin decreases within bounded families. They imply convergence conditional on vanishing continuum loss; they do not assert convergence of a training algorithm or conversion of sampled loss into continuum loss. We keep the notation of \cref{sec:ttc-regularity}.

\subsection{Value function, HJB equation, and the loss}\label{S:HJB-setup}\label{S:HJB-strong-PINNs}

\paragraph{The incoming collision surface.}  Let
\begin{equation}\label{E:Sigma-minus}
\Sigma^-\Def\{(t,0,\Lambda):\ 0\le t<\timehor,\ \Lambda<0\}
\end{equation}
be the \emph{incoming} part of the collision plane $\partial\mathcal Q=\{\Xi=0\}$.

\paragraph{Value function.}  Set $\bm V\Def\tau$ on $\mathbb D$ and $\bm V\Def0$ on
$\Sigma^-$.  By \eqref{E:DPP} (optimality principle), for $z=(s,\bm y)\in\mathbb D$ and
$0\le\kappa\le\bm V(s,\bm y)$,
\begin{equation}\label{E:DPP-V}
\bm V\bigl(s+\kappa,Z(\kappa;z)\bigr)=\bm V(s,\bm y)-\kappa ,
\end{equation}
where at $\kappa=\bm V(s,\bm y)$ the left side is evaluated at the collision point
\[
(s+\tau(z),0,\Lambda(\tau(z);z)),
\]
which lies in $\Sigma^-$ because
$\Lambda(\tau(z);z)<0$ and $s+\tau(z)<\timehor$ by definition of $\mathbb D$, and equals $0$.

\paragraph{HJB equation.}  Differentiating \eqref{E:DPP-V} at $\kappa=0$, formally,
gives the transport equation
\begin{equation}\label{E:HJB}
\partial_s\bm V+\mathsf G(s,\bm y)\cdot\nabla_{\bm y}\bm V=-1\ \text{ on }\mathbb D,\qquad
\bm V=0\ \text{ on }\Sigma^- .
\end{equation}
The senses in which $\bm V$ solves \eqref{E:HJB} are made precise in \cref{P:HJB-senses}:
it satisfies the integrated form \eqref{E:DPP-V} exactly along every characteristic; it
satisfies \eqref{E:HJB} classically at every state for all times outside a fixed null set
determined by $\sa_\mL$, hence a.e.\ on $\mathbb D$ (\cref{lem:V-C1}\eqref{it:Lip-time});
and it is a classical $C^1$ solution on every time interval on which $\sa_\mL$ is
continuous (\cref{lem:V-C1}\eqref{it:C1-joint}).  Writing \eqref{E:HJB} as
$\partial_s\bm V+H(s,\bm y,\nabla_{\bm y}\bm V)=0$, the Hamiltonian
$H(s,\bm y,p)=p\cdot\mathsf G(s,\bm y)+1$ is affine in $p$, so the characteristics are the
predicted OVM trajectories, independent of $p$, and they do not cross.  This is what makes
the equation integrable along characteristics from the collision set, which is all the
proofs below use: a solution that is absolutely continuous along characteristics---for
instance a locally Lipschitz a.e.\ solution---is thereby determined by its values on
$\Sigma^-$ (\cref{P:HJB-senses}(d),(e)), whereas a merely continuous a.e.\ solution is
not (\cref{R:viscosity}).  Whether $\bm V$ is also the unique viscosity solution of
\eqref{E:HJB} is discussed in \cref{R:viscosity}; since $\sa_\mL$ is only measurable in
time, that question requires the $L^1$-in-time framework. The estimates below rely on the
characteristic identity \eqref{E:DPP-V}.

\paragraph{Topological conventions.}  All sets of states lie in the half-space
$\mathcal H\Def[0,\timehor)\times\RR^2$; ``open'', ``neighbourhood'' and ``interior''
refer to the relative topology of $\mathcal H$ unless $\operatorname{int}_{\RR^3}$ is
written.  The face $\{s=0\}$ is $\mathcal L^3$-null, so almost-everywhere statements are
unaffected by it; Rademacher's theorem and the area and co-area formulas are applied on
the $\RR^3$-open set $\operatorname{int}_{\RR^3}K=\operatorname{int}K\setminus\{s=0\}$ of
a compact neighbourhood $K$, which contains $\bm D$ up to a null set.  At states with
$s=0$, derivatives in $s$ are one-sided.  (The prediction flow is defined for absolute
times in $[0,\timehor+1]$ only, \cref{L:flow}; there is no negative-time extension, and
none is used.)

\paragraph{Forward tube, endpoint map, boundary set.}  Fix a compact
$\bm D\Subset\mathbb D$ and set
\begin{equation}\label{E:Dplus}
\begin{aligned}
\bm D^+\Def\bigl\{&(s+\kappa,Z(\kappa;s,\bm y)):\ (s,\bm y)\in\bm D,\\
&0\le\kappa<\bm V(s,\bm y)\bigr\}\subset\mathbb D,
\end{aligned}
\end{equation}
where the inclusion follows from \cref{L:D-open}(b). For $z=(s,\bm y)\in\mathbb D$, define
\begin{equation}\label{E:Psi-def}
\begin{aligned}
\Psi(z)&\Def\bigl(s+\tau(z),\,Z(\tau(z);z)\bigr)\in\Sigma^-,\\
\tilde\Psi(z)&\Def\bigl(s+\tau(z),\,\Lambda(\tau(z);z)\bigr)\in\RR^2,
\end{aligned}
\end{equation}
and set
\begin{equation}\label{E:QD-def}
\mathcal Q_{\bm D}\Def\Psi(\bm D^+)=\Psi(\bm D)\subset\Sigma^-\subset\partial\mathcal Q .
\end{equation}
$\partial\mathcal Q=\{\Xi=0\}$ is a plane in $(t,\Xi,\Lambda)$-space, and
$\tilde\Psi$ is $\Psi$ written in the chart $(t,\Lambda)$ of that plane; the chart is an
isometry, so $\abs{\Psi(z_1)-\Psi(z_2)}=\abs{\tilde\Psi(z_1)-\tilde\Psi(z_2)}$ and
$\mathcal H^2\llcorner\partial\mathcal Q$ corresponds to Lebesgue measure on $\RR^2$.
Note $\Psi(\bm D^+)=\Psi(\bm D)$ because points on one characteristic share their endpoint.

\begin{lemma}[Continuity of the zero extension on the incoming surface]\label{L:V-cont-Sigma}
$\Sigma^-\subset\cl[\mathbb D]\setminus\mathbb D$, and $\bm V$ is continuous on
$\mathbb D\cup\Sigma^-$.  For every compact $\bm D\Subset\mathbb D$ one has
$\mathcal Q_{\bm D}\subset\Sigma^-$ (indeed $\Lambda\le-\gamma_{\bm D}$ on $\mathcal Q_{\bm D}$)
and $\cl[\bm D^+]=\bm D^+\cup\mathcal Q_{\bm D}$.
\end{lemma}
\begin{proof}
Let $p=(t^\ast,0,\Lambda^\ast)\in\Sigma^-$ and $\gamma\Def-\Lambda^\ast/2>0$.  The
prediction from $p$ exists on $[0,\timehor+1]$ (\cref{L:flow}); choose $h\in(0,1)$ with
$t^\ast+h<\timehor$ and $\Lambda^{\rm abs}(r;p)\le-\tfrac32\gamma$ for
$r\in I\Def[\max\{0,t^\ast-h\},\,t^\ast+h]$ (the prediction is defined for absolute times
in $[0,\timehor+1]$ only, so the interval is clipped at $0$; $t^\ast=0$ is allowed).
Let $B$ be the closed unit ball around $p$ intersected with $[0,\timehor)\times\RR^2$, and
$L_B$ as in \cref{L:flow}(b); pick $\rho\in(0,h/2)$ with $L_B\rho\le\gamma/2$ and
$\rho\le\gamma h/2$.  For $z=(s,X,\Lambda)\in\Om$ with $\abs{z-p}<\rho$ (so $s\ge0$ and
$s\in I$), \eqref{E:flow-lip-z} gives $\Lambda^{\rm abs}(r;z)\le-\gamma$ on $I\ni s$, hence
$\Xi^{\rm abs}(r;z)=X+\int_s^r\Lambda^{\rm abs}\le X-\gamma(r-s)$ for $r\in[s,t^\ast+h]$.
So the first zero of the predicted gap lies in $[s,s+X/\gamma]\subset[s,t^\ast+h]$, where
$\Lambda\le-\gamma<0$: it is transversal and occurs before $\timehor$.  Thus
$z\in\mathbb D$, which shows $p\in\cl[\mathbb D]$ (and $p\notin\mathbb D\subset\Om$), and
$\tau(z)\le X/\gamma\to0$ as $z\to p$.  Since $\bm V=0$ on $\Sigma^-$, $\bm V$ is
continuous at $p$; continuity on $\mathbb D$ is \cref{L:unif-transv}(a).  The inclusion
$\mathcal Q_{\bm D}\subset\Sigma^-$ is \cref{L:unif-transv}(b) together with
$s+\tau(z)=\sigma(z)<\timehor$.  Finally $\bm D^+\cup\mathcal Q_{\bm D}$ is the image of the
compact set $\{(\kappa,z):z\in\bm D,\ 0\le\kappa\le\tau(z)\}$ under the continuous map
$(\kappa,z)\mapsto(s+\kappa,Z(\kappa;z))$, hence compact, and every point of
$\mathcal Q_{\bm D}$ is the limit of $(s+\kappa,Z(\kappa;z))$ as $\kappa\uparrow\tau(z)$.
\end{proof}

\begin{remark}[Incoming data and extension at other boundary points]\label{R:no-global-ext}
The zero datum in \eqref{E:HJB} is specific to the incoming surface $\Sigma^-$.
At opening points it can conflict with the limiting TTC. For example, with $\alpha=0.1$,
$\timehor=3$ and the smooth admissible leader data $v_\mL=u+q'$,
$\sa_\mL=u'+q''$, where $q(t)=0.01t(2-t)$ and
$u'=\alpha(\mathcal V(q)-u)$, $u(0)=0.5$, the prediction from
$z_0=(0,0,0.02)$ is $(\Xi,\Lambda)=(q,q')$. The perturbed states
$z_\delta=(0,\delta,0.02)$ belong to $\mathbb D$ for small $\delta>0$ and
$\tau(z_\delta)\to2$. Thus the zero extension fails at this opening point, although
extension by the value $2$ is continuous there.

At grazing, discontinuity of the capped TTC on the ambient domain $\Om$ must be
distinguished from nonextendability of its restriction to $\mathbb D$.
Consider a nondegenerate interior graze as in \cref{L:window}, and suppose that its
continuation has no later zero up to and including $\timehor$. After shrinking the
neighbourhood $U$, only the hitting side belongs to the first-contact domain:
$\mathbb D\cap U=U^-$. On that side set $\mu=-m(z)>0$. The lower acceleration bound
in \cref{L:window}(a) gives
\[
\mu=\int_{\sigma(z)}^{s^\ast(z)}|\Lambda(r;z)|\,dr
\ge\frac{c_0}{2}\bigl(s^\ast(z)-\sigma(z)\bigr)^2.
\]
Consequently $\sigma(z)\to s^\ast(z_0)$ as $z\to z_0$ within $\mathbb D$, so the
restriction $\bm V$ extends there by the grazing flight time
$s^\ast(z_0)-\bar t_0$. The ambient capped TTC can nevertheless jump to the miss-side
value, whose states are outside $\mathbb D$. If the miss side instead has a later
transversal first contact before $\timehor$, approaches from the two sides within
$\mathbb D$ have different limiting contact times and no continuous extension exists
at that grazing point. Thus no universal nonextension statement at grazing is needed.
All candidate classes below require continuity on $\mathbb D\cup\Sigma^-$ only,
and boundary data are imposed on $\Sigma^-$, not on the entire collision plane or on
the grazing boundary.
\end{remark}

\paragraph{Candidate class and continuum loss.}\label{def:augmented-loss}
We allow $\bm V_\theta\in C^{0,1}_{\mathrm{loc}}(\mathbb D)\cap C(\mathbb D\cup\Sigma^-)$:
the candidate is locally Lipschitz in the interior and has a continuous incoming trace.
Fix a Borel representative of $\sa_\mL$. The candidate partial derivatives are taken
almost everywhere; choose Borel representatives, arbitrarily defined on the null set of
nondifferentiability. The HJB residual and the
two continuum loss components are
\begin{equation}\label{E:loss}
r_\theta\Def\partial_s\bm V_\theta+\mathsf G\cdot\nabla_{\bm y}\bm V_\theta+1,\qquad
J_{\mathrm{HJB}}\Def\int_{\bm D^+}\abs{r_\theta}^2\,d\mathcal L^3,\qquad
J_{\partial\mathcal Q}\Def\int_{\mathcal Q_{\bm D}}\abs{\bm V_\theta}^2\,d\mathcal H^2 ,
\end{equation}
and their weighted sum is
\begin{equation}\label{eq:J-total}
J_{\mathrm{total}}(\theta)\Def\lambda_{\mathrm{HJB}}\,J_{\mathrm{HJB}}+\lambda_{\partial\mathcal Q}\,J_{\partial\mathcal Q},\qquad\lambda_{\mathrm{HJB}},\lambda_{\partial\mathcal Q}>0 .
\end{equation}
The residual $r_\theta$ is measurable and locally essentially bounded; even for a $C^1$
candidate it is generally not continuous under $L^\infty$ forcing, since it contains
$\sa_\mL(s)\,\partial_\Lambda\bm V_\theta$. Its a.e. representative does not affect the
loss or the characteristic identities valid for a.e. initial point below. Interior
$C^1$ candidates additionally admit the pointwise improper identity in
\cref{L:error-char}. The constants refer to the unweighted Lebesgue and surface measures
in \eqref{E:loss}. Monte-Carlo losses under other sampling laws require the corresponding
density weights to approximate these integrals; equivalence of measures alone does not
identify their values. The additional penalties and sampling choices used in
\cref{S:num-learning} are stated separately.

\paragraph{Bounded reference set and constants.}  Fix a compact neighbourhood
$K\Subset\mathbb D$ of $\bm D$ (it exists since $\mathbb D$ is open), with forward tube
$K^+$ defined as in \eqref{E:Dplus}.  Let $B$ be the closed $1$-neighbourhood of the
compact set $\cl[K^+]$, intersected with $[0,\timehor)\times\RR^2$, and let $m_B,L_B$ be
the constants of \cref{L:flow}(b) for $B$.  Recall from \cref{L:flow} that
the prediction flow is defined, forward and backward, on all of $[0,\timehor+1]$, and
that \eqref{E:flow-consts} holds for all absolute times, in particular for \emph{backward}
evaluation.  Let $\tau_{\max}=\tau_{\max}(\bm D)$ and $\widehat L_\tau(\bm D)$ be as in
\cref{T:tau-Lipschitz}.

\subsection{The error identity along characteristics}\label{S:error-identity}

\begin{lemma}[Error along a characteristic]\label{L:error-char}
In this lemma the residual is computed from the classical derivatives of the $C^1$ candidate at every interior point. Let $\bm V_\theta\in C^1(\mathbb D)\cap C(\mathbb D\cup\Sigma^-)$, $e\Def\bm V_\theta-\bm V$,
$z=(s,\bm y)\in\bm D$, $\nu\Def\bm V(z)$, and $c(\kappa)\Def(s+\kappa,Z(\kappa;z))$ for
$\kappa\in[0,\nu]$.
\begin{enumerate}[label=\tup{\alph*}]
\item The improper integral
$\int_0^{\nu^-}r_\theta(c(\kappa))\,d\kappa\Def\lim_{\nu'\uparrow\nu}\int_0^{\nu'}r_\theta(c(\kappa))\,d\kappa$
exists, and
\begin{equation}\label{E:error-char}
e(s,\bm y)=\bm V_\theta\bigl(\Psi(z)\bigr)-\int_0^{\nu^-}r_\theta\bigl(s+\kappa,Z(\kappa;z)\bigr)\,d\kappa .
\end{equation}
\item For every $z\in\bm D$, as an inequality in $[0,\infty]$,
\begin{equation}\label{E:pointwise-bound}
\abs{e(s,\bm y)}^2\le2\,\abs{\bm V_\theta(\Psi(z))}^2+2\nu\int_0^{\nu}\abs{r_\theta(s+\kappa,Z(\kappa;z))}^2\,d\kappa .
\end{equation}
\item If the integral on the right of \eqref{E:pointwise-bound} is finite, then
$r_\theta\circ c\in L^1(0,\nu)$ and the integral in \eqref{E:error-char} is an ordinary
Lebesgue integral.  This is the case for $\mathcal L^3$-a.e.\ $z\in\bm D$ whenever
$J_{\mathrm{HJB}}<\infty$ (by Tonelli and \cref{P:I2}, whose proof does not use this
lemma), and for every $z\in\bm D$ if $\bm V_\theta$ is $C^1$ on a neighbourhood of
$\cl[\bm D^+]$ (then $r_\theta\circ c\in L^\infty(0,\nu)$).
\end{enumerate}
\end{lemma}

\begin{proof}
The curve $c$ is Lipschitz, lies in $\bm D^+\subset\mathbb D$ for $\kappa<\nu$, and ends at
$\Psi(z)\in\mathcal Q_{\bm D}\subset\Sigma^-$ (\cref{L:V-cont-Sigma}).  Since $\bm V_\theta\in C^1(\mathbb D)$,
$\bm V_\theta\circ c$ is absolutely continuous on $[0,\nu']$ for every $\nu'<\nu$, with
$(\bm V_\theta\circ c)'(\kappa)=\partial_s\bm V_\theta+\nabla_{\bm y}\bm V_\theta\cdot\partial_\kappa Z
=\partial_s\bm V_\theta+\nabla_{\bm y}\bm V_\theta\cdot\mathsf G=r_\theta(c(\kappa))-1$ for a.e.\
$\kappa$ (the chain rule for a $C^1$ function along a Lipschitz curve, and
$\partial_\kappa Z=\mathsf G(s+\kappa,Z)$ a.e.).  Also $\bm V(c(\kappa))=\nu-\kappa$ by
\eqref{E:DPP-V}.  Hence $(e\circ c)'=r_\theta\circ c$ a.e.\ on $[0,\nu']$ and
\begin{equation}\label{E:error-char-trunc}
e(c(\nu'))-e(z)=\int_0^{\nu'}r_\theta(c(\kappa))\,d\kappa\qquad(0\le\nu'<\nu).
\end{equation}
As $\nu'\uparrow\nu$, $c(\nu')\to\Psi(z)$, so $\bm V_\theta(c(\nu'))\to\bm V_\theta(\Psi(z))$
by continuity of $\bm V_\theta$ on $\mathbb D\cup\Sigma^-$ at $\Psi(z)$, and $\bm V(c(\nu'))=\nu-\nu'\to0$.
Thus the left side of \eqref{E:error-char-trunc} converges, the improper integral exists,
and (a) follows.  Note that interior $C^1$ regularity and continuity at the endpoint do
\emph{not} by themselves make $r_\theta\circ c$ Lebesgue integrable on $[0,\nu]$: the
candidate $\bm V_\theta=\Xi\sin(\Xi^{-2})$ (extended by $0$ at $\Xi=0$) is in
$C^1(\mathbb D)\cap C(\mathbb D\cup\Sigma^-)$, and along a characteristic its exact
residual is
$r_\theta(c(\kappa))=1+\Lambda(\kappa)\bigl[\sin(\Xi(\kappa)^{-2})-2\Xi(\kappa)^{-2}\cos(\Xi(\kappa)^{-2})\bigr]$.
Near a transversal endpoint $\Lambda$ is bounded and bounded away from $0$, and
$\Xi$ is strictly decreasing with $\Xi'=\Lambda$; substituting $x=\Xi(\kappa)$, the
leading term contributes $2\int_0^a x^{-2}\abs{\cos(x^{-2})}\,dx=\int_{a^{-2}}^\infty
u^{-1/2}\abs{\cos u}\,du=\infty$ (with $u=x^{-2}$), while the remaining terms have finite
integral; so $r_\theta\circ c\notin L^1(0,\nu)$, although the truncated integrals
converge.

(b) By Cauchy--Schwarz on $[0,\nu']$ and monotonicity of the integral,
$\abs{\int_0^{\nu'}r_\theta\circ c}^2\le\nu'\int_0^{\nu'}\abs{r_\theta\circ c}^2\le\nu\int_0^{\nu}\abs{r_\theta\circ c}^2\in[0,\infty]$
for every $\nu'<\nu$; let $\nu'\uparrow\nu$ in \eqref{E:error-char-trunc} and use
$(a+b)^2\le2a^2+2b^2$.

(c) If $\int_0^\nu\abs{r_\theta\circ c}^2<\infty$ then $r_\theta\circ c\in L^2(0,\nu)\subset L^1(0,\nu)$
and the improper integral equals the Lebesgue integral (dominated convergence).  If
$J_{\mathrm{HJB}}<\infty$, \cref{P:I2} gives
\[
\int_{\bm D}\int_0^{\bm V(z)}\abs{r_\theta(c_z(\kappa))}^2d\kappa\,dz<\infty ,
\]
so the inner integral is finite for a.e.\ $z\in\bm D$ (the integrand is Borel on
$[0,\tau_{\max}]\times\bm D$ once a Borel representative of $\sa_\mL$ is fixed; $r_\theta$
is measurable but, under $L^\infty$ forcing, not continuous, since it contains
$\sa_\mL(s)\partial_\Lambda\bm V_\theta$).  If $\bm V_\theta$ is $C^1$ on a neighbourhood of
the compact set $\cl[\bm D^+]$, then $\nabla\bm V_\theta$ is bounded there and, since
$\sa_\mL\in L^\infty$, $r_\theta$ is essentially bounded on $\bm D^+$; so
$r_\theta\circ c\in L^\infty(0,\nu)$.
\end{proof}

\begin{lemma}[Locally Lipschitz candidates]\label{L:error-char-Lip}
Let $\bm V_\theta\in C^{0,1}_{\mathrm{loc}}(\mathbb D)\cap C(\mathbb D\cup\Sigma^-)$,
and interpret $r_\theta$ using Borel representatives of the a.e. derivatives.
Then \eqref{E:error-char} and \eqref{E:pointwise-bound} hold for
$\mathcal L^3$-a.e. $z\in\bm D$. If $J_{\mathrm{HJB}}<\infty$, the characteristic
integral in \eqref{E:error-char} is an ordinary Lebesgue integral for a.e. such $z$.
These conclusions are unchanged by modifying the derivative representatives on a null set.
\end{lemma}
\begin{proof}
Let $N_\theta$ be a Borel null set containing the interior nondifferentiability points
of $\bm V_\theta$, the time-zero face, and the null set where the chosen derivatives
may differ from its differential. For fixed $\kappa\in[0,\tau_{\max}]$, put
$\bm D_\kappa=\{z\in\bm D:\kappa<\tau(z)\}$ and
$\mathcal F_\kappa(s,\bm y)=(s+\kappa,Z(\kappa;s,\bm y))$. This map is injective, and
its inverse on $\mathcal F_\kappa(\bm D_\kappa)\subset B$ is
\[
(t,\bm x)\longmapsto
\bigl(t-\kappa,Z^{\rm abs}(t-\kappa;(t,\bm x))\bigr).
\]
The inverse is Lipschitz by \eqref{E:flow-consts}: move the evaluation time along one
trajectory and then compare initial data at a common absolute time. Consequently
$\{z\in\bm D_\kappa:\mathcal F_\kappa(z)\in N_\theta\}$ is null.
The map $(\kappa,z)\mapsto\mathcal F_\kappa(z)$ is continuous on its admissible domain,
so the corresponding exceptional set in $[0,\tau_{\max}]\times\bm D$ is Borel.
Tonelli therefore shows that, for a.e. $z\in\bm D$, the characteristic
$c_z(\kappa)=\mathcal F_\kappa(z)$ avoids $N_\theta$ for a.e. $\kappa<\tau(z)$.

Fix such a $z$. On every truncated interval $[0,\nu']$, $\nu'<\tau(z)$,
$\bm V_\theta\circ c_z$ is Lipschitz. At the a.e. parameter values where
$c_z$ is differentiable and $c_z(\kappa)\notin N_\theta$, the chain rule gives
$(\bm V_\theta\circ c_z)'=r_\theta\circ c_z-1$.
Integrating this equality and taking the continuous incoming trace gives
\eqref{E:error-char}; Cauchy--Schwarz on the truncated intervals gives
\eqref{E:pointwise-bound}, exactly as in \cref{L:error-char}.
Finally \cref{P:I2}, whose proof uses only the flow and a nonnegative measurable
integrand, gives $r_\theta\circ c_z\in L^2(0,\tau(z))\subset L^1(0,\tau(z))$ for a.e.
$z$ when $J_{\mathrm{HJB}}<\infty$. The same null-set argument proves independence of
all representative choices.
\end{proof}

\begin{remark}[Fixed-signal slices of an augmented network]\label{R:rough-candidates}
For a fixed acceleration profile, let
$A(s)=\int_0^s\sa_\mL(r)\,dr$ and prescribe the leader speed as
$v_\mL(s)=w_0+A(s)$. A smooth network with current leader speed as an input induces
\[
W_\theta^{w_0}(s,X,\Lambda)
=N_\theta(s,X,\Lambda,w_0+A(s)).
\]
Since $A$ is Lipschitz, each slice is locally Lipschitz, with
$\partial_sW_\theta^{w_0}=\partial_sN_\theta+\sa_\mL(s)\partial_wN_\theta$ a.e.
If the profile has finitely many jumps and is continuous between them, the slice is
$C^1$ on the intervening time strips but need not be $C^1$ across the jumps.
Its residual is the restriction to $w=w_0+A(s)$ of
\[
\partial_sN_\theta+\Lambda\partial_XN_\theta
+\bigl[\sa_\mL(s)-\alpha(\mathcal V(X)-w+\Lambda)\bigr]\partial_\Lambda N_\theta
+\sa_\mL(s)\partial_wN_\theta+1.
\]
In particular the jumping coefficient multiplies
$(\partial_\Lambda+\partial_w)N_\theta$. The a.e. chain rule and
\cref{L:error-char-Lip} apply to these slices, provided the prescribed leader signal
satisfies the model assumptions. For finitely many forcing jumps the same identity
can also be obtained by integrating on each strip and using continuity at its endpoints.
Thus a smooth augmented network can induce the time-derivative jumps excluded by a
fixed-signal $C^1$ ansatz. Integration over a family of $w_0$ values requires a specified
parameter measure in addition to the three-dimensional losses for each slice.
\end{remark}

Integrating \eqref{E:pointwise-bound} (pointwise for $C^1$ candidates and almost
everywhere for locally Lipschitz candidates by \cref{L:error-char-Lip}) over $\bm D$,
and using $\nu\le\tau_{\max}$,
\begin{equation}\label{E:L2-split}
\norm{e}_{L^2(\bm D)}^2\ \le\ 2I_1+2\tau_{\max}I_2,
\end{equation}
where
\[
\begin{aligned}
I_1&\Def\int_{\bm D}\abs{\bm V_\theta(\Psi(z))}^2dz,\\
I_2&\Def\int_{\bm D}\int_0^{\bm V(z)}\abs{r_\theta(s+\kappa,Z(\kappa;z))}^2d\kappa\,dz .
\end{aligned}
\]

\subsection{The endpoint map: Lipschitz bound, slice estimate, Jacobian}\label{S:endpoint}

\begin{lemma}[$\tilde\Psi$ is Lipschitz on compacts]\label{L:Psi-Lip}
For $z_1,z_2\in K$,
\[
\abs{\tilde\Psi(z_1)-\tilde\Psi(z_2)}\le C_\Psi(K)\abs{z_1-z_2},\qquad
C_\Psi(K)\Def\bigl(1+\widehat L_\tau(K)\bigr)+\bigl(m_B\widehat L_\tau(K)+L_B+m_B\bigr).
\]
\end{lemma}
\begin{proof}
Write $z_i=(s_i,\bm y_i)$, $\tau_i\Def\tau(z_i)$ and $t_i\Def s_i+\tau_i\in[0,\timehor)$
(the predicted contact time; $t_i<\timehor$ because $z_i\in\mathbb D$).  First component:
$\abs{t_1-t_2}\le\abs{s_1-s_2}+\abs{\tau_1-\tau_2}\le(1+\widehat L_\tau(K))\abs{z_1-z_2}$ by
\cref{T:tau-Lipschitz} on $K$.  Second component: $\Lambda(\tau_i;z_i)=\Lambda^{\rm abs}(t_i;z_i)$
is the predicted relative velocity of $z_i$ at the \emph{absolute} time $t_i$.  We compare
in absolute time, inserting $\Lambda^{\rm abs}(t_2;z_1)$, which is defined since
$t_2\in[0,\timehor)\subset[0,\timehor+1]$ (\cref{L:flow}(a)):
\begin{align*}
\abs{\Lambda^{\rm abs}(t_1;z_1)-\Lambda^{\rm abs}(t_2;z_2)}
&\le\abs{\Lambda^{\rm abs}(t_1;z_1)-\Lambda^{\rm abs}(t_2;z_1)}+\abs{\Lambda^{\rm abs}(t_2;z_1)-\Lambda^{\rm abs}(t_2;z_2)}\\
&\le m_B\abs{t_1-t_2}+L_B\abs{z_1-z_2}
\end{align*}
by \eqref{E:flow-consts} ($z_1,z_2\in K\subset B$; all absolute times in $[0,\timehor)$).
With the first-component bound this is $\le(m_B\widehat L_\tau(K)+L_B+m_B)\abs{z_1-z_2}$.
Adding the two components gives the claim.  (Comparing at a common \emph{elapsed} time
$\tau_2$ instead would evaluate the prediction from $z_1$ at the absolute time
$s_1+\tau_2=t_2+(s_1-s_2)$, which may exceed $\timehor+1$ and leave the domain of the
flow; this is why \cref{L:elapsed} carries its domain condition.)
\end{proof}

\begin{remark}
The endpoint-map Lipschitz constant $C_\Psi(K)$ is not controlled uniformly in $K$:
through $\widehat L_\tau(K)$ it depends on the transversality margin $\gamma_K$ (the local
term $2(L_B+m_B)/\gamma_K$ of \eqref{E:tau-fd-lip}) and on the secant radius $\varrho_0(K)$
of \cref{P:tau-fd-lip}, and our bound for it is not uniform as $K$ approaches
$\partial\mathbb D$.  It is used only to guarantee that $\tilde\Psi$ is Lipschitz, hence
differentiable a.e.\ (Rademacher) and admissible in the co-area formula; it does not enter
the final constant.
\end{remark}

\begin{lemma}[The flow map at fixed start time is bi-Lipschitz]\label{L:flowmap-bilip}
Fix $s\in[0,\timehor)$ and let $\mathcal F_s(\kappa,\bm y)\Def(s+\kappa,Z(\kappa;s,\bm y))$.
For $(\kappa_i,\bm y_i)$ with $(s,\bm y_i)\in B$, $0\le\kappa_i$, $s+\kappa_i<\timehor$, and
$\mathcal F_s(\kappa_i,\bm y_i)\in B$:
\begin{equation}\label{E:flowmap-bilip}
\abs{\mathcal F_s(\kappa_1,\bm y_1)-\mathcal F_s(\kappa_2,\bm y_2)}\ \ge\ \frac{1}{1+L_B}\,\abs{(\kappa_1,\bm y_1)-(\kappa_2,\bm y_2)} .
\end{equation}
\end{lemma}
\begin{proof}
Write $p_i=(t_i,\bm x_i)\Def\mathcal F_s(\kappa_i,\bm y_i)$.  By uniqueness of the
prediction flow, $\bm y_i=Z^{\rm abs}(s;p_i)$, the state at absolute time $s$ of the
trajectory through $p_i$ (backward evaluation), and $\kappa_i=t_i-s$.  Hence, by
\eqref{E:flow-consts} with absolute time $s$,
\[
\abs{(\kappa_1,\bm y_1)-(\kappa_2,\bm y_2)}\le\abs{t_1-t_2}+\abs{Z^{\rm abs}(s;p_1)-Z^{\rm abs}(s;p_2)}
\le(1+L_B)\abs{p_1-p_2}.
\]
\end{proof}

\begin{lemma}[Slice estimate]\label{L:slice}
For $z_i=(s,\bm y_i)\in K$ with the \emph{same} start time $s$,
\begin{equation}\label{E:slice}
\abs{\tilde\Psi(s,\bm y_1)-\tilde\Psi(s,\bm y_2)}\ \ge\ c_B\,\abs{\bm y_1-\bm y_2},\qquad c_B\Def\frac1{1+L_B}.
\end{equation}
\end{lemma}
\begin{proof}
$\Psi(s,\bm y_i)=\mathcal F_s(\tau_i,\bm y_i)$ with $\tau_i=\tau(s,\bm y_i)$, and all points
lie in $\cl[K^+]\subset B$.  By \eqref{E:flowmap-bilip} and the chart isometry,
$\abs{\tilde\Psi(s,\bm y_1)-\tilde\Psi(s,\bm y_2)}=\abs{\Psi(s,\bm y_1)-\Psi(s,\bm y_2)}
\ge c_B\abs{(\tau_1,\bm y_1)-(\tau_2,\bm y_2)}\ge c_B\abs{\bm y_1-\bm y_2}$.
\end{proof}

\begin{remark}[Transversality is used, but not quantitatively]\label{R:slice-meaning}
\eqref{E:slice} says that two states at the same time with different relative states
have endpoints at least $c_B\abs{\bm y_1-\bm y_2}$ apart on the collision plane: the two
endpoints cannot coincide because the characteristics through them are distinct and do
not cross.  The margin $\gamma_{\bm D}$ does not appear; it enters only through the
Lipschitz regularity of $\tau$, which is needed to make sense of $D\tilde\Psi$ at all.
\end{remark}

\begin{proposition}[Jacobian lower bound]\label{P:J2-bound}\label{lem:J2Psi-bound}
$\tilde\Psi$ is differentiable at $\mathcal L^3$-a.e.\ point of $\operatorname{int}_{\RR^3}K\supset\bm D\setminus\{s=0\}$,
and at every such point
\begin{equation}\label{E:J2-lower}
J_2\tilde\Psi(z)\Def\sqrt{\det\bigl(D\tilde\Psi(z)\,D\tilde\Psi(z)^\top\bigr)}\ \ge\ c_\ast\Def c_B^2=\frac{1}{(1+L_B)^2}.
\end{equation}
In particular $D\tilde\Psi(z)$ has rank $2$ there.
\end{proposition}
\begin{proof}
Differentiability a.e.\ on the open set $\operatorname{int}_{\RR^3}K$ is Rademacher's theorem
(\cref{L:Psi-Lip}).  Let $z=(s,\bm y)\in\operatorname{int}_{\RR^3}K$ be a point of
differentiability and $A\Def D\tilde\Psi(z)\in\RR^{2\times3}$; write $A=[a_t\ A_{\bm y}]$
with $A_{\bm y}\in\RR^{2\times2}$ the derivative in the state directions.  For $v\in\RR^2$
and small $h$, $(s,\bm y+hv)\in K$, so \eqref{E:slice} gives
$\abs{\tilde\Psi(s,\bm y+hv)-\tilde\Psi(s,\bm y)}\ge c_Bh\abs v$; dividing by $h$ and letting
$h\to0$, $\abs{A_{\bm y}v}\ge c_B\abs v$.  Hence both singular values of $A_{\bm y}$ are
$\ge c_B$ and $\abs{\det A_{\bm y}}\ge c_B^2$.  By the Cauchy--Binet formula,
$\det(AA^\top)=\sum_{S}\det(A_S)^2$, the sum over the three $2\times2$ minors $A_S$ of $A$;
keeping only $S=\{\text{state columns}\}$, $\det(AA^\top)\ge(\det A_{\bm y})^2\ge c_B^4$.
\end{proof}

\begin{lemma}[Fibres are arcs of characteristics]\label{L:fibre}
Let $p=(t^\ast,\Lambda^\ast)\in\RR^2$ and $E_p\Def\tilde\Psi^{-1}(p)\cap\bm D$.  Set
$r_p\Def\min\{t^\ast,\tau_{\max}(\bm D)\}$ if $E_p\ne\varnothing$ and $r_p\Def0$
otherwise.  If $E_p\ne\varnothing$, then $t^\ast\in(0,\timehor)$, and with
\[
\gamma_p:[0,r_p]\to[0,\timehor)\times\RR^2,\qquad
\gamma_p(r)\Def\bigl(t^\ast-r,\ Z^{\rm abs}(t^\ast-r;(t^\ast,0,\Lambda^\ast))\bigr)
\]
(the backward prediction from the collision point) we have $E_p\subset\gamma_p([0,r_p])$,
and $\gamma_p$ is Lipschitz with constant $\sqrt{1+m_B^2}$.  Consequently, for every
$p\in\RR^2$,
\begin{equation}\label{E:fibre-length}
\mathcal H^1\bigl(\tilde\Psi^{-1}(p)\cap\bm D\bigr)\ \le\ \sqrt{1+m_B^2}\;r_p\ \le\
\mL_{\max}\Def\sqrt{1+m_B^2}\ \tau_{\max}(\bm D).
\end{equation}
\end{lemma}
\begin{proof}
Let $z=(s,\bm y)\in E_p$.  Then $t^\ast=s+\tau(z)$ with $0<\tau(z)\le\tau_{\max}(\bm D)$,
and $t^\ast=\sigma(z)<\timehor$ since $z\in\mathbb D$; so $t^\ast\in(0,\timehor)$ and
$\tau(z)=t^\ast-s\le\min\{t^\ast,\tau_{\max}(\bm D)\}=r_p$.  The prediction from $z$ is at
the state $(0,\Lambda^\ast)$ at absolute time $t^\ast$; by uniqueness of the global
prediction flow, forward and backward (\cref{L:flow}(a)), it coincides on $[0,\timehor+1]$
with the prediction from the datum $q_p\Def(t^\ast,0,\Lambda^\ast)$.  Hence
$\bm y=Z^{\rm abs}(s;q_p)$, i.e.\ $z=\gamma_p(\tau(z))$ with $\tau(z)\in[0,r_p]$.  This
proves $E_p\subset\gamma_p([0,r_p])$.

The datum $q_p=\Psi(z)$ lies in $\cl[\bm D^+]\subset B$, and for $r\in[0,r_p]$ the absolute
time $t^\ast-r$ lies in $[0,t^\ast]\subset[0,\timehor+1]$.  Since \eqref{E:flow-consts} holds
for data in $B$ at all absolute times in $[0,\timehor+1]$, whether or not the trajectory
remains in $B$ (it is proved from the a priori bound $\abs{Z}\le R_B$ along the whole
trajectory, \cref{L:flow}(b)),
\[
\abs{\gamma_p(r_1)-\gamma_p(r_2)}^2=\abs{r_1-r_2}^2+\abs{Z^{\rm abs}(t^\ast-r_1;q_p)-Z^{\rm abs}(t^\ast-r_2;q_p)}^2
\le(1+m_B^2)\abs{r_1-r_2}^2 .
\]
A Lipschitz image of an interval of length $r_p$ has $\mathcal H^1$-measure at most
$\mathrm{Lip}(\gamma_p)\,r_p$, which gives \eqref{E:fibre-length}; if $E_p=\varnothing$ the
bound is trivial.  (The parametrisation must stop at $r_p$: for $r>t^\ast$ the absolute
time $t^\ast-r$ is negative and the prediction is not defined.)
\end{proof}

\subsection{The two bounds and the theorem}\label{S:two-bounds}

\begin{proposition}[Boundary term]\label{P:I1}
\[
\begin{aligned}
I_1&\le\frac{\mL_{\max}}{c_\ast}\,J_{\partial\mathcal Q}\\
&=(1+L_B)^2\sqrt{1+m_B^2}\ \tau_{\max}(\bm D)\ J_{\partial\mathcal Q}.
\end{aligned}
\]
\end{proposition}
\begin{proof}
Let $f:\RR^3\to\RR^2$ be a Lipschitz extension of $\tilde\Psi|_K$ (McShane, componentwise);
$f=\tilde\Psi$ on $\operatorname{int}_{\RR^3}K\supset\bm D\setminus\{s=0\}$, so $Df=D\tilde\Psi$ and $J_2f=J_2\tilde\Psi$
at every differentiability point of $\operatorname{int}_{\RR^3}K$.  Let $E\subset\bm D$ be the
null set where $\tilde\Psi$ is not differentiable, and set
$g\Def\mathbf 1_{\bm D\setminus E}\,\abs{\bm V_\theta\circ\Psi}^2/J_2\tilde\Psi\ge0$, measurable
and well defined by \eqref{E:J2-lower}.  The co-area formula
\cite[Thm.~3.2.12]{Federer}, \cite[Thm.~3.11]{EvansGariepy} for Lipschitz
$f:\RR^3\to\RR^2$ gives
\[
\begin{aligned}
I_1&=\int_{\bm D\setminus E}\abs{\bm V_\theta(\Psi(z))}^2dz
=\int_{\RR^3}g\,J_2f\,dz\\
&=\int_{\RR^2}\int_{f^{-1}(p)}g\,d\mathcal H^1\,dp\\
&=\int_{\RR^2}\abs{\bm V_\theta(p)}^2
\int_{\tilde\Psi^{-1}(p)\cap(\bm D\setminus E)}\frac{d\mathcal H^1}{J_2\tilde\Psi}\,dp ,
\end{aligned}
\]
where $\bm V_\theta(p)$ denotes the value at the point of $\partial\mathcal Q$ with chart
coordinates $p$ (constant on the fibre).  By \eqref{E:J2-lower} and \eqref{E:fibre-length}
the inner integral is $\le\mL_{\max}/c_\ast$, and the outer integral runs over
$\tilde\Psi(\bm D)$, which in the chart is $\Psi(\bm D)\subset\mathcal Q_{\bm D}$ with
$dp=d\mathcal H^2$.  Hence $I_1\le(\mL_{\max}/c_\ast)\int_{\mathcal Q_{\bm D}}\abs{\bm V_\theta}^2d\mathcal H^2$.
\end{proof}

\begin{proposition}[Interior term]\label{P:I2}
\[
\begin{aligned}
I_2&\le\frac{e^{\alpha\tau_{\max}(\bm D)}-1}{\alpha}\,J_{\mathrm{HJB}}\\
&\le\tau_{\max}(\bm D)\,e^{\alpha\tau_{\max}(\bm D)}\,J_{\mathrm{HJB}}.
\end{aligned}
\]
\end{proposition}
\begin{proof}
By Tonelli, $I_2=\int_0^{\tau_{\max}}\int_{\bm D_\kappa}\abs{r_\theta(\mathcal F_\kappa(z))}^2dz\,d\kappa$
with $\bm D_\kappa\Def\{z\in\bm D:\tau(z)>\kappa\}$ and
$\mathcal F_\kappa(s,\bm y)\Def(s+\kappa,Z(\kappa;s,\bm y))$, which is defined on
$\bm D_\kappa$ since there $s+\kappa<s+\tau(z)<\timehor$.  For fixed
$\kappa\in[0,\tau_{\max}]$, $\mathcal F_\kappa$ is injective (uniqueness of the flow) and
Lipschitz on $\bm D_\kappa$ with constant $\le1+L_B+m_B$ (\cref{L:elapsed}, whose domain
condition holds); its inverse $w=(t,\bm x)\mapsto(t-\kappa,Z^{\rm abs}(t-\kappa;w))$ is
Lipschitz with constant $\le1+L_B+m_B$ on $\mathcal F_\kappa(\bm D_\kappa)\subset\bm D^+\subset B$:
for $w_1,w_2$ there, the start times are $s_i=t_i-\kappa$ and
$\abs{Z^{\rm abs}(s_1;w_1)-Z^{\rm abs}(s_2;w_2)}\le m_B\abs{s_1-s_2}+L_B\abs{w_1-w_2}$ by
\eqref{E:flow-consts} (first move the evaluation time along trajectory~1, then compare
the two trajectories at the common absolute time $s_2$).

\emph{The Jacobian of $\mathcal F_\kappa$ is explicit.}  At every point of
differentiability of $\mathcal F_\kappa$ (a.e.\ on $\bm D_\kappa$, by Rademacher), its
Jacobian matrix is block lower-triangular,
\[
D\mathcal F_\kappa(s,\bm y)=\begin{pmatrix}1&0\\ \partial_sZ(\kappa;s,\bm y)&Y(\kappa;s,\bm y)\end{pmatrix},
\]
with $Y=\partial_{\bm y}Z$ the state linearisation of \cref{lem:clean-linearization}, which
exists everywhere.  Since $\kappa\mapsto Y(\kappa;s,\bm y)$ is absolutely continuous and
solves \eqref{E:var-clean} with $\operatorname{tr}\partial_{\bm y}\mathsf G\equiv-\alpha$
(see \eqref{E:dG-dy}: the trace is independent of the state and of the forcing), Jacobi's
formula $(\det Y)'=\operatorname{tr}(\partial_{\bm y}\mathsf G)\det Y$ a.e.\ gives
Liouville's identity
\begin{equation}\label{E:Liouville}
\det Y(\kappa;s,\bm y)=e^{-\alpha\kappa}\qquad\text{for all }(s,\bm y)\text{ and }\kappa .
\end{equation}
Hence $J\mathcal F_\kappa=\abs{\det D\mathcal F_\kappa}=e^{-\alpha\kappa}$ a.e.\ on
$\bm D_\kappa$, and the area formula for the injective Lipschitz map $\mathcal F_\kappa$
\cite[Thm.~3.9]{EvansGariepy} gives
\[
\int_{\bm D_\kappa}\abs{r_\theta\circ\mathcal F_\kappa}^2dz
=e^{\alpha\kappa}\int_{\bm D_\kappa}\abs{r_\theta\circ\mathcal F_\kappa}^2\,J\mathcal F_\kappa\,dz
=e^{\alpha\kappa}\int_{\mathcal F_\kappa(\bm D_\kappa)}\abs{r_\theta(w)}^2\,dw
\le e^{\alpha\kappa}J_{\mathrm{HJB}} .
\]
Integrating in $\kappa$ over $[0,\tau_{\max}]$ gives the claim.  (Bounding instead the
Jacobian of the Lipschitz inverse by $(1+L_B+m_B)^3$ gives the cruder
$I_2\le\tau_{\max}(1+L_B+m_B)^3J_{\mathrm{HJB}}$, which is also valid; \eqref{E:Liouville}
removes the flow constants from the interior term altogether.)
\end{proof}

\begin{theorem}[Continuum residual stability]\label{T:two-sided}\label{thm:two-sided-convergence}
Let $\bm D\Subset\mathbb D$ be compact and let
$\bm V_\theta\in C^{0,1}_{\mathrm{loc}}(\mathbb D)\cap C(\mathbb D\cup\Sigma^-)$.
The residual is interpreted almost everywhere as in \eqref{E:loss}.
Then
\begin{equation}\label{E:two-sided}
\norm{\bm V_\theta-\bm V}_{L^2(\bm D)}^2\ \le\ C_{\bm D}^2\bigl(J_{\mathrm{HJB}}+J_{\partial\mathcal Q}\bigr),
\end{equation}
with
\[
C_{\bm D}^2\Def2\,\tau_{\max}(\bm D)\,\max\Bigl\{(1+L_B)^2\sqrt{1+m_B^2},\ \frac{e^{\alpha\tau_{\max}(\bm D)}-1}{\alpha}\Bigr\}.
\]
The interior coefficient depends only on $\alpha$ and $\tau_{\max}(\bm D)\le\timehor$; the
boundary coefficient depends on $\bm D$ only through the bounded set $B\supset\cl[K^+]$
(through $m_B,L_B$).  In particular the constant does \emph{not} depend on the
transversality margin $\gamma_{\bm D}$ and remains bounded as $\bm D$ approaches
$\partial\mathbb D$ within a family admitting a common bounded reference set.
\end{theorem}
\begin{proof}
\eqref{E:L2-split} with \cref{P:I1,P:I2}; if $J_{\mathrm{HJB}}=\infty$ there is nothing to
prove.  The estimate is one-directional: it bounds the error by the loss.  The reverse
inequality is false in general (perturb $\bm V$ on a smooth interior patch by
$n^{-1/2}\chi(z)\sin(ns)$ with $\chi$ smooth and compactly supported: the $L^2$ error tends
to zero while the residual grows like $n^{1/2}$), which is why we do not call the estimate
``two-sided''.
\end{proof}

\begin{remark}[Where $\bm D\Subset\mathbb D$ is used]\label{R:where-compact}
In the co-area proof, the positive distance of $\bm D$ to $\partial\mathbb D$ (the margin
$\gamma_K$ and the secant radius $\varrho_0(K)$ of a compact neighbourhood $K$) is used in
\cref{L:Psi-Lip}: $\tilde\Psi$ must be Lipschitz on a neighbourhood of $\bm D$ for
Rademacher and the co-area formula.  The Lipschitz constant $C_\Psi(K)$ there is bounded
only in terms of $\gamma_K$ and $\varrho_0(K)$ and may be large, but it never enters
\eqref{E:two-sided}.  Boundedness of $\bm D$ enters only through $m_B,L_B$, and since
\eqref{E:flow-consts} for a bounded set implies it for every subset, $m_B,L_B$ may be
replaced by the constants of any bounded $B'\supset B$.  Hence, if
$\bm D_\infty\subset\mathbb D$ is bounded and open (or bounded and $\sigma$-compact, so
that it admits a compact exhaustion $\bm D_n\uparrow\bm D_\infty$; an arbitrary bounded
subset need not even be measurable) and $B_\infty$ is a bounded set containing the
closed $1$-neighbourhood of the set reached on $[0,\timehor)$ by predictions from the
closed $1$-neighbourhood of $\bm D_\infty$ (bounded by the Gr\"onwall bound in
\cref{L:flow}(b)), then \eqref{E:two-sided} holds for every compact
$\bm D\subset\bm D_\infty$ with the single constant obtained from
$m_{B_\infty},L_{B_\infty}$ and $\tau_{\max}\le\timehor$; since
$\bm D_\infty^+=\bigcup_n\bm D_n^+$ and $\Psi(\bm D_\infty)=\bigcup_n\Psi(\bm D_n)$, monotone
convergence along the exhaustion gives the estimate on $\bm D_\infty$ itself, with the
losses taken over $\bm D_\infty^+$ and $\Psi(\bm D_\infty)$.  Thus the transversality margin
governs the regularity of $\tau$ (\cref{T:tau-Lipschitz}) and the constant $C_\Psi$, but
not the stability estimate.
\end{remark}

\subsection{An explicit estimate from endpoint and age coordinates}\label{S:flux-stability}

The co-area argument separates the rank of the endpoint map from the length of its
fibres. A direct change of variables retains their exact volume factor and removes the
flow Lipschitz constants from the boundary coefficient as well.

\begin{proposition}[Boundary pullback by terminal endpoint and age]\label{P:flux-boundary}
Let $\bm D\Subset\mathbb D$ be nonempty and compact. Set
\[
T=\tau_{\max}(\bm D),\qquad
\widetilde{\mathcal Q}_{\bm D}=\tilde\Psi(\bm D),\qquad
M_*=\sup_{(t^*,\Lambda_*)\in\widetilde{\mathcal Q}_{\bm D}}|\Lambda_*|,
\qquad A_T=\frac{e^{\alpha T}-1}{\alpha}.
\]
Define the compact set of actual first-contact endpoint/age pairs by
\[
E_{\bm D}=\{(\tilde\Psi(z),\tau(z)):z\in\bm D\}.
\]
For every nonnegative Borel function $f$ on $\bm D$,
\begin{equation}\label{E:flux-area}
\int_{\bm D}f(z)\,dz
=\int_{E_{\bm D}}f\bigl(H(p,\kappa)\bigr)|\Lambda_*|e^{\alpha\kappa}\,dp\,d\kappa,
\end{equation}
where $p=(t^*,\Lambda_*)$ and
\begin{equation}\label{E:flux-map}
H(p,\kappa)=
\bigl(t^*-\kappa,Z^{\rm abs}(t^*-\kappa;(t^*,0,\Lambda_*))\bigr).
\end{equation}
In particular, writing $b(p)=\bm V_\theta(t^*,0,\Lambda_*)$,
\begin{equation}\label{E:flux-boundary}
I_1=\int_{E_{\bm D}}|b(p)|^2|\Lambda_*|e^{\alpha\kappa}\,dp\,d\kappa
\le M_*A_T J_{\partial\mathcal Q}.
\end{equation}
\end{proposition}
\begin{proof}
Continuity of $\tilde\Psi$ and $\tau$ gives compactness of $E_{\bm D}$. Every pair in
this set satisfies $0<\kappa\le t^*<\timehor$, so the backward evaluation in
\eqref{E:flux-map} is admissible. The absolute-time flow bounds in \cref{L:flow}(b)
show that $H$ is Lipschitz on bounded admissible parameter sets, including a relative
neighbourhood of $E_{\bm D}$. By uniqueness of the prediction flow,
$H(\tilde\Psi(z),\tau(z))=z$. Conversely, an initial point determines its first-contact
endpoint and flight time uniquely, so $H|_{E_{\bm D}}$ is injective and maps onto $\bm D$.
No injectivity is asserted on the full endpoint/age rectangle, which can contain
recrossings of the collision plane.

We compute the a.e. Jacobian. Put $q_p=(t^*,0,\Lambda_*)$, $s=t^*-\kappa$,
$\bm y=Z^{\rm abs}(s;q_p)$, and
$Y_{\rm b}=\partial_{\bm\eta}Z^{\rm abs}(s;(t^*,\bm\eta))|_{\bm\eta=(0,\Lambda_*)}$.
The backward state derivative exists by the same variational argument as in
\cref{lem:clean-linearization}; equivalently it is the inverse of the forward state
derivative from $(s,\bm y)$ to time $t^*$. Hence \eqref{E:Liouville} gives
$\det Y_{\rm b}=e^{\alpha\kappa}$.
Let $N$ be the common null set of non-Lebesgue times of $\sa_\mL$.
At Lebesgue times the initial-time derivative of the flow follows from the cocycle
identity and the startup expansion used in \cref{lem:V-C1}(ii): at fixed end time,
$\partial_{t^*}Z^{\rm abs}(s;q_p)=-Y_{\rm b}\mathsf G(q_p)$.
The evaluation-time derivative is $\mathsf G(s,\bm y)$ when $s\notin N$.
Thus, at a.e. differentiability point of the Lipschitz map $H$ with
$t^*\notin N$ and $t^*-\kappa\notin N$, its three columns are
\begin{align*}
\partial_{t^*}H&=(1,\mathsf G(s,\bm y)-Y_{\rm b}\mathsf G(q_p)),\\
\partial_{\Lambda_*}H&=(0,Y_{\rm b}\bm e_2),\\
\partial_\kappa H&=(-1,-\mathsf G(s,\bm y)).
\end{align*}
The exceptional sets $t^*\in N$ and $t^*-\kappa\in N$ are product-null by Fubini.
The faces $s=0$ and $\kappa=0$ are also null. Adding the first and third columns,
and using that the first component of $\mathsf G(q_p)$ is $\Lambda_*$, gives
\begin{equation}\label{E:flux-jacobian}
|\det DH(p,\kappa)|=|\Lambda_*|\det Y_{\rm b}
=|\Lambda_*|e^{\alpha\kappa}\quad\text{a.e.}
\end{equation}
One may extend $H$ componentwise from a bounded admissible parameter neighbourhood
before applying the area formula; its differential on the interior is unchanged, and
the parameter boundary faces and their Lipschitz images are null. The area formula for
the injective restriction to $E_{\bm D}$ \cite[Thm.~3.9]{EvansGariepy} now yields
\eqref{E:flux-area}.

For $z=H(p,\kappa)$ with $(p,\kappa)\in E_{\bm D}$, the first-contact endpoint is $p$.
Applying \eqref{E:flux-area} to $f(z)=|\bm V_\theta(\Psi(z))|^2$ proves the equality in
\eqref{E:flux-boundary}. Each age fibre of $E_{\bm D}$ is a measurable subset of
$[0,T]$, so Tonelli and $|\Lambda_*|\le M_*$ bound this integral by
$M_*\int_0^T e^{\alpha\kappa}d\kappa\int_{\widetilde{\mathcal Q}_{\bm D}}|b(p)|^2dp$.
The endpoint chart is an isometry, so the last integral is
$J_{\partial\mathcal Q}$, which proves the inequality.
\end{proof}

\begin{theorem}[Residual stability with an explicit flux coefficient]\label{T:flux-stability}
Under the candidate assumptions of \cref{T:two-sided}, and with $T,M_*,A_T$ as in
\cref{P:flux-boundary},
\begin{equation}\label{E:flux-stability}
\norm{\bm V_\theta-\bm V}_{L^2(\bm D)}^2
\le A_T(M_*+T)\bigl(J_{\partial\mathcal Q}+J_{\mathrm{HJB}}\bigr).
\end{equation}
More generally, if $\kappa_0\le\tau(z)\le T$ on $\bm D$ for some $\kappa_0\ge0$, then
\begin{equation}\label{E:flux-stability-age}
\norm{\bm V_\theta-\bm V}_{L^2(\bm D)}^2
\le\left[
M_*\frac{e^{\alpha T}-e^{\alpha\kappa_0}}{\alpha}+TA_T
\right]\bigl(J_{\partial\mathcal Q}+J_{\mathrm{HJB}}\bigr).
\end{equation}
\end{theorem}
\begin{proof}
If $J_{\mathrm{HJB}}=\infty$ the claim is immediate. Otherwise
\cref{L:error-char-Lip,P:I2} give the ordinary error identity for a.e. $z\in\bm D$.
Minkowski, followed by Cauchy--Schwarz on each characteristic, yields
\[
\norm{e}_{L^2(\bm D)}
\le\sqrt{I_1}+\sqrt{TI_2}
\le\sqrt{M_*A_T}\sqrt{J_{\partial\mathcal Q}}
+\sqrt{TA_T}\sqrt{J_{\mathrm{HJB}}}.
\]
Cauchy--Schwarz in $\RR^2$ proves \eqref{E:flux-stability}.
If the lower age bound holds, the fibres of $E_{\bm D}$ lie in $[\kappa_0,T]$, so
\eqref{E:flux-boundary} instead gives
$I_1\le M_*(e^{\alpha T}-e^{\alpha\kappa_0})J_{\partial\mathcal Q}/\alpha$.
The same two inequalities then prove \eqref{E:flux-stability-age}.
\end{proof}

\begin{remark}[Uniformity and the role of grazing]\label{R:flux-uniform}
For physical initial states, \cref{L:prediction-bounds}(a) gives
$M_*\le v_{\max}$. For arbitrary bounded analytical initial families, the prediction
bounds likewise control $M_*$ uniformly. Thus \eqref{E:flux-stability} extends to a
bounded $\sigma$-compact evaluation region by the compact exhaustion in
\cref{R:where-compact}, with $T\le\timehor$ and a uniform upper bound for $M_*$.
The normal flux $|\Lambda_*|$ in \eqref{E:flux-jacobian}, rather than its reciprocal,
is the reason the integrated estimate remains bounded while the pointwise hitting-time
gradient can grow like the inverse transversality margin. The flux estimate does not
change the need for the correct incoming boundary trace or turn sampled losses into
continuum guarantees.
\end{remark}

\subsection{The HJB equation: what is proved and what is cited}\label{S:HJB-senses}

\begin{proposition}[Senses in which $\bm V$ solves \eqref{E:HJB}]\label{P:HJB-senses}\label{rem:HJB-regularity}
\begin{enumerate}[label=\tup{\alph*}]
\item \textbf{(Characteristic form.)}  $\bm V$ satisfies \eqref{E:DPP-V} on $\mathbb D$;
  this is the integrated form of \eqref{E:HJB} along characteristics and holds exactly.
\item \textbf{(Almost everywhere classical.)}  $\bm V$ is Lipschitz on compacts of
  $\mathbb D$, Fr\'echet-differentiable at every point where $\partial_s\bm V$ exists, and
  \eqref{E:HJB} holds at every $\bm y$ for every $s$ outside the null set $N$ of
  non-Lebesgue points of $\sa_\mL$, hence a.e.\ on $\mathbb D$
  (\cref{lem:V-C1}\eqref{it:Lip-time}); the exceptional set of times is $N$ and does not
  depend on $\bm y$.
\item \textbf{(Classical where the forcing is continuous.)}  On
  $\mathbb D\cap(J\times\RR^2)$, $J$ an interval of continuity of $\sa_\mL$, $\bm V\in C^1$
  and \eqref{E:HJB} holds everywhere (\cref{lem:V-C1}\eqref{it:C1-joint}).
\item \textbf{(Uniqueness along characteristics.)}
  If $W\in C(\mathbb D\cup\Sigma^-)$ is absolutely continuous on compact subintervals of
  every predicted characteristic in $\mathbb D$ with
  $\tfrac{d}{d\kappa}W(s+\kappa,Z(\kappa;z))=-1$ a.e.\ in $\kappa$, and $W=0$ on
  $\Sigma^-$, then $W=\bm V$.  (Integrate from the collision point; the derivative is
  bounded, so no endpoint-integrability issue arises.)
\item \textbf{(Uniqueness among locally Lipschitz a.e.\ solutions.)}  If
  $W\in C(\mathbb D\cup\Sigma^-)$ is locally Lipschitz on $\mathbb D$, satisfies
  $\partial_sW+\mathsf G\cdot\nabla_{\bm y}W=-1$ a.e.\ on $\mathbb D$ and $W=0$ on
  $\Sigma^-$, then $W=\bm V$. This follows from \cref{P:zero-residual}, whose candidate
  class is locally Lipschitz, by exhausting $\mathbb D$ with compact sets
  $\bm D=\cl[\operatorname{int}\bm D]$.
\end{enumerate}
\end{proposition}

\begin{remark}[On the viscosity characterisation]\label{R:viscosity}
One expects $\bm V$ to be the unique viscosity solution of \eqref{E:HJB} on $\mathbb D$,
by the standard argument from \eqref{E:DPP-V}; we do not prove this, and nothing in this
paper depends on it.  Two points deserve care.  The Hamiltonian
$H(s,\bm y,p)=p\cdot\mathsf G(s,\bm y)+1$ (convention $\partial_s\bm V+H=0$) is only
measurable in $s$ (through $\sa_\mL$), so the classical theory
\cite{crandall1992user,bardi1997optimal} does not apply directly and the relevant framework is the
$L^1$-in-time theory of Ishii \cite{Ishii85} and Nunziante \cite{Nunziante90}; moreover a
pointwise definition would depend on the representative of $\sa_\mL$, whereas the flow,
$\tau$ and the losses do not.  And the datum is prescribed only on the incoming surface
$\Sigma^-$, so a comparison principle requires that characteristics exit $\mathbb D$ only
through $\Sigma^-$, which is \cref{L:D-open}(b) with \cref{L:V-cont-Sigma}.  When
$\sa_\mL$ is continuous, $\bm V\in C^1(\mathbb D)$ by \cref{lem:V-C1}\eqref{it:C1-joint}
and is then a classical, hence viscosity, solution; we make no claim in the measurable
case.  The stability estimate uses only \eqref{E:DPP-V}, and \cref{P:zero-residual}
proves $\bm V_\theta=\bm V$ directly.

The structural reason is (d),(e) above: $H$ is affine in $p$, the characteristics are
the predicted trajectories and do not cross, so a solution that is absolutely continuous
along characteristics (for instance a locally Lipschitz a.e.\ solution) is already
determined by its values on $\Sigma^-$.  Absolute continuity cannot be dropped: if $C$ is
a nonconstant Cantor function with $C(0)=0$ and $C'=0$ off a null set $N_C$, then
$W\Def\tau+C(\tau)$ on $\mathbb D$, $W\Def0$ on $\Sigma^-$, is continuous on
$\mathbb D\cup\Sigma^-$, vanishes on $\Sigma^-$ and satisfies \eqref{E:HJB} a.e.\ on
$\mathbb D$ (where $\tau$ is differentiable and $\tau(z)\notin N_C$,
$DW=(1+C'(\tau))D\tau=D\tau$; and $\tau^{-1}(N_C)\cap\bm D$ is $\mathcal L^3$-null for every
compact $\bm D$ by the co-area formula, because $\abs{\nabla_{\bm y}\tau}\ge
e^{-G_1\tau_{\max}}/m_B>0$ there by \eqref{eq:grad-tau} and \eqref{E:Y-two-sided}), but
$W\ne\bm V$.
\end{remark}

\begin{definition}[Strong fixed point]\label{D:strong-fp}
$W\in C(\mathbb D\cup\Sigma^-)$ is a \emph{strong fixed point} on $\mathbb D$ if for every
$(s,\bm y)\in\mathbb D$ and every $\kappa\in[0,\bm V(s,\bm y)]$,
\begin{equation}\label{E:strong-fp}
W\bigl(s+\kappa,Z(\kappa;s,\bm y)\bigr)=W(s,\bm y)-\kappa .
\end{equation}
\end{definition}
By \eqref{E:DPP-V} and
\cref{L:V-cont-Sigma}, $\bm V$ is a strong fixed point; by \cref{P:HJB-senses}(d) it is
the only one vanishing on $\Sigma^-$.  Continuity is required on $\mathbb D\cup\Sigma^-$
only: the prescribed zero extension need not be continuous on the opening part of the
collision plane, and extension at grazing depends on the accessible contact branches
(\cref{R:no-global-ext}).

\subsection{Zero total loss implies \texorpdfstring{$\bm V_\theta=\bm V$}{V\_theta = V}}\label{S:zero-residual}\label{ss:zero-residual}

\begin{proposition}[Exact fixed point at zero total loss]\label{P:zero-residual}\label{prop:zero-residual-fp}
Let $\bm D\Subset\mathbb D$ be compact with $\bm D=\cl[\operatorname{int}\bm D]$, and let
$\bm V_\theta\in C^{0,1}_{\mathrm{loc}}(\mathbb D)\cap C(\mathbb D\cup\Sigma^-)$ satisfy
$J_{\mathrm{HJB}}=0$ and
$J_{\partial\mathcal Q}=0$ (losses over $\bm D^+$ and $\mathcal Q_{\bm D}$ as in \eqref{E:loss}).
Then
\begin{equation}\label{E:char-identity}
\bm V_\theta\bigl(\mathcal F_\kappa(z)\bigr)=\tau(z)-\kappa\qquad\text{for every }z\in\bm D
\text{ and every }\kappa\in[0,\tau(z)],
\end{equation}
where $\mathcal F_\kappa(s,\bm y)\Def(s+\kappa,Z(\kappa;s,\bm y))$.  In particular
$\bm V_\theta=\bm V$ on $\bm D^+\cup\mathcal Q_{\bm D}\supset\bm D$, and $\bm V_\theta$
satisfies \eqref{E:strong-fp} for every $(s,\bm y)\in\bm D$.
\end{proposition}

\begin{proof}
Throughout, $B$, $m_B$, $L_B$ are the bounded reference set and the flow constants fixed
in \cref{S:HJB-setup} (``Bounded reference set and constants''), so that
$\cl[\bm D^+]\subset\cl[K^+]\subset B$; $\bm D_\kappa\Def\{z\in\bm D:\tau(z)>\kappa\}$, and
$\Psi(z)=\mathcal F_{\tau(z)}(z)\in\mathcal Q_{\bm D}\subset\Sigma^-$ (\cref{L:V-cont-Sigma}).

\emph{Step 0: the space--time flow map is bi-Lipschitz.}  On
$A\Def\{(\kappa,z)=(\kappa,s,\bm y):\ z\in B,\ 0\le\kappa\le\tau_{\max},\ s+\kappa<\timehor\}$
let $\Phi(\kappa,z)\Def(\kappa,\mathcal F_\kappa(z))\in\RR^4$.  By \eqref{E:flow-consts},
$\Phi$ is Lipschitz on $A$; it is injective by uniqueness of the prediction flow, and its
inverse $\Phi^{-1}(\kappa,t,\bm x)=\bigl(\kappa,\,t-\kappa,\,Z^{\rm abs}(t-\kappa;(t,\bm x))\bigr)$
is Lipschitz on $\Phi(A)\cap(\RR\times B)$ with constant $\le1+2m_B+L_B$ (as in the proof
of \cref{P:I2}: move the evaluation time along one trajectory, then compare trajectories
at a common absolute time).  Hence, on the target region $\Phi(A)\cap(\RR\times B)$,
$\Phi^{-1}$ maps $\mathcal L^4$-null sets to $\mathcal L^4$-null sets, and preimages under
$\Phi$ of Lebesgue measurable subsets of that region are Lebesgue measurable (a measurable
set is a Borel set up to a null set, and $\Phi$ is continuous).  This is all that Step~1
needs, since the null set used there lies in $[0,\tau_{\max}]\times\bm D^+\subset\RR\times B$.

\emph{Step 1: the residual vanishes along a.e.\ characteristic.}  The residual $r_\theta$
has a Borel representative on $\mathbb D$ by the convention following \eqref{E:loss}.
Let $N_{\rm diff}$ be a Borel null set containing the nondifferentiability points of
$\bm V_\theta$ in the interior, the time-zero face, and any points where the chosen derivative representatives differ from its actual differential. The candidate derivatives are
only a.e. defined, and the forcing can be merely measurable; continuity of the residual
is not required.  The set
$\bm D^+$ is Borel (the image of the $\sigma$-compact set
$\{(\kappa,z):z\in\bm D,\ 0\le\kappa<\tau(z)\}$ under the continuous map
$(\kappa,z)\mapsto\mathcal F_\kappa(z)$). Let
$N\Def\{w\in\bm D^+:r_\theta(w)\ne0\}\cup(N_{\rm diff}\cap\bm D^+)$;
$\mathcal L^3(N)=0$ because $J_{\mathrm{HJB}}=0$ and $N_{\rm diff}$ is null.  Set
\begin{align*}
\tilde N&\Def\{(\kappa,z)\in[0,\tau_{\max}]\times\bm D:\ \kappa<\tau(z),\ \mathcal F_\kappa(z)\in N\}\\
&=\Phi^{-1}\bigl([0,\tau_{\max}]\times N\bigr)\cap\{(\kappa,z):z\in\bm D,\ \kappa<\tau(z)\}\subset A .
\end{align*}
Since $[0,\tau_{\max}]\times N$ is $\mathcal L^4$-null, Step~0 shows that $\tilde N$ is
measurable and $\mathcal L^4(\tilde N)=0$.  By Fubini, for $\mathcal L^3$-a.e.\ $z\in\bm D$
the slice $\tilde N_z=\{\kappa\in[0,\tau(z)):\mathcal F_\kappa(z)\in N\}$ is
$\mathcal L^1$-null.  Let $E_1\subset\bm D$ be the set of such $z$;
$\mathcal L^3(\bm D\setminus E_1)=0$.

\emph{Step 2: integration along such a characteristic.}  Let $z=(s,\bm y)\in E_1$ and
$c(\kappa)\Def\mathcal F_\kappa(z)$, $\kappa\in[0,\tau(z)]$.  The curve $c$ is Lipschitz,
lies in $\bm D^+\subset\mathbb D$ for $\kappa<\tau(z)$, and $c(\tau(z))=\Psi(z)\in\Sigma^-$.
Since $\bm V_\theta$ is locally Lipschitz, $g\Def\bm V_\theta\circ c$ is locally
Lipschitz on $[0,\tau(z))$. By Step~1, $c(\kappa)$ avoids $N_{\rm diff}$ for a.e.
$\kappa$, so the chain rule at differentiability points gives
$g'(\kappa)=\partial_s\bm V_\theta(c(\kappa))+\nabla_{\bm y}\bm V_\theta(c(\kappa))\cdot\mathsf G(c(\kappa))
=r_\theta(c(\kappa))-1$ for a.e.\ $\kappa$, since
$\partial_\kappa Z=\mathsf G(s+\kappa,Z)$ a.e.  By Step~1, $g'=-1$ a.e.\ on
$[0,\tau(z))$, so $g(\kappa)=g(0)-\kappa$ for $0\le\kappa<\tau(z)$.  Letting
$\kappa\uparrow\tau(z)$ and using the continuity of $\bm V_\theta$ on
$\mathbb D\cup\Sigma^-$ at $\Psi(z)$,
\begin{equation}\label{E:V-theta-along}
\bm V_\theta\bigl(\mathcal F_\kappa(z)\bigr)=\bm V_\theta\bigl(\Psi(z)\bigr)+\tau(z)-\kappa
\qquad\text{for all }\kappa\in[0,\tau(z)],\ z\in E_1 .
\end{equation}

\emph{Step 3: the boundary value vanishes at a.e.\ endpoint.}  Let
$N_\partial\Def\{p\in\mathcal Q_{\bm D}:\bm V_\theta(p)\ne0\}$, relatively open in the
compact set $\mathcal Q_{\bm D}$ (continuity of $\bm V_\theta$ on $\Sigma^-$), hence Borel;
$\mathcal H^2(N_\partial)=0$ because $J_{\partial\mathcal Q}=0$.  The set
$\Psi^{-1}(N_\partial)\cap\bm D$ is Borel ($\Psi$ is continuous on $\mathbb D$).  Apply the
co-area formula as in the proof of \cref{P:I1}, with
$g\Def\mathbf 1_{(\Psi^{-1}(N_\partial)\cap\bm D)\setminus E}$, $E$ the null set of
non-differentiability of $\tilde\Psi$:
\[
\int_{\Psi^{-1}(N_\partial)\cap\bm D}J_2\tilde\Psi\,dz
=\int_{\tilde\Psi(\Psi^{-1}(N_\partial))}\mathcal H^1\bigl(\tilde\Psi^{-1}(p)\cap\bm D\bigr)\,dp=0 ,
\]
since the outer domain of integration is the chart image of $N_\partial$, which is
$\mathcal L^2$-null.  As $J_2\tilde\Psi\ge c_\ast>0$ a.e.\ on $\bm D$ (\cref{P:J2-bound}),
$\mathcal L^3(\Psi^{-1}(N_\partial)\cap\bm D)=0$.  Put
$E\Def E_1\setminus\Psi^{-1}(N_\partial)$.  Then $\mathcal L^3(\bm D\setminus E)=0$ and, by
\eqref{E:V-theta-along},
\begin{equation}\label{E:char-identity-ae}
\bm V_\theta\bigl(\mathcal F_\kappa(z)\bigr)=\tau(z)-\kappa\qquad\text{for all }\kappa\in[0,\tau(z)],\ z\in E .
\end{equation}

\emph{Step 4: from almost every to every point of $\bm D$, by density.}  $E$ is dense in
$\bm D$: if $U\subset\RR^3$ is open with $U\cap\bm D\ne\varnothing$, then
$U\cap\operatorname{int}\bm D\ne\varnothing$ because $\bm D=\cl[\operatorname{int}\bm D]$,
and this nonempty (relatively) open set has positive measure, so it meets $E$.  Fix
$z\in\bm D$ and choose $z_n\in E$ with $z_n\to z$.

Let $0\le\kappa<\tau(z)$.  Since $\tau$ is continuous on $\mathbb D$
(\cref{L:unif-transv}(a)), $\kappa<\tau(z_n)$ for all large $n$, so
$\bm V_\theta(\mathcal F_\kappa(z_n))=\tau(z_n)-\kappa$ by \eqref{E:char-identity-ae}.  As
$n\to\infty$, $\mathcal F_\kappa(z_n)\to\mathcal F_\kappa(z)\in\bm D^+\subset\mathbb D$
(\cref{L:flow}(b)) and $\bm V_\theta$ is continuous on $\mathbb D$; hence
$\bm V_\theta(\mathcal F_\kappa(z))=\tau(z)-\kappa$.

Let $\kappa=\tau(z)$.  By joint continuity of $(\kappa,z)\mapsto\mathcal F_\kappa(z)$ and of
$\tau$,
\[
\Psi(z_n)=\mathcal F_{\tau(z_n)}(z_n)\to\mathcal F_{\tau(z)}(z)=\Psi(z).
\]
All these points lie in $\mathcal Q_{\bm D}\subset\Sigma^-$, $\bm V_\theta(\Psi(z_n))=0$ by
\eqref{E:char-identity-ae} at $\kappa=\tau(z_n)$, and $\bm V_\theta$ is continuous on
$\mathbb D\cup\Sigma^-$; hence $\bm V_\theta(\Psi(z))=0=\tau(z)-\tau(z)$.  This proves
\eqref{E:char-identity}.

\emph{Step 5: conclusions.}  Taking $\kappa=0$ in \eqref{E:char-identity} gives
$\bm V_\theta=\tau=\bm V$ on $\bm D$.  Every point of $\bm D^+$ is $\mathcal F_\kappa(z)$ with
$z\in\bm D$, $\kappa<\tau(z)$, and
$\tau(z)-\kappa=\tau(\mathcal F_\kappa(z))=\bm V(\mathcal F_\kappa(z))$ by \eqref{E:DPP}; every
point of $\mathcal Q_{\bm D}$ is some $\Psi(z)$, where $\bm V=0$.  Thus $\bm V_\theta=\bm V$ on
$\bm D^+\cup\mathcal Q_{\bm D}$.  Finally, for $(s,\bm y)\in\bm D$ and
$\kappa\in[0,\bm V(s,\bm y)]$, \eqref{E:char-identity} reads
$\bm V_\theta(s+\kappa,Z(\kappa;s,\bm y))=\bm V_\theta(s,\bm y)-\kappa$, which is
\eqref{E:strong-fp}.
\end{proof}

\begin{remark}[What the proof does not use]\label{R:zero-residual-noneed}
The proof uses neither continuity of $r_\theta$ (which fails under $L^\infty$ forcing),
nor continuity of $\bm V_\theta$ beyond $\mathbb D\cup\Sigma^-$, nor any claim that the
characteristics through almost every point of $\bm D$ fill almost all of $\bm D^+$.  The
latter claim is false in general: a null
transversal subset of $\bm D$ can sweep out positive volume under the flow (in flow-box
coordinates $(p,r)\in\RR^2\times\RR$, take $\bm D=\{\abs r\le\operatorname{dist}(p,C)\}$ with
$C$ compact of positive area and empty interior; then $\bm D=\cl[\operatorname{int}\bm D]$,
$C\times\{0\}$ is $\mathcal L^3$-null, but its forward saturation contains
$C\times(0,a)$, of positive measure).  This is why Step~4 argues by density.  The hypothesis
$\bm D=\cl[\operatorname{int}\bm D]$ enters only in Step~4.
\end{remark}

\begin{remark}[Zero total loss and a smooth-class obstruction]\label{R:zero-residual-scope}
The locally Lipschitz class of \cref{P:zero-residual} includes the exact value function
even for $L^\infty$ forcing. A restriction to $C^1(\mathbb D)$ candidates can, however,
prevent attainment of zero total loss. Suppose the forcing has a genuine jump at
$t_1$ with the one-sided regularity in \cref{R:not-C1-in-time}, and some
$(t_1,\bm y)\in\operatorname{int}(\bm D^+)$ satisfies
$\partial_\Lambda\bm V(t_1,\bm y)\ne0$. Then the time derivative of $\bm V$ has a jump
there. Since zero total loss forces equality on $\bm D^+$ by \cref{P:zero-residual},
no $C^1$ candidate can attain it. The sufficient positivity conditions in
\cref{R:not-C1-in-time} supply examples of this obstruction.

This concerns the combination of interior residual and incoming boundary error.
The $C^1$ function $W(s,\bm y)=\timehor-s$ has zero interior residual for every forcing,
but fails the incoming zero datum. Moreover, nonattainment of zero total loss does not
imply a positive infimum or an unavoidable $L^2$ error. We do not prove an approximation
rate or characterize that infimum for a specified neural-network class. In particular,
this obstruction does not predict a numerical bias under rough forcing: the fixed-signal
slices in \cref{R:rough-candidates} are generally locally Lipschitz rather than $C^1$,
and can have time-derivative jumps. The stability estimates assert that a sequence of
admissible candidates converges in $L^2(\bm D)$ whenever its continuum total loss tends
to zero; they do not assert that a training procedure produces such a sequence.
\end{remark}

\section{Numerical experiments}
\label{sec:numerics}

The analysis of \Cref{sec:ttc-regularity,sec:HJ_collision} has two numerical counterparts.
\Cref{S:num-learning} studies a learned map and the relation between its continuum residual
and approximation error. The principal checkpoint was trained with an additional
horizon-value penalty on the dynamic boundary, beyond the loss \eqref{eq:J-total}; this
penalty conflicts with the hitting-side trace on the grazing part. We retain this experiment
as a diagnostic of boundary-error transport and compare three smaller-network objectives
with their actual loss weights stated below. On a separately defined compact characteristic tube we estimate the
continuum losses and compare their ratio to the analytical bound. \Cref{S:num-regularity}
checks gradient scaling and forcing-jump formulas with a reference solver, and distinguishes
variational-flow sensitivity from sensitivity at actual first contact. The closed loop driven
by a learned map is not simulated, and no property of such a closed loop is claimed. The
repository is identified in the code-availability statement.

\subsection{Learning the TTC map}\label{S:num-learning}

\paragraph{Setting}
We take $\alpha=2$, $\timehor=1$ and the optimal velocity function
\eqref{E:compact-two-dynamics-maintion}, whose supremum is $v_{\max}=1+\tanh2$. Because the
leader is prescribed, the TTC of a state at time $s$ depends on $(s,\Xi,\Lambda)$, on the
acceleration profile, and on the leader speed $v_\mL$ at time $s$, which fixes the leader
trajectory on the prediction window. A single network
$N_\theta(s,\Xi,\Lambda,w;\sa_\mL)$, whose fourth input $w$ is current leader speed,
therefore represents a family of value functions. To identify a fixed-signal problem, put
$A(s)=\int_0^s\sa_\mL(r)\,dr$ and fix an initial leader speed $w_0$. The corresponding slice is
\begin{equation}\label{E:num-fixed-signal-slice}
W_\theta^{w_0}(s,\Xi,\Lambda)
=N_\theta(s,\Xi,\Lambda,w_0+A(s);\sa_\mL),\qquad w_0=w-A(s).
\end{equation}
We also write $\bm V_\theta$ for this slice when the signal is fixed. The thirteen profiles of
\Cref{tab:profiles} are encoded by a one-hot vector. For continuous acceleration the slice
is $C^1$; for the piecewise-constant profile it is generally only locally Lipschitz and
piecewise $C^1$, even though the four-input network is smooth. Its a.e. chain rule is
\[
\partial_s W_\theta^{w_0}
=\partial_sN_\theta+\sa_\mL(s)\partial_wN_\theta,
\]
so the augmented residual used in training is
$r_\theta=\partial_sN_\theta+\mathsf G\cdot\nabla_{\bm y}N_\theta
+\sa_\mL(s)\partial_wN_\theta+1$.
This agrees with \eqref{E:loss} on the fixed-signal slice; its locally Lipschitz interpretation
is covered by \Cref{L:error-char-Lip,R:rough-candidates}. In particular a forcing jump
multiplies $(\partial_\Lambda+\partial_w)N_\theta$ in the augmented residual. On $5{,}243$
interior states of the accelerating and braking profiles at distance at least $0.01$ from
the boundary layers, the reported root-mean-square augmented residual is $4\times10^{-3}$;
omitting the moving-input term gives $0.32$.

The leader speed is clipped to $[0,v_{\max}]$ during training. Where the nominal speed would
leave this range, the effective signal differs from the prescribed profile and
\eqref{E:accel_ell} does not describe the clipped trajectory. On $59$ states where the two
reference conventions differ by more than $0.02$, the checkpoint is closer to the clipped
value in every case. The independent accuracy evaluation below uses the nominal prescribed flow with an
exact check of the future leader-speed bounds on $[s,1]$. The fresh residual diagnostic uses
a different sampler and includes nominal speed-bound violations, as specified below.
Pooled statistics over $w$ involve a family of fixed-signal slices. The compact-tube
comparison fixes $w_0$ and satisfies the physical speed bounds on the whole horizon.

\begin{table}[htbp!]
\centering\small
\caption{Leader-acceleration profiles used for training and evaluation, $t\in[0,1]$.}
\label{tab:profiles}
\begin{tabular}{ll}
\toprule
Profile & Leader acceleration $\sa_\mL(t)$\\
\midrule
Constant velocity & $0$\\
Pulsed & $3.5\,e^{-(\ln100)t}\sin(8\pi t)$\\
Constant acceleration & $\pm0.1,\ \pm0.5$\\
Fast acceleration & $0.7,\ 1,\ 1.2$\\
Hard braking & $-0.7,\ -1,\ -1.2$\\
Piecewise constant & $0.5,\,-1.2,\,0.8,\,-0.4,\,0$ on $[0,0.2),\dots,[0.8,1]$\\
\bottomrule
\end{tabular}
\end{table}

\paragraph{Network and training}
The normalized input $(s,\Xi,\Lambda,v_\mL)$ and the one-hot profile pass through a
$256$-unit $\tanh$ input layer and six hidden layers of widths $128$, $256$, $512$, $256$, $128$,
$64$ with activations $\tanh$, swish, $\tanh$, $\tanh$, swish, $\tanh$ and residual connections.
The scalar output is $\sigma(u_\theta)\max\{\timehor-s,0\}$, with $\sigma$ the logistic
function, so that $0\le\bm V_\theta\le\timehor-s$. The model has $710{,}529$ trainable
parameters. The residual is computed by automatic differentiation. Profiles are balanced at the proposal stage; rejection changes their proportions in the
accepted interior batch. Training uses an approximate boundary $\Xi_\partial(s,\Lambda,v_\mL)$ between contact and no contact,
located by a gap line search. Interior collocation points are drawn with
$10^{-5}<\Xi<\Xi_{\partial}-10^{-5}$ ($131{,}072$ per step, $s$ uniform, $(\Lambda,v_\mL)$
in the feasible triangle $\{-v_{\max}\le\Lambda\le0,\ 0\le v_\mL\le v_{\max}+\Lambda\}$,
then $\Xi$ uniform on the admissible interval). The boundary batch ($16{,}384$ per step) consists
of points on the incoming surface $\Sigma^-=\{\Xi=0,\ -v_{\max}\le\Lambda<0\}$ with target
$0$, and of points on the dynamic boundary $\Xi=\Xi_\partial$ with target $\timehor-s$. The
second kind of point supplies an additional penalty absent from \eqref{eq:J-total}. Its target
is correct on the horizon part of the boundary. On a grazing boundary whose miss side has no
later contact before the horizon, the hitting-side trace is the grazing flight time, smaller
than $\timehor-s$, so that target is incorrect (\Cref{R:no-global-ext}). The smaller comparison below
uses different boundary targets and incoming-boundary weights. The weights of the interior and the combined boundary
batches are $10$ and $1$, respectively; these batches are empirical averages under their
sampling laws, not the continuum integrals of \eqref{E:loss}. The reported optimizer is Adam
with gradient clipping at norm one and a peak learning rate of $3\times10^{-4}$, halved after
eight validation checks without improvement down to $10^{-6}$; the network is evaluated at the
$500{,}000$-step checkpoint. The configuration specifies validation and test sets of
$4{,}096$ states each; the final log reports test RMSE $0.0413$ and $R^2=0.9818$.
The principal checkpoint predates the boundary-objective comparison below; retraining the
full network with the incoming-boundary datum alone is the planned next step.

\paragraph{Reference values}
All reference values below are computed with the solver of \Cref{S:num-regularity} (order-eight
Runge--Kutta with the first contact located on monotone pieces of the predicted gap), which reduces missed shallow near-grazing contacts by detecting extrema of the gap. The
finite-sample checks below do not establish a universal event-detection guarantee. The reported
training-boundary construction uses fourth-order Runge--Kutta with event detection at $\Xi=0$.
On $254$ random lines $(s,\Lambda,v_\mL)$ across the twelve continuous profiles, interpolation
of the stored boundary cache agrees with the boundary found by
bisection on the reference solver to a root-mean-square difference of $1.3\times10^{-3}$ and a
maximal difference of $8.9\times10^{-3}$ and a small sample mean difference
$6\times10^{-5}$. This comparison assesses the stored boundary table on those sampled lines.

\paragraph{Accuracy on a bounded evaluation region}
For each profile, the reported evaluation draws use $s\in[0,0.9]$, $\Xi\in[10^{-3},1]$,
$\Lambda\in[-v_{\max},0]$, a feasible current leader speed with no clipping on the prediction
window, and rejection of states outside the transversal domain. For a fixed $w_0$, write
\begin{align*}
\bm D_{0.05}(w_0)&\Def
\{z\in\mathbb D^{w_0}:0\le s\le0.9,\ 10^{-3}\le\Xi\le1,\\
&\hspace{39mm}-v_{\max}\le\Lambda\le0,\ \gamma(z)\ge0.05\},\\
\gamma(z)&\Def-\Lambda(\tau(z);z),
\end{align*}
with the physical/no-clipping restrictions of the evaluation sampler understood.
We abbreviate this family of evaluation regions by $\bm D_{0.05}$. The lower bound on
$\gamma$ does not make these sets compactly contained in $\mathbb D$: contacts may approach the terminal horizon. On each fixed-signal slice they are
bounded $\sigma$-compact regions: the box, margin, and physical-speed inequalities
are closed conditions within the open transversal domain. Thus \Cref{R:where-compact}
applies by compact exhaustion with uniform bounded-set constants. The pooled metrics in
\Cref{tab:profilemetrics}, however, are errors under the evaluation sampler over current
leader speeds. Conditional uniform sampling of $w$ on a state-dependent interval is not
joint Lebesgue-uniform sampling without a density correction. We therefore do not identify
these pooled RMSEs with normalized three-dimensional $L^2$ norms. The explicit continuum
comparison below uses a fixed signal and its characteristic-volume measure.

\begin{figure}[htbp!]
\centering
\includegraphics[width=\textwidth]{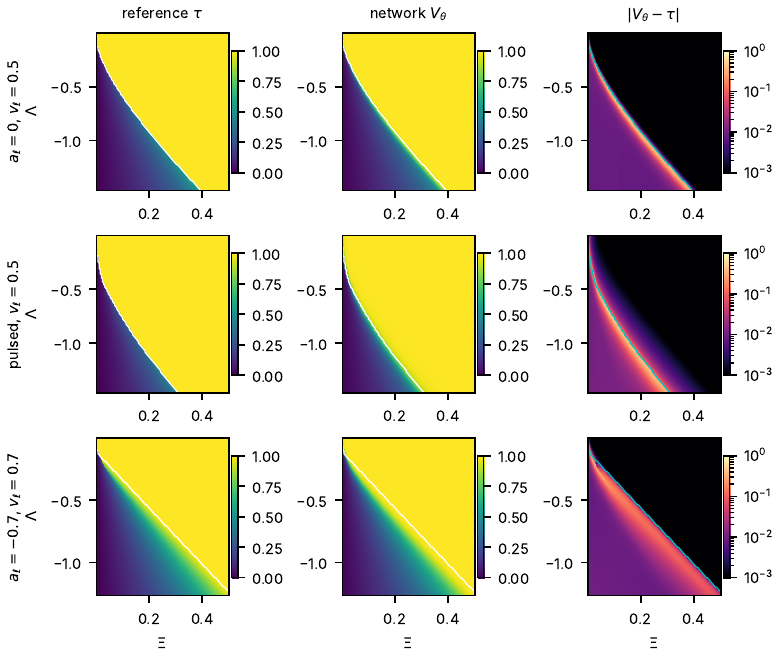}
\caption{Reference TTC (left), network (middle) and absolute error (right, logarithmic) at $s=0$
on the closing half-plane, for the constant-velocity profile with $v_\mL=0.5$ (top), the pulsed
profile with $v_\mL=0.5$ (middle) and $\sa_\mL=-0.7$ with $v_\mL=0.7$ (bottom); $\Lambda$ is
restricted to follower speeds in $[0,v_{\max}]$. The white (cyan) curve is the boundary of the
transversal domain; above it $\tau=\timehor-s=1$.}
\label{fig:maps}
\end{figure}

The saved sample has $1{,}500$ states per profile after the margin filter. Exact integration
of the prescribed leader signal detects two rough-profile states with slightly negative
future leader speed; removing them leaves $1{,}498$ states in that row and changes its
RMSE from $0.188468$ to $0.188314$. All continuous-profile samples pass the same check.
Proposal-retention fractions are $11\%$--$24\%$; these depend on the conditional sampler
and are not volume fractions. The fixed-speed sections in \Cref{fig:maps} display a thin
contact wedge. On $\bm D_{0.05}$ the RMSE is $0.034$--$0.068$ for the
twelve continuous-acceleration profiles and the mean absolute error is $0.017$--$0.043$.
The error increases as the margin decreases: the RMSE is $0.008$--$0.018$ for $\gamma\ge0.8$,
$0.03$--$0.14$ for $0.2\le\gamma<0.4$, and $0.17$--$0.45$ for
$0.05\le\gamma<0.1$. This is consistent with the worsening regularity bound in
\Cref{T:tau-Lipschitz}, although that theorem does not predict the error of a particular
trained network. The pulsed and fast-acceleration profiles have the largest errors in this sample. Values of $R^2$ are low when the TTC variance is small; for fast
acceleration the mean TTC on this sample is approximately $0.13$. They are not directly comparable with the pooled log-reported test $R^2$.

For the piecewise-constant profile the RMSE is $0.188$, with mean overprediction about $0.16$,
and the dependence on the margin is weak. The time derivative of the exact fixed-signal map
jumps at the four forcing jumps (\Cref{R:not-C1-in-time}), but
\Cref{R:zero-residual-scope} does not predict a positive loss infimum or this bias. Moreover,
\eqref{E:num-fixed-signal-slice} need not be $C^1$ across those jumps, so the smooth fixed-signal
nonattainment statement does not directly apply to this architecture. The observed error may
reflect the boundary target, sampling, optimization, approximation capacity, or clipping in
training; these experiments do not isolate its cause. The ranges quoted for the twelve continuous profiles exclude this row. The small-margin
RMSEs use only $7$--$26$ states per continuous profile and $12$ rough-profile states; they
are descriptive sample statistics rather than precise population estimates.

\begin{table}[htbp!]
\centering\small
\caption{Independent evaluation of the $500{,}000$-step checkpoint. The columns use
$1{,}500$ reference-solver states per continuous profile and $1{,}498$ for the piecewise
profile after exact future-speed postfiltering: RMSE, mean absolute error, $R^2$, mean TTC,
and RMSE in two margin bands. These are statistics under the conditional sampler, not
normalized continuum norms. The last column is incoming-boundary RMS (true value $0$), with $s$ uniform on
$[0,1-10^{-4}]$ and $(\Lambda,w)$ area-uniform on
$\{-v_{\max}\le\Lambda\le0,\ 0\le w\le v_{\max}+\Lambda\}$.}
\label{tab:profilemetrics}
\footnotesize\setlength{\tabcolsep}{3pt}
\begin{tabular}{lrrrrrrr}
\toprule
Profile & RMSE & MAE & $R^2$ & mean $\tau$ & $\gamma<.1$ & $\gamma\ge.8$ & $\Sigma^-$\\
\midrule
Constant $0$        & .039 & .017 & .95 & .20 & .24 & .009 & .13\\
Pulsed              & .068 & .043 & .83 & .20 & .21 & .018 & .15\\
Constant $+0.1$     & .044 & .019 & .93 & .20 & .27 & .009 & .13\\
Constant $-0.1$     & .044 & .020 & .93 & .20 & .33 & .009 & .13\\
Constant $+0.5$     & .046 & .018 & .89 & .16 & .28 & .008 & .14\\
Constant $-0.5$     & .040 & .020 & .95 & .21 & .27 & .011 & .12\\
Brake $-0.7$   & .040 & .019 & .94 & .19 & .27 & .010 & .11\\
Brake $-1.0$   & .037 & .019 & .95 & .19 & .21 & .011 & .10\\
Brake $-1.2$   & .034 & .020 & .96 & .19 & .17 & .013 & .09\\
Fast $+0.7$ & .054 & .020 & .79 & .14 & .29 & .008 & .15\\
Fast $+1.0$ & .063 & .021 & .66 & .13 & .45 & .008 & .15\\
Fast $+1.2$ & .057 & .021 & .67 & .13 & .38 & .008 & .16\\
Piecewise           & .188 & .159 & $<0$ & .19 & .19 & .17 & .18\\
\bottomrule
\end{tabular}
\end{table}

\paragraph{Where the error is}
\Cref{fig:maps} shows, for three profiles and one fixed leader speed, the reference map, the
network, and their difference at $s=0$ on the closing half-plane. The transversal domain is the
wedge below the white curve (the dynamic boundary $\Xi_\partial$); above it the prediction does
not reach the collision set and $\tau=\timehor-s$, which the network reproduces to
$0.006$--$0.018$. Inside the wedge the error is concentrated in two layers: along the grazing
boundary, where $\tau$ jumps from the grazing time to $\timehor-s$ and the continuous network
interpolates across the jump over a width of about $0.01$--$0.03$ in $\Xi$, and along the
collision surface $\Xi=0$ at small closing speeds. The second layer is the boundary loss: the
root-mean-square value of $\bm V_\theta$ on $\Sigma^-$ is $0.09$--$0.16$ for the twelve
continuous-acceleration profiles (\Cref{tab:profilemetrics}), and it is concentrated at small $\abs\Lambda$ (pooled over the
profiles under a separate sampler uniform in $\Lambda$ and conditionally uniform in
leader speed, $0.03$ for $\abs\Lambda>1$, $0.05$ for $0.3<\abs\Lambda<1$, $0.20$ for
$0.05<\abs\Lambda<0.3$ and $0.40$ for $\abs\Lambda<0.05$), that is, where $\Sigma^-$ meets the
grazing set and $\tau$ is discontinuous nearby. A fresh residual diagnostic uses $131{,}072$ proposals and retains $130{,}734$ points
across the twelve continuous profiles. Its unweighted mean squared residual is
$3.4799\times10^{-3}$; the union of the two gap layers of width $0.01$ and the
small-closing-speed layer $\abs\Lambda<0.05$ carries $99.60\%$ of that sampled loss,
and the mean outside this union is $2.56\times10^{-5}$. Here $(\Lambda,w)$ is uniform on
the feasible triangle, but the gap is conditionally uniform on
$[m,\Xi_\partial-m]$, $m=10^{-5}$. The joint density therefore contains
$(\Xi_\partial-2m)^{-1}$: this is an empirical residual mean, not a Lebesgue-volume loss.
Moreover $19{,}556$ retained states fail the nominal future leader-speed bounds, so the
diagnostic is not an evaluation confined to physical fixed-signal domains. Its discrepancy
from the logged training-batch mean $9.37\times10^{-5}$ does not by itself identify a
sampling mechanism or estimate the continuum loss. The separately defined tube below
uses a fixed physical signal and the correct characteristic-volume weights.

\paragraph{Continuum comparison on a compact characteristic tube}
We now fix a signal, rather than pooling current leader speeds. For the constant-velocity and
pulsed profiles take $w_0=0.6$; for $\sa_\mL=-1.2$ take $w_0=1.25$. These prescribed leader
speeds remain in $[0,v_{\max}]$ on $[0,1]$. Let
$Q=[0.35,0.95]\times[-1,-0.3]$ parameterize incoming points
$p=(t^*,0,\Lambda_*)$. For $0.04\le\kappa\le0.25$, let $z(p,\kappa)$ be the prediction state
at time $t^*-\kappa$ obtained by integrating backward from $p$. Retain the pairs for which
the follower speed is physical and the gap is positive on the segment before $p$, and denote
their image by $\bm D$.

Here the retained tube is compactly contained in $\mathbb D$; this uses more than forward
speed invariance. To verify it, write $r$ for backward age and $v(r)$ for follower speed.
On a retained physical segment, $0\le\Xi(r)\le v_{\max}r\le v_{\max}/4$ and therefore
\begin{align*}
d&\Def\mathcal V(v_{\max}/4)\approx0.0572673775,\\
v'(r)&=2\bigl(v(r)-\mathcal V(\Xi(r))\bigr),
\qquad v(r)\ge d+(v(0)-d)e^{2r}.
\end{align*}
The endpoint follower speed is at least $0.9$ for the constant and pulsed signals and at least
$0.41$ for braking. For the pulse, direct integration gives
$0\le A(t)\le0.2106$ on $[0,1]$ (the maximum is approximately $0.2105049201$).
Consequently, up to the physical-speed exit,
\begin{align*}
\Lambda(r)&\le-0.3 &&\text{for the constant signal},\\
\Lambda(r)&\le-0.0894 &&\text{for the pulse},\\
\Lambda(r)&\le-0.3+1.2r-(0.41-d)(e^{2r}-1)\le-0.2288
&&\text{for braking},\qquad 0\le r\le0.25.
\end{align*}
Thus $\Xi(r)\ge0.0894r$, and every retained initial point has
$\Xi\ge0.00357$. The backward follower speed is strictly increasing while it is physical,
so its exit at $v_{\max}$ is transverse. Continuous dependence, together with the uniform
positive-gap bound, shows that the retained parameter set is closed in
$Q\times[0.04,0.25]$. Its image is compact, its contact is first and transversal with
$|\Lambda_*|\ge0.3$, and $t^*\le0.95<\timehor$. In particular,
\begin{align*}
\bm D&=\{z(p,\kappa):p\in Q,\ 0.04\le\kappa\le\kappa_{\max}(p)\},\\
\mathcal Q_{\bm D}&=\{p\in Q:\kappa_{\max}(p)\ge0.04\},
\end{align*}
where $\kappa_{\max}$ is the physical-speed exit age capped at $0.25$, and
$\tau(z(p,\kappa))=\kappa$. The forward tube $\bm D^+$ uses the same endpoints with
$0<\kappa\le\kappa_{\max}(p)$.

By \Cref{P:flux-boundary}, the parameterization has Jacobian
$|\Lambda_*|e^{\alpha\kappa}$. Hence the squared error and interior loss are estimated by
Monte Carlo integration on the appropriate retained parts of $Q\times[0,0.25]$ using this
weight, while the boundary loss uses the endpoint-plane measure. The reported estimates use $3{,}000$ endpoint proposals, with rejected endpoints
contributing zero and conditional age-interval lengths included in the volume weights;
the exact TTC of a retained state is its age $\kappa$. The implementation approximates
the physical exit age from below with bisection tolerance $0.002$, so the integrals refer
to its numerical tube. Monte Carlo standard errors do not include this domain truncation
or ODE error.
The same endpoint/age argument gives an explicit analytical coefficient for this tube.
Since $|\Lambda_*|\le1$ and the boundary pullback on $\bm D$ uses ages at least $0.04$,
\[
\|\bm V_\theta\circ\Psi\|_{L^2(\bm D)}^2
\le M_{\rm tube}J_{\partial\mathcal Q},\qquad
M_{\rm tube}\Def\frac{e^{0.5}-e^{0.08}}{2}.
\]
The residual integral has squared norm at most
$0.25A_{0.25}J_{\mathrm{HJB}}$, where $A_{0.25}=(e^{0.5}-1)/2$.
Minkowski's inequality followed by Cauchy--Schwarz in $\RR^2$ therefore yields
\begin{equation}\label{E:num-tube-coefficient}
\begin{aligned}
\|\bm V_\theta-\bm V\|_{L^2(\bm D)}^2
&\le\bigl(M_{\rm tube}+0.25A_{0.25}\bigr)
\bigl(J_{\partial\mathcal Q}+J_{\mathrm{HJB}}\bigr)\\
&\le0.3638072604\bigl(J_{\partial\mathcal Q}+J_{\mathrm{HJB}}\bigr),
\end{aligned}
\end{equation}
where the final decimal is rounded upward. This is the age-restricted version of
\Cref{T:flux-stability}, specifically \eqref{E:flux-stability-age}; it requires no
estimated flow suprema.

The ratios $0.10$--$0.14$ in \Cref{tab:exp1} are approximately $2.6$--$3.6$ times below
the explicit coefficient $0.3638072604$. Since the integrals are estimated,
the comparison is empirical; a certified numerical inequality would also require validated
quadrature or confidence bounds for the integrals.

\begin{figure}[htbp!]
\centering
\includegraphics[width=\textwidth]{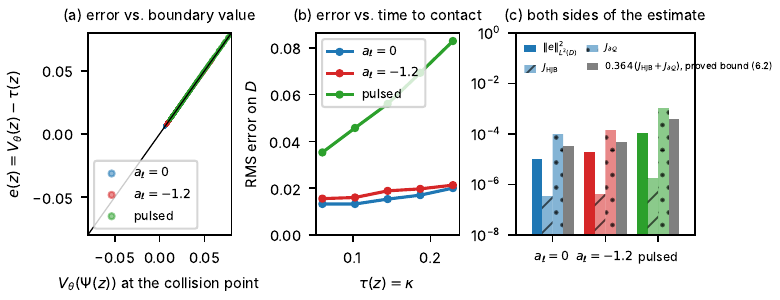}
\caption{The estimate on the backward tube $\bm D$. (a) Error $e(z)$ against the network value
at the collision point $\Psi(z)$ of $z$, for the samples of $\bm D$: the error closely
tracks the transported boundary value. (b) Sampled age-band RMS against time to contact $\tau(z)=\kappa$ under the
endpoint/conditional-age sampler. (c) Estimated squared error, interior and boundary losses, and the bound obtained
by multiplying the estimated loss sum by the proved coefficient $0.3638072604$
(displayed as $0.364$ in the legend);
this coefficient is analytical, while the integrals are numerical estimates.}
\label{fig:exp1}
\end{figure}

\begin{figure}[htbp!]
\centering
\includegraphics[width=\textwidth]{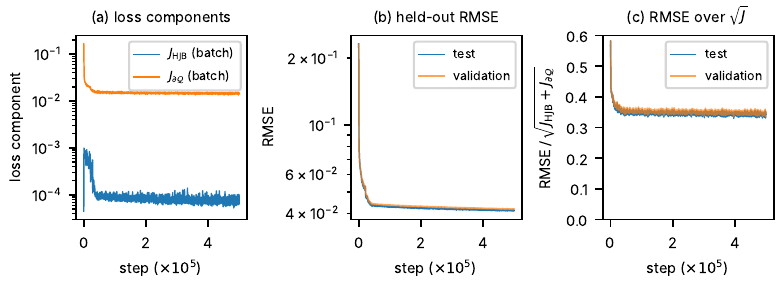}
\caption{Training history. (a) The two loss components (training-batch means). (b) RMSE on the
validation and test sets. (c) Ratio of RMSE to the square root of the summed logged batch
loss components. The loss symbols in this plot refer to sampled training quantities, including
the extra dynamic-boundary target, rather than the continuum integrals in \eqref{E:loss}.}
\label{fig:training}
\end{figure}

The boundary term dominates the measured loss by a factor $10^2$--$10^3$. The scatter of
$e(z)$ against $\bm V_\theta(\Psi(z))$ lies close to the diagonal
(\Cref{fig:exp1}(a)), consistent with the characteristic identity \eqref{E:error-char} and a
small residual integral. The factor near $0.1$ has a geometric interpretation:
\[
\int_{\bm D}\bm V_\theta(\Psi(z))^2\,dz
=\int_{\mathcal Q_{\bm D}}\bm V_\theta(p)^2
\int_{0.04}^{\kappa_{\max}(p)}|\Lambda_*|e^{\alpha\kappa}\,d\kappa\,dp
\approx0.1J_{\partial\mathcal Q}.
\]
The sampled RMS difference between $e(z)$ and $\bm V_\theta(\Psi(z))$ is below
$10^{-3}$ in every age band; this is not a pointwise bound. Thus the growth with age in \Cref{fig:exp1}(b) is consistent with selection
of endpoints having longer admissible backward tubes and larger boundary errors, rather than
large accumulation of the residual. During training, the ratio of test RMSE to the square root
of the summed logged batch losses stays approximately between $0.33$ and $0.36$ from step $5\times10^4$
onward (\Cref{fig:training}(c)). This is an empirical correlation under the training and test
samplers; it does not establish an optimization or network-size convergence rate and is not a
direct evaluation of the continuum estimate on $\bm D$.

\begin{table}[htbp!]
\centering\footnotesize
\caption{Estimated error and losses on the numerical backward tube associated with
$Q=[0.35,0.95]\times[-1,-0.3]$, $\kappa\in[0.04,0.25]$ ($3{,}000$ endpoint proposals;
Monte Carlo standard errors below $5\%$; estimated volumes $0.053$, $0.060$, $0.039$
in row order). The initial leader speed is $0.6$ (constant, pulsed) and $1.25$ (braking).
Here $\norm e^2=\norm e^2_{L^2(\bm D)}$ and the ratio is
$\norm e^2/(J_{\mathrm{HJB}}+J_{\partial\mathcal Q})$. The proved upper coefficient
in \eqref{E:num-tube-coefficient} is $0.3638072604$ in every row; numerical integration
and domain-truncation errors are separate from Monte Carlo error.}
\label{tab:exp1}
\setlength{\tabcolsep}{3pt}
\begin{tabular}{lrrrrr}
\toprule
Profile & RMSE & $\norm{e}^2$ & $J_{\mathrm{HJB}}$ & $J_{\partial\mathcal Q}$ & ratio\\
\midrule
$\sa_\mL=0$    & $0.014$ & $1.0\times10^{-5}$ & $3.3\times10^{-7}$ & $9.4\times10^{-5}$ & $0.11$\\
$\sa_\mL=-1.2$ & $0.018$ & $1.9\times10^{-5}$ & $4.1\times10^{-7}$ & $1.3\times10^{-4}$ & $0.14$\\
Pulsed          & $0.052$ & $1.1\times10^{-4}$ & $1.7\times10^{-6}$ & $1.0\times10^{-3}$ & $0.10$\\
\bottomrule
\end{tabular}
\end{table}

\paragraph{Comparison of boundary objectives}
For the constant-velocity profile with $v_\mL=1$, a smaller $4\times128$ $\tanh$ network
uses $20{,}000$ Adam steps, $4{,}096$ interior points and $1{,}024$ points per boundary
term per step, with two seeds for each objective. Variant (A) reuses its incoming batch
for both boundary terms. Write $L_r,L_-,L_{\rm open},L_{\rm dyn}$ for empirical
mean-square penalties for the residual, zero incoming trace, zero opening-plane trace,
and horizon value on the dynamic boundary. The three objectives are
\[
\text{(A)}\quad 10L_r+2L_-,\qquad
\text{(B)}\quad 10L_r+L_-+L_{\rm open},\qquad
\text{(C)}\quad 10L_r+L_-+L_{\rm dyn}.
\]
The recorded results in \Cref{tab:ablation} favor (A). Relative to (A), the mean interior
RMSE is about $15\%$ larger for (B) and a factor $3.1$ larger for (C); the incoming-boundary
RMS is a factor $4.6$ larger for (C). These are comparisons of differently weighted
objectives: (A) doubles the incoming penalty, so the differences do not isolate the added
target. The dynamic target conflicts with the grazing trace by \Cref{R:no-global-ext},
but its causal contribution to these errors requires a matched-weight comparison.
No full-scale retraining is reported.

\begin{table}[htbp!]
\centering\small
\caption{Recorded comparison of the three weighted objectives above (constant-velocity
profile, $v_\mL=1$, $4\times128$ network,
$20{,}000$ steps, two seeds). RMSE against the reference solver on $2{,}500$ states of
$\bm D_{0.05}$, and root-mean-square value on $4{,}000$ points of $\Sigma^-$ (true value $0$),
with $w=1$ and independent uniform draws of $s\in[0,0.99]$ and
$\Lambda\in[1-v_{\max},-10^{-3}]$.}
\label{tab:ablation}
\begin{tabular}{llrr}
\toprule
Variant & Boundary terms & RMSE on $\bm D_{0.05}$ & RMS on $\Sigma^-$\\
\midrule
(A) & $2L_-$ & $0.043$ / $0.044$ & $0.038$ / $0.039$\\
(B) & $L_-+L_{\rm open}$ & $0.051$ / $0.048$ & $0.046$ / $0.043$\\
(C) & $L_-+L_{\rm dyn}$ & $0.133$ / $0.131$ & $0.179$ / $0.177$\\
\bottomrule
\end{tabular}
\end{table}

\paragraph{Empirical and population loss}
Redrawing collocation points does not remove sampling or generalization error. The
training losses, independent conditional-sampler diagnostics, and fixed-signal volume
integrals are different quantities. Only the characteristic-tube comparison uses the
measure in the continuum estimate. Training and evaluation also share the same thirteen
profiles: a one-hot input provides no encoding of unseen forcing. Broader claims would
require a parametric forcing representation and control of the conversion from sampled
to continuum loss.

\subsection{Checks of the regularity theory with a reference solver}\label{S:num-regularity}
The experiments of \Cref{S:num-learning} concern the learned map. The results of
\Cref{sec:ttc-regularity} concern the exact map $\tau$, and several of them make quantitative
predictions that can be checked directly with an accurate solver for the prediction flow, with
no network involved. This subsection reports two such checks: the dependence of the gradient of
$\tau$ on the transversality margin, and the behavior of $\partial_s\tau$ under
piecewise-constant forcing. The samplers restrict nominal leader speed on the
prediction window $[s,1]$ (from the earliest stencil time in the jump experiment),
not on its backward continuation to time zero. They use the prescribed acceleration
without clipping. The fixed-state time trace described below also uses the extended
state domain and is not confined to physical follower speeds.

\paragraph{Reference solver}
The prediction flow \eqref{E:prediction} is integrated together with the variational
equation \eqref{E:var-clean} by an explicit Runge--Kutta method of order eight (DOP853) with
relative and absolute tolerances $10^{-11}$ and $10^{-13}$ and maximal step $0.01$; the
integration is restarted at each jump of $\sa_\mL$. The first contact is located as follows.
Since $\Xi'=\Lambda$, isolated sign-changing zeros of $\Lambda$ locate extrema of the
predicted gap. The procedure detects these events and partitions the dense output at them;
on an interval between consecutive extrema the gap is monotone, and a sign change from
positive to nonpositive is refined by Brent's method. Zeros of $\Lambda$ are not universally
transversal: at an equilibrium, for example, $\Lambda\equiv0$. The stated event search and
finite-sample validation do not certify detection for every tangential or nonisolated zero,
or for arbitrary measurable forcing. The gradient is evaluated by \eqref{eq:grad-tau},
$\nabla_{\bm y}\tau=-e_1^{\top}Y(\tau)/\Lambda(\tau)$. This contact search matters near the
grazing part of $\partial\mathbb D$: a sign-change event on $\Xi$ alone misses predictions that
dip below zero and return within a single integration step. Such missed contacts can occur
near grazing and can cause an incorrect no-contact classification. In a comparison with a
brute-force search on a time grid of step $2\times10^{-6}$, the procedure above misclassified
none of $160$ random states, and its TTC agreed with the brute-force value to the grid step.
Loose and tight tolerance runs agree on all $300$ additional contact classifications.
On twelve transversal states (six pulsed and six piecewise-constant), centered spatial
differences with step $10^{-6}$ agree with \eqref{eq:grad-tau} to a maximal relative error
$8.5\times10^{-11}$. The validation script reproduces these finite checks.
The same solver provides all reference values of \Cref{S:num-learning}.

\paragraph{Gradient and transversality margin}
By \eqref{eq:grad-tau} and the two-sided bound \eqref{E:Y-two-sided} on the state
linearization,
\begin{equation}\label{eq:num-grad-band}
e^{-G_1\tau}\ \le\ \gamma\,\abs{\nabla_{\bm y}\tau}\ \le\ e^{G_1\tau},
\qquad \gamma\Def-\Lambda(\tau;z)>0,\quad G_1=\sqrt{1+2\alpha^2},
\end{equation}
at every state of $\mathbb D$. Thus the gradient of $\tau$ blows up exactly like $1/\gamma$ as
the margin vanishes, and the factor $1/\gamma_{\bm D}$ in the Lipschitz constant of
\Cref{T:tau-Lipschitz} cannot be removed. We test \eqref{eq:num-grad-band} with $\alpha=2$
(so $G_1=3$), $\timehor=1$ and the twelve profiles of \Cref{tab:profiles} with continuous
acceleration, on two samples per
profile: $1{,}500$ accepted states obtained by drawing $s\in[0,0.9]$ uniformly,
leader speed conditionally uniformly on its admissible future-speed interval, follower
speed uniformly in $[0,v_{\max}]$, and gap uniformly in $[10^{-3},1]$, with rejection
outside $\mathbb D$; and $150$ lines of constant
$(s,\Lambda,v_\mL)$ on which the boundary gap $\Xi_b$ of $\mathbb D$ is located by bisection
and the states $\Xi_b-\delta$, $\delta\in[10^{-6},10^{-1}]$ (twelve log-spaced values), are
evaluated, discarding gaps at most $10^{-6}$ and nontransversal contacts. For the lines,
$s$ is uniform on $[0,0.8]$ and speeds follow the same conditional
rule; the seeds are $1000+k$ for profile index $k=0,\ldots,11$. The accepted states are
not jointly volume-uniform, which does not affect this pointwise bound check. This gives $37{,}284$ states with margins from $3.6\times10^{-4}$ to $1.88$ and
gradient norms up to $2.9\times10^{3}$.

\begin{figure}[htbp!]
\centering
\includegraphics[width=\textwidth]{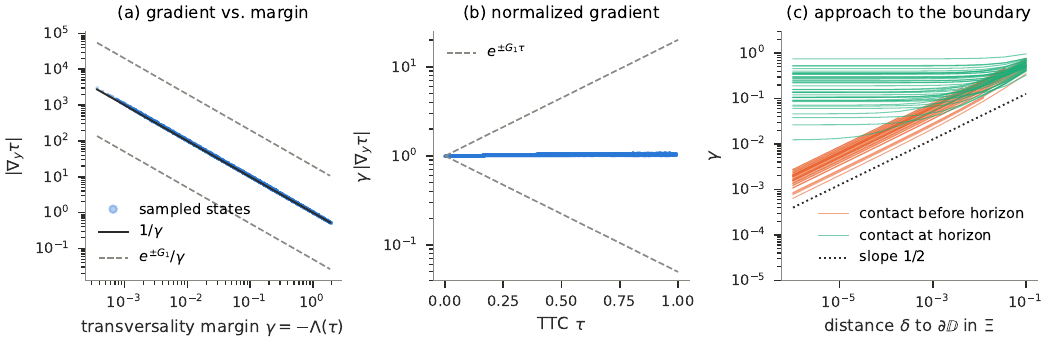}
\caption{Gradient of the TTC and the transversality margin ($\alpha=2$, $\timehor=1$, twelve
profiles, $37{,}284$ states). (a) $\abs{\nabla_{\bm y}\tau}$ against $\gamma$, with $1/\gamma$
and the uniform envelope $e^{\pm G_1}/\gamma$ of \eqref{eq:num-grad-band}. (b) The normalized
gradient $\gamma\abs{\nabla_{\bm y}\tau}$ against $\tau$, with the band $e^{\pm G_1\tau}$.
(c) The margin $\gamma$ at the states $\Xi_b-\delta$ on lines approaching the boundary of
$\mathbb D$ (forty lines of each numerical type shown): grazing-type ($\gamma\sim\delta^{1/2}$
in this sample) and horizon-type (nearly constant $\gamma$), after the numerical
minimum-location check described in the text.}
\label{fig:exp2}
\end{figure}

The computed values satisfy \eqref{eq:num-grad-band} at every sampled state, and a least-squares fit of
$\log\abs{\nabla_{\bm y}\tau}$ against $\log\gamma$ has slope $-1.002$
(\Cref{fig:exp2}(a)). The band is far from tight here: $\gamma\abs{\nabla_{\bm y}\tau}$ lies in
$[1.000,1.049]$ across the sampled states of the twelve profiles, whereas \eqref{eq:num-grad-band} allows
$[e^{-3},e^{3}]$ at $\tau=1$ (\Cref{fig:exp2}(b)). The reason is that $G_1$ bounds the norm of
the whole linearization, while $\gamma\abs{\nabla_{\bm y}\tau}=\abs{e_1^{\top}Y(\tau)}$ involves
only its first row, which starts at $e_1^{\top}$ and, for short horizons, satisfies
$\abs{e_1^{\top}Y(\kappa)}^2=1+(1-\alpha\mathcal V'(\Xi))\kappa^2+O(\kappa^3)$, with $\Xi$ the
initial gap. All sampled states of $\mathbb D$ have $\Xi\le0.78$, hence
$\alpha\mathcal V'(\Xi)\le0.59<1$, and the first row grows slowly. The constant $G_1$
therefore enters \Cref{T:tau-Lipschitz} as a worst case; the dependence on $1/\gamma$ is the
sharp part of the estimate. Since the gradients here are formed using \eqref{eq:grad-tau},
the band comparison is a consistency check on the variational calculation; the finite-difference
comparison reported above provides a separate check of the gradient formula.

The lines toward $\partial\mathbb D$ separate the two parts of the boundary
(\Cref{fig:exp2}(c)). Classification by a terminal-time tolerance of $10^{-3}$ alone labels one near-horizon
grazing line as horizon-type. Checking the location of the minimum gap along each boundary
trajectory corrects that label and gives $1{,}208$ grazing-type and $592$ horizon-type lines.
These are numerical classifications, not a general nondegeneracy certificate. On the
grazing-type sample the observed scaling is $\gamma\sim\delta^{1/2}$ (median fitted exponent $0.502$ for
$\delta\le10^{-3}$, interquartile range $[0.501,0.503]$): near a nondegenerate grazing state
the predicted gap has a quadratic minimum, and moving the initial gap by $\delta$ moves the
contact by $O(\delta^{1/2})$ along a parabola, so that the closing speed at contact is of
order $\delta^{1/2}$. By \eqref{eq:num-grad-band}, $\abs{\nabla_{\bm y}\tau}\sim\delta^{-1/2}$
under this nondegenerate local geometry. This is consistent with the one-sided
square-root behavior and possible jump across grazing described in \cref{R:discontinuity};
it is not a claim about every degenerate grazing point. The horizon-type sample has
median exponent $5\times10^{-4}$, consistent with a nonvanishing contact speed and
local Lipschitz continuation at transversal terminal contacts. This provides a geometric explanation consistent with the
dependence of the network error on the margin in \Cref{tab:profilemetrics}: between the bands
$\gamma\ge0.8$ and $0.05\le\gamma<0.1$ the error grows by a factor of approximately $12$ to $58$, of the
same order as the ratio of the Lipschitz bounds $\sim1/\gamma$ of the function being
approximated. The regularity estimate alone does not determine approximation error; the
boundary objective and sampling also affect these observed ratios.

\paragraph{Piecewise-constant forcing}
\Cref{R:not-C1-in-time} predicts that at a jump time $t_1$ of $\sa_\mL$, the time derivative of
$\tau$ jumps by
\begin{equation}\label{eq:num-jump}
\partial_s\tau(t_1^+,\bm y)-\partial_s\tau(t_1^-,\bm y)
=-\bigl[\sa_\mL(t_1^+)-\sa_\mL(t_1^-)\bigr]\,\partial_\Lambda\tau(t_1,\bm y),
\qquad \partial_\Lambda\tau=\frac{Y_{12}(\tau)}{\abs{\Lambda(\tau)}},
\end{equation}
while $\nabla_{\bm y}\tau$ is continuous across $t_1$, and that $\partial_\Lambda\tau>0$ whenever
$\tau<\pi/\omega_\alpha$. We use the piecewise-constant profile of \Cref{tab:profiles}
($\sa_\mL=0.5,\,-1.2,\,0.8,\,-0.4,\,0$ on the successive intervals of length $0.2$ of $[0,1]$),
so that the jumps are $-1.7$, $2.0$, $-1.2$ and $0.4$ at $t_1=0.2,0.4,0.6,0.8$, with $\alpha=2$
and $\timehor=1$. For each $t_1$ we draw $250$ states in $\mathbb D$ with $\gamma\ge0.05$ and
$\tau\ge0.02$. For each sampled state, the stencil follows one fixed leader trajectory
through $t_1$. We measure the two
one-sided derivatives of $s\mapsto\tau(s,\bm y)$ at $t_1$ by second-order one-sided finite
differences with step $10^{-4}$ (relative and absolute solver tolerances
$10^{-12}$ and $10^{-14}$). Seeds are $\operatorname{int}(1000t_1)$. Leader speed is
sampled for admissibility from $t_1-3\times10^{-4}$ to $1$; exact piecewise integration
confirms the future-speed bounds for every accepted state. As a control we repeat the
test at $t=0.3$, where $\sa_\mL$ is continuous.

\begin{figure}[htbp!]
\centering
\includegraphics[width=\textwidth]{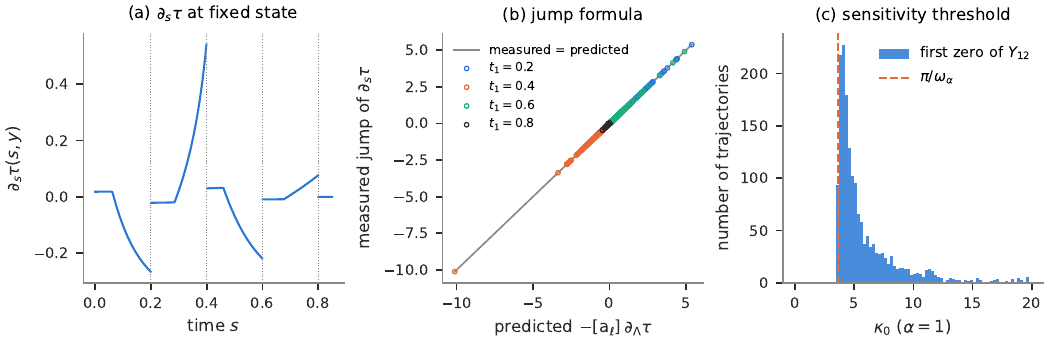}
\caption{Piecewise-constant forcing ($\alpha=2$, $\timehor=1$ in (a), (b)). (a) $\partial_s\tau$
estimated off the forcing jumps at fixed $(\Xi,\Lambda)=(0.1,-1)$ and initial
leader speed $0.9$; dotted lines mark the jumps of $\sa_\mL$. This extended-domain
trace briefly has follower speed above $v_{\max}$ (maximum $2$). The kinks between jumps occur where the contact time
$s+\tau$ crosses a jump time. (b) Measured jump of $\partial_s\tau$ against the prediction
\eqref{eq:num-jump}, $250$ states for each jump time. (c) First zero $\kappa_0$ of $Y_{12}$ on
$2{,}000$ prediction trajectories with $\alpha=1$ and horizon $20$ ($1{,}728$ have a zero), with
the bound $\pi/\omega_\alpha$ of \Cref{R:not-C1-in-time}. Panel (c) concerns the
global variational flow, including post-contact continuation; its zero times are not
necessarily first-contact times.}
\label{fig:exp3}
\end{figure}

The measured jumps range in absolute value from $0.007$ to $10.1$ and agree with
\eqref{eq:num-jump} with a median relative error of $1.6\times10^{-7}$ and a maximal relative
error of $3.1\times10^{-4}$ over the $1{,}000$ states (\Cref{fig:exp3}(b), \Cref{tab:jumps}); the
deviations are finite-step effects. Some stencils cross a later forcing change at
the contact time, so second-order accuracy is not uniform across those stencils;
targeted step refinement substantially reduces their errors. At the control time the
measured jump is at most $2.4\times10^{-5}$ in absolute value. On each side of every jump the classical HJB equation
$\partial_s\tau+\mathsf G(t_1^\pm,\bm y)\cdot\nabla_{\bm y}\tau+1=0$ holds with the one-sided
derivative. The median, over the $1{,}000$ jump states, of the larger absolute
residual on the two sides is $2.4\times10^{-8}$; the largest residual is approximately
$1.5\times10^{-3}$,
where the one-sided derivative is largest, and is again at the level of the finite-difference
error. The sensitivity $\partial_\Lambda\tau$ is positive at every tested state (minimum
$0.012$), as predicted, since $\tau<1<\pi/\omega_2=\pi$. \Cref{fig:exp3}(a) shows the
corresponding extended-domain illustration at one fixed $\bm y$: its implied
follower speed reaches $2>v_{\max}$, so it is not an entirely physical-state trace.
The formula predicts jumps of $\partial_s\tau$ at each $t_1$. Between these times, the trace
also shows continuous kinks where the contact time $s+\tau$ crosses a later forcing change.
Thus $\tau$ is Lipschitz but not $C^1$ in time, and the HJB equation holds classically on each
side of a jump but not at it.

\begin{table}[htbp!]
\centering\small
\caption{Jump of $\partial_s\tau$ at the jump times of a piecewise-constant leader
acceleration, $250$ states per time ($\alpha=2$, $\timehor=1$); ``rel.\ error'' compares the
measured jump with \eqref{eq:num-jump}, and the control time $0.3$ has no jump of $\sa_\mL$.}
\label{tab:jumps}
\begin{tabular}{rrrrr}
\toprule
$t_1$ & $[\sa_\mL](t_1)$ & range of $\abs{\text{measured jump}}$ & max.\ rel.\ error & max.\ HJB residual\\
\midrule
$0.2$ & $-1.7$ & $0.029$--$5.37$ & $3.1\times10^{-4}$ & $5.3\times10^{-4}$\\
$0.4$ & $2.0$ & $0.024$--$10.1$ & $1.6\times10^{-4}$ & $1.5\times10^{-3}$\\
$0.6$ & $-1.2$ & $0.019$--$4.90$ & $4.9\times10^{-5}$ & $2.0\times10^{-4}$\\
$0.8$ & $0.4$ & $0.007$--$0.461$ & $5.2\times10^{-6}$ & $9.0\times10^{-7}$\\
$0.3$ & $0$ & $\le2.4\times10^{-5}$ & -- & $3.9\times10^{-4}$\\
\bottomrule
\end{tabular}
\end{table}

\paragraph{The variational-flow sensitivity threshold}
The sign of $\partial_\Lambda\tau$ is the sign of $Y_{12}(\tau)$, and the proof in
\Cref{R:not-C1-in-time} shows that $Y_{12}(\kappa)>0$ for $0<\kappa<\pi/\omega_\alpha$ along
every prediction trajectory, and for all $\kappa>0$ if $\alpha\ge4$. The bound
$\pi/\omega_\alpha$ comes from the comparison $\mathcal V'\le1$ and is attained only by a
trajectory that stays where $\mathcal V'=1$, that is, at $\Xi=2$. With the unit horizon used
above the threshold is never reached, so we test the statement itself, on long prediction
trajectories with a constant leader speed and horizon $20$ (the prediction flow is global, so
trajectories are followed through $\Xi=0$), from $2{,}000$ initial states with
$\Xi\in[0,4]$ and both speeds in $[0,v_{\max}]$. This separate DOP853 integration uses
relative tolerance $10^{-10}$, absolute tolerance $10^{-12}$, and maximum step $0.05$;
the seeds are $1$ and $2$ for $\alpha=1$ and $4.5$, respectively. For $\alpha=1$, where
$\pi/\omega_1=2\pi/\sqrt3\approx3.6276$, $Y_{12}$ vanishes on $1{,}728$ of the trajectories,
never before $\pi/\omega_1$; the earliest zero, at $\kappa_0=3.6300$, is $0.07\%$ above the
bound and occurs for an initial gap $\Xi=2.14$, near the point where the bound is attained
(\Cref{fig:exp3}(c)). For $\alpha=4.5$, $Y_{12}$ has no zero on any of the $2{,}000$
trajectories. These observations support near-sharpness of the positivity bound for the global
variational flow. They do not establish optimal sharpness for TTC sensitivity: the recorded
zero need not occur at first contact, and the equilibrium $\Xi=2$ at which
$\mathcal V'=1$ never reaches the collision set. The first-contact experiment below addresses
the separate question of whether TTC sensitivity can vanish or become negative after the
sufficient threshold.

\paragraph{Sensitivity at actual first contact}
A separate reproducible calculation evaluates the variational solution at the first
transversal contact itself. Take $\alpha=0.2$, $\timehor=20$, initial gap $\Xi_0=2$, and
$v_\mL(0)=v(0)=\tanh2$. Hold the leader speed until time $h$, brake linearly to zero over
$0.1$ time units, and then keep the leader stopped. This defines a physical prescribed
leader speed with bounded piecewise-constant acceleration. The prediction remains at the
equilibrium gap until braking begins; after braking its gap decreases to a first transversal
contact. Integrating the prediction together with the column $(Y_{12},Y_{22})$ of its state
linearization gives \Cref{tab:first-contact-sensitivity}. The calculation uses DOP853,
relative and absolute tolerances $2\times10^{-12}$ and $2\times10^{-14}$, maximum step $0.01$,
and restarts at both acceleration changes. A scalar Brent solve in $h$ locates the near-zero
row. The script and generated records are available in the companion repository under
\texttt{reproducibility/actual\_first\_contact\_sensitivity.py}.

\begin{table}[htbp!]
\centering\small
\caption{Variational sensitivity evaluated at actual first contact for the hold--brake--stop
leader signal. Here $\alpha=0.2$, $\Xi_0=2$, initial speeds equal $\tanh2$, and
$\pi/\omega_{0.2}\approx7.2073078415$. The middle row is a numerical root, not an exact
symbolic zero.}
\label{tab:first-contact-sensitivity}
\begin{tabular}{rrrr}
\toprule
Hold $h$ & First-contact $\tau$ & $Y_{12}(\tau)$ & $\partial_\Lambda\tau$\\
\midrule
$0$ & $2.3959704553$ & $1.7481387877$ & $2.5189797785$\\
$4.9932259521$ & $7.3891964074$ & $\approx0$ & $\approx0$\\
$8$ & $10.3959704553$ & $-0.9233070010$ & $-1.3304387965$\\
\bottomrule
\end{tabular}
\end{table}

The contact relative velocity is approximately $-0.6939868286$ in all three cases. Thus the
computed zero and negative sensitivities occur at transversal first contacts after the
sufficient positivity threshold. This example demonstrates the phenomenon at the relevant
endpoint, but does not establish an optimal threshold. These are trajectories of the
prediction flow; they are not collisions of the closed-loop system.

\subsection{Scope of the experiments}\label{S:num-scope}

The checks of \Cref{S:num-regularity} are consistent with quantitative predictions about
the exact map $\tau$ on the sampled states, subject to the stated solver tolerances. The learning experiments of
\Cref{S:num-learning} report pooled empirical accuracy and estimate the continuum losses
on one compact tube for each of three fixed signals, but a finite set of scenarios does not prove the theorems
above, and the reported behavior is consistent with, not evidence for, the analytical
guarantees. Several limitations should be stated. First, \Cref{thm:two-sided-convergence} bounds
the $L^2$ error of the learned value, and the impact-prevention theorem concerns the exact
coefficient $\tau$; neither gives a near-contact bound on the reciprocal gain
$1/\bm V_\theta^2$, and a uniformly small absolute TTC error does not by itself preserve
impact-freeness: with a stationary leader and $\bm V_\theta=\tau+\varepsilon$ the feedback gain
is bounded by $\alpha+\beta/\varepsilon^2$ near contact, and a follower started close enough to
the leader reaches it with negative relative velocity. The nonzero boundary traces measured in \Cref{S:num-learning} illustrate why a small
absolute TTC error does not control this reciprocal gain; their RMS value is not itself a
closed-loop collision test. A safety statement for a controller driven by the learned map would need a
relative near-contact estimate (of the type $0<\bm V_\theta\le CX/\abs V$ on the approach
region) or a separate safety filter; neither is established here. Second, training and
evaluation use the same thirteen profiles, encoded as a one-hot input, so the results speak to
accuracy on these profiles only; a statement about unseen forcing would require a parametric
encoding of the profile. Third, Theorems~\ref{T:two-sided} and~\ref{T:flux-stability} bound the error by continuum
losses on a fixed-signal domain, whereas the training and fresh residual diagnostics use other sampling laws.
A bound connecting these quantities is not established. Fourth, the principal checkpoint
uses an additional dynamic-boundary target that conflicts with the grazing trace. The
smaller comparison changes the incoming-boundary weight as well, so it does not isolate
the effect of that target. None of these checks is a learned-feedback safety certificate.

\paragraph{Reproducibility}
The companion repository provides the training and continuation programs, the principal
checkpoint, the original validation and test sets, and the evaluation and seeded regularity
scripts with reference arrays and software versions. The two $4{,}096$-state held-out sets
regenerate from their recorded seeds, and checkpoint evaluation recovers the logged metrics
to the documented numerical tolerances. These labels use the speed-clipped training reference;
the independent evaluations use the prescribed flow without clipping. The ablation code and
summaries are included, but the six original ablation checkpoints and raw evaluation arrays
are not archived; exact replay of the long training history also requires the historical
step-$500$ starting state.

\section{Conclusions}\label{S:conclusion}
We have analyzed a TTC feedback whose collision-time coefficient is generated by a nominal OVM prediction with bounded measurable leader acceleration. The prediction excludes the additional feedback, so its collision time describes a counterfactual nominal trajectory. Its grazing geometry can make the capped TTC, and hence the feedback field, discontinuous. For the prescribed leader--follower problem, Filippov solutions exist for physical initial data and every contact has zero closing speed. Under the stated nondegeneracy hypotheses, the feedback moves trajectories transversally through the grazing set and uniqueness persists. These conclusions do not exclude a zero-speed touch or establish uniqueness at degenerate grazing states.

On the transversal first-contact domain, the TTC is locally Lipschitz, has a jointly continuous state gradient, and satisfies its Hamilton--Jacobi equation at every state for almost every time, with a common exceptional time set. The quantitative estimates explain the dependence of pointwise sensitivity on the terminal closing speed. A jump of the prescribed acceleration produces the stated time-derivative jump whenever the relative-velocity sensitivity is nonzero. The sufficient positivity horizon follows from a variational comparison; optimal sharpness for actual first-contact sensitivity is not established.

The approximation analysis controls the continuum $L^2$ error by the incoming boundary loss and the interior residual. It applies to locally Lipschitz candidates with a continuous incoming trace, including smooth augmented networks evaluated at the evolving leader speed under rough forcing. Besides the coarea estimate, an endpoint/age change of variables gives the explicit coefficient in \cref{T:flux-stability}. Its normal-flux factor explains why integrated residual stability remains uniform within bounded families even when pointwise gradients become large near grazing. Zero total loss yields the stated pointwise identification, and a sequence whose continuum loss vanishes converges in $L^2$. These results are conditional on the loss, without a theorem about convergence of the optimizer or a network-size approximation rate.

The numerical study distinguishes boundary objectives and sampling measures. The principal checkpoint includes a dynamic-boundary target that conflicts with the hitting-side trace on grazing portions; the smaller objective comparison also changes the incoming-boundary weight and therefore does not isolate the extra penalty. The characteristic-tube experiment gives empirical error-to-loss ratios near $0.10$--$0.14$, compared with the explicit analytical squared-error coefficient $0.3638072604$. The losses in this comparison are estimated integrals. The observed bias under piecewise-constant forcing is not predicted by the smooth-class nonattainment result. Reference calculations support the gradient and forcing-jump formulas, while the separate reproducible first-contact example demonstrates post-threshold sensitivity sign changes without proving an optimal threshold.

Further work includes full-scale training with the correct incoming-boundary objective, quantitative control of sampling and optimization errors, and analysis of a closed loop driven by an approximate TTC. A uniform absolute error in TTC alone does not control the reciprocal feedback gain near contact. Safety of a learned implementation therefore requires an appropriate near-contact estimate or an additional safety mechanism. Extensions to controlled platoons, uncertain leader forecasts, delay, and more general prediction models also require hypotheses beyond those proved for the present prescribed leader--follower system.

\section*{Code availability}
Code, trained checkpoints, reference data, and reproduction instructions are provided at
\url{https://github.com/hosein-matin/TTC-HJB-Feedback-Regularity}.
The scope of reproducibility is stated at the end of Section~\ref{sec:numerics}.

\bibliographystyle{siamplain}
\bibliography{bibliography}

\end{document}